\documentclass[10pt,french]{article}

\usepackage{latexsym}
\usepackage[all]{xy}
\usepackage{amsmath}
\usepackage{amsfonts}

\newtheorem{thm}{Th\'eor\`eme}[subsection]
\newtheorem{prop}[thm]{Proposition}
\newtheorem{lem}[thm]{Lemme}
\newtheorem{df}[thm]{D\'efinition}
\newtheorem{cor}[thm]{Corollaire}
\newtheorem{conj}[thm]{Conjecture}

\begin{document}

\title{Champs affines}

\author{Bertrand To\"en \\
Laboratoire Emile Picard UMR CNRS 5580 \\
Universit\'e Paul Sabatier, 118, route de Narbonne\\
31062 Toulouse Cedex 04, France\\
email: toen@math.ups-tlse.fr}
\date{Avril 2006}

\maketitle

\textbf{Keywords}: Stacks, schematic homotopy type, 
rational homotopy theory, p-adic homotopy theory. \\

\textbf{MSC}:14F25 \\

\begin{abstract}
The purpose of this work is to introduce a notion of
\emph{affine stacks}, which is a homotopy version of the notion
of affine schemes, and to give several applications in the context of
algebraic topology and algebraic geometry.

As a first application we show how affine stacks can be used
in order to give a new point of view (and new proofs) on
rational and $p$-adic homotopy theory. This gives a first solution to
A. Grothendieck's \emph{schematization problem} described in \cite{gr}.

We also use affine stacks in order to introduce a notion
of \emph{schematic homotopy types}. We show that schematic homotopy types give
a second solution to the schematization problem, which also allows us to
go beyound rational and $p$-adic homotopy theory for spaces with
arbitrary fundamental groups. The notion of schematic homotopy types is
also used in order to construct various homotopy types of algebraic varieties
corresponding to various cohomology theories (Betti, de Rham, $l$-adic, \dots),
extending the well known constructions of the various fundamental groups.

Finally, as algebraic stacks are obtained by gluing affine schemes we
define \emph{$\infty$-geometric stacks} as a certain gluing
of affine stacks. Example of $\infty$-geometric stacks in the context
of algebraic topology (moduli spaces of dga structures up to quasi-isomorphisms) and Hodge theory
(non-abelian periods) are given.
\end{abstract}

\newpage

\tableofcontents

\newpage

\setcounter{section}{-1}

\section{Introduction}

Le but principal de ce travail est de pr\'esenter une notion
de \textit{champ\footnote{Dans ce texte l'expression \emph{champ}
signiefira toujours \emph{champ en $\infty$-groupoides}.
Un mod\`ele pour la th\'eorie des champs (que nous utiliserons) est celui des
pr\'efaisceaux simpliciaux de A. Joyal et R. Jardine (voir \cite{jo,j}).} affine},
g\'en\'eralisation de nature homotopique
de la notion de sch\'ema affine, et d'en donner plusieurs applications
que nous allons commencer par r\'esumer.

\begin{itemize}

\item Dans un premier temps on montrera comment la notion de champs affines peut-\^etre utilis\'ee
afin de r\'einterpr\'eter g\'eom\'etriquement, et d'\'eclairer, les th\'eories de l'homotopie
rationnelle et $p$-adique (par exemple au sens de \cite{q,bo,bg,s,k,g,m1}).
Ceci passe par l'existence d'un \emph{foncteur d'affination}, associant
\`a tout type d'homotopie un champ affine v\'erifiant une propri\'et\'e universelle.
Nous proposerons cette construction comme une premi\`ere solution au
\emph{probl\`eme de la sch\'ematisation} de A. Grothendieck de \cite{gr} (voir
aussi l'appendice $A$).
\item Les champs affines seront utilis\'es pour d\'efinir
une notion de \textit{type d'homotopie sch\'ematique} que
nous proposons comme mod\`eles pour construire des types d'homotopie de vari\'et\'es
alg\'ebriques. Nous d\'emontrerons l'existence d'un \textit{foncteur de sch\'ematisation},
associant \`a un type d'homotopie (au sens usuel)
un type d'homotopie sch\'ematique v\'erifiant une propri\'et\'e universelle.
Cette construction donne une seconde solution au
probl\`eme de la sch\'ematisation, et permet de plus d'aller
au-del\`a des th\'eories de l'homotopie
rationnelle et $p$-adique pour des espaces \`a groupes
fondamentaux arbitraires.
\item On construit, pour les th\'eories cohomologiques usuelles de la g\'eom\'etrie
alg\'e\-brique (Betti, de Rham, Hodge,
$l$-adique et cristalline) des th\'eories homotopiques associ\'ees dont les valeurs
sont des types d'homotopie sch\'ematiques. On esp\`ere que ces th\'eories
d'homotopie capturent des informations arithm\'e\-tiques et/ou g\'eom\'e\-triques
int\'eressantes. Dans le cas des vari\'et\'es complexes projectives et de
la th\'eorie de Betti, on explique comment
le type d'homotopie associ\'e poss\`ede une \textit{d\'ecomposition
de Hodge} g\'en\'eralisant les structures de Hodge sur la coholomogie, le groupe
fondamental et sur le type d'homotopie rationnel. De m\^eme, dans le cas
de la cohomologie cristalline le type d'ho\-mo\-to\-pie associ\'e poss\`ede
une structure de $F$-isocristal (due a M. Olsson, voir \cite{ol}), et pour le cas
de la cohomologie $l$-adique une action continue du groupe de Galois du corps
de base. On propose aussi une construction d'une variante
non-ab\'elienne des applications d'Abel Jacobi.
\item Tout comme un sch\'ema (ou plus g\'en\'eralement un champ
alg\'ebrique) est obtenu par recollement de sch\'emas affines, on peut
recoller des champs affines pour obtenir une notion
de \textit{champs $\infty$-g\'eom\'etriques}. Cette notion, qui est une g\'en\'eralisation
de la notion de champs alg\'ebriques, permet de r\'esoudre de nouveaux
\textit{probl\`emes de modules} qui ne peuvent
\^etre raisonnablement r\'esolus \`a l'aide de la notion usuelle de champs alg\'ebriques.
On construit par exemple le champ $\infty$-g\'eo\-m\'e\-trique
qui classifie les structures
d'alg\`ebre diff\'erentielle gradu\'ee \`a quasi-isomorphismes pr\`es.
On construit aussi le champ $\infty$-g\'eo\-m\'e\-trique classifiant les filtrations
sur une alg\`ebre diff\'erentielle gradu\'ee commutative, que l'on utilise
pour d\'efinir une \emph{application des p\'eriodes non-ab\'eliennes}
qui contr\^ole la variation de la filtration de Hodge sur le type
d'homotopie rationnel d'une famille de vari\'et\'es complexes lisses et projectives.

\end{itemize}

Certains des points pr\'ec\'edents seront trait\'es en d\'etails, mais d'autres
feront l'objet d'articles \`a parts enti\`eres et ne seront que survol\'es
dans ce travail. \\

D\'ecrivons \`a pr\'esent le contenu math\'ematique de cet article. \\

\begin{center} \textit{Champs} \end{center}

Un mod\`ele ad\'equat pour une th\'eorie des \textit{champs en $\infty$-groupoides}
est la th\'eorie des pr\'efaisceaux simpliciaux
(voir \cite{jo,j,j2,bla,du}). Ainsi, le mot \textit{champ}
signifiera pour nous un objet de la cat\'egorie homotopique des pr\'efaisceaux simpliciaux. Pour un anneau $k$, nous
consid\'ererons la cat\'egorie des sch\'emas affines sur $Spec\, k$, munie de la topologie fid\`element plate
et quasi-compacte. La cat\'egorie homotopique des pr\'efaisceaux simpliciaux sur ce site sera notre cat\'egorie
des champs sur $Spec\, k$.

D'apr\`es les travaux fondamentaux \cite{jo,j} la th\'eorie homotopique des champs
admet une structure de mod\`eles. Une cons\'equence importante, et non triviale, de l'existence de
cette structure de mod\`eles est l'existence de constructions standards tel que
les limites et colimites homotopiques. De fa\c{c}on plus pr\'ecise la cat\'egorie
de mod\`eles des pr\'efaisceaux simpliciaux est un \textit{topos de mod\`eles}, au
sens de \cite{tv}. Ceci implique en particulier des propri\'et\'es d'exactitudes
additionelles de la th\'eorie des champs, comme par exemple l'existence
de Hom internes. Il est bon de garder \`a l'esprit que la th\'eorie des champs
fonctionne de fa\c{c}on tout \`a fait similaire \`a celle des faisceaux. Une pr\'esentation
de quelques constructions et propri\'et\'es standards de la th\'eorie des champs
fait l'objet du chapitre \S 1. \\

\begin{center} \textit{Champs affines} \end{center}

La notion de \textit{champs affines} sur un anneau $k$ est une
\textit{version d\'eriv\'ee de la notion de sch\'emas affines}, et g\'en\'eralise celle-ci dans
un cadre homotopique.
Pour d\'efinir cette notion nous utiliserons
une structure de cat\'egorie de mod\`eles sur la cat\'egorie des $k$-alg\`ebres co-simpliciales. Pour $A$ une telle
alg\`ebre co-simpliciale, nous consid\'ererons le sch\'ema affine simplicial $Spec\, A$, ainsi
que le pr\'efaisceau simplicial sur le site des $k$-sch\'emas affines qu'il repr\'esente. Nous d\'efinissons ainsi un
foncteur
$$Spec  : (k-Alg^{\Delta})^{op} \longrightarrow SPr(k),$$
de la cat\'egorie des $k$-alg\`ebres co-simpliciales,
vers la cat\'egorie des pr\'efaisceaux simpliciaux sur le site des $k$-sch\'emas affines pour la topologie
fid\`element plate et quasi-compacte. Notre premi\`ere observation, qui est au centre du pr\'esent travail est
la proposition suivante.

\begin{prop}{(Voir Prop.  \ref{p2} et Cor. \ref{c1})}\label{pi}
Le foncteur $Spec$ est un foncteur de Quillen \`a droite. De plus, le foncteur d\'eriv\'e
$$\mathbb{R}Spec : Ho(k-Alg^{\Delta})^{op} \longrightarrow Ho(SPr(k))$$
est pleinement fid\`ele.
\end{prop}

Avant d'aller plus loin dans la d\'efinition des champs affines il me semble important de faire
deux remarques.

\begin{itemize}

\item
Lorsque la cat\'egorie des sch\'emas affines sur $k$ est munie de la topologie
triviale, la proposition pr\'ec\'edente est plus ou moins une trivialit\'e. En contre partie, lorsque l'on
travaille avec la topologie fid\`element plate, la preuve du fait que $Spec$ est un foncteur de Quillen \`a droite
utilise un r\'esultat r\'ecent de D. Dugger, S. Hollander et D. Isaksen
qui caract\'erise les pr\'efaisceaux simpliciaux
fibrants (voir \cite{du}). Il est important de noter ici que l'utilisation de la topologie fid\`element plate
est cruciale pour les r\'esultats qui vont suivre, qui sont de toute \'evidence faux pour la topologie triviale.
Par exemple, une des propri\'et\'es fondamentales de cette topologie que nous utiliserons est le fait que la sous-cat\'egorie
pleine de la cat\'egorie des faisceaux en groupes, form\'ee des faisceaux repr\'esentables par des sch\'emas en groupes
affines, est stable par conoyaux et limites projectives (du moins lorsque $k$ est un corps). Cette propri\'et\'e
sera tr\`es importante pour les th\'eor\`emes fondamentaux \ref{ti} et \ref{ti2}.

\item

La proposition \ref{pi} implique donc que la th\'eorie homotopique des $k$-alg\`ebres co-simpliciales se plongent dans celle
des champs sur $Spec\, k$ pour la topologie fid\`element plate. Le foncteur $\mathbb{R}Spec$ permet ainsi de \textit{g\'eom\'etriser}
les objets alg\'ebriques que sont les $k$-alg\`ebres co-simpliciales. De plus, lorsque $k$ est de caract\'eristique
nulle, les th\'eories homotopiques des $k$-alg\`ebres co-simpliciales et des $k$-alg\`ebres diff\'erentielles gradu\'ees
commutatives et en degr\'es positifs sont \'equivalentes (voir \cite{hisc}). Ainsi, nous obtenons aussi un plongement de la cat\'egorie homtopique des
$k$-alg\`ebres diff\'erentielles gradu\'ees commutatives (et en degr\'es positifs)
dans la cat\'egorie des champs sur $Spec\, k$.

\end{itemize}

Nous d\'efinissons la cat\'egorie des champs affines sur $k$ commme la sous-cat\'egorie pleine des champs sur
$Spec\, k$ image essentielle du foncteur $\mathbb{R}Spec$. Ainsi, en caract\'eristique nulle, la cat\'egorie
des champs affines sur $k$ est anti-\'equivalente \`a la cat\'egorie homotopique des $k$-alg\`ebres diff\'erentielles
gradu\'ees commutatives et en degr\'es positifs.
Il se trouve que la cat\'egorie des champs affines ainsi d\'efinie poss\`ede plusieurs descriptions
\'equivalentes. Nous commencerons par montrer que les champs affines sont les
\textit{objets locaux pour la th\'eorie cohomologique
repr\'esent\'ee par $\mathbb{G}_{a}$, et de taille raisonnable} (voir
Thm. \ref{t4}). On diposera ainsi d'une interpr\'etation
de nature cohomologique de la notion de champs affines. Nous montrerons
aussi que la cat\'egorie des champs affines
est aussi la plus petite
sous-cat\'egorie pleine de la cat\'egorie des champs sur $Spec\, k$ contenant les champs de la forme 
$K(\mathbb{G}_{a},n)$ et qui
est stable par limites homotopiques (voir la remarque qui suit le th\'eor\`eme \ref{t4}). Ceci montre en particulier
que la notion de champs affines est intrins\`eque \`a la th\'eorie des champs sur $Spec\, k$.
Enfin, lorsque $k$ est un corps, la notion de champs affines se simplifie quelque peu, et nous donnerons alors
la caract\'erisation suivante des champs affines, point\'es et connexes.

\begin{thm}{(Voir Thm. \ref{t5} et Thm. \ref{t5'})}\label{ti}
Les champs affines, point\'es et connexes sur un corps $k$ sont exactement les pr\'efaisceaux simpliciaux
point\'es et connexes $F$, tels que pour tout $i>0$ le faisceau $\pi_{i}(F,*)$ soit repr\'esentable par un
sch\'ema en groupes affine et unipotent.
\end{thm}

En corollaire de ce th\'eor\`eme on obtient qu'il existe un \'equivalence entre la th\'eorie
homotopique des $k$-alg\`ebres co-simpliciales (ou encore
$k$-alg\`ebres di\-ff\'e\-ren\-ti\-elles gradu\'ees si $k$ est de caract\'eristique nulle)
augment\'ees et connexes et celle des champs point\'es v\'erifiant les conditions
du th\'eor\`eme \ref{ti}. Cette intepr\'etation g\'eom\'etrique des $k$-alg\`ebres
co-simpliciales est \`a ma connaissance un r\'esultat nouveau, tout particuli\`erement
dans le cas o\`u $k$ est de caract\'eristique positive. C'est aussi un th\'eor\`eme non-trivial car il permet
par exemple de retrouver (et ce de fa\c{c}on purement formelle)
les th\'eor\`emes fondamentaux de l'homotopie rationnelle et $p$-adique (voir
Thm. \ref{t8} et Cor. \ref{ct8}).  \\

\begin{center} \textit{Affination des types d'homotopie} \end{center}

Les champs affines forment une cat\'egorie de champs qui donne une solution \textit{probl\`eme de la
sch\'ematisation de A. Grothendieck} tel que je l'ai personnellement compris (voir
l'appendice A). Dans ce cadre il s'interpr\'etera de la fa\c{c}on suivante.

Si $k$ est un anneau, on peut consid\'erer le foncteur d\'eriv\'e des sections globales
$\mathbb{R}\Gamma$, de la cat\'egorie des champs affines sur $k$ vers la cat\'egorie homotopique des ensembles simpliciaux.
Notre premier r\'esultat d'existence est le suivant, et est une cons\'equence imm\'ediate de la proposition \ref{pi}.

\begin{cor}{(voir Cor. \ref{caff})}
Le foncteur $\mathbb{R}\Gamma$ restreint aux champs affines poss\`ede un adjoint \`a gauche $X \mapsto (X\otimes k)^{uni}$.
De plus, si $k^{X}$ est la $k$-alg\`ebre co-simpliciale de cohomologie de $X$ \`a valeurs dans $k$ on a
$$(X\otimes k)^{uni} \simeq \mathbb{R}Spec\, (k^{X}).$$
Le champ affine $(X\otimes k)^{uni}$ est appel\'e une \textit{affination} de $X$ sur $k$.
\end{cor}

Il est important de remarquer que le fait que les champs $K(\mathbb{G}_{a},n)$ soient des champs affines implique
qu'il existe une propri\'et\'e de conservation de la cohomologie
$H^{*}(X,k)\simeq H^{*}((X\otimes k)^{uni},\mathbb{G}_{a})$
(qui est une des propri\'et\'es ch\`eres \`a A. Grothendieck dans \cite{gr}).
Ainsi, comme les champs affines sont des objets locaux pour la th\'eorie cohomologique repr\'esent\'es par
$\mathbb{G}_{a}$, le morphisme d'adjonction $X \longrightarrow (X\otimes k)^{uni}$ peut
aussi s'interpr\'eter comme un morphisme de \textit{localisation} (dans le sens que A.K. Bousfield a donn\'e \`a ce terme
dans \cite{bo}).

En utilisant le th\'eor\`eme \ref{ti}, il
est possible de d\'ecrire le champ $(X\otimes k)^{uni}$, du moins lorsque $X$ est nilpotent de type
fini et $k$ est un corps. Dans ce cas, nous montrerons que $(X\otimes \mathbb{Q})^{uni}$ est un mod\`ele
pour le type d'homotopie
rationnel de $X$, et
$(X\otimes \overline{\mathbb{F}}_{p})^{uni}$ un mod\`ele pour son type d'homotopie $p$-adique (voir Thm. \ref{t8} et Cor.
\ref{ct8}). Ces r\'esultats sont des cons\'equences formelle de la propri\'et\'e universelle
satisfaite par $(X\otimes k)^{uni}$ et du th\'eor\`eme \ref{ti}.
La notion de champ affine permet ainsi de
construire des mod\`eles \textit{g\'eom\'etrico-alg\'ebriques} aux types d'homotopie.

Enfin, signalons que lorsque $X$ est un ensemble simplicial fini,
le foncteur $X \mapsto (X\otimes k)^{uni}$ est compatible aux changements d'anneaux. De ce fait,
le champ affine $(X\otimes \mathbb{Z})^{uni}$ permet d'unifier les types d'homotopie rationnel et $p$-adique
de $X$. Nous donnerons une description conjecturale de ce champ lorsque $X$ est $1$-connexe et de type fini au
paragraphe \S 2.3 (voir Conj. \ref{conj2}). \\

\begin{center} \textit{$\infty$-Gerbes affines et types d'homotopie sch\'ematiques} \end{center}

De par leur d\'efinition, les types d'homotopie repr\'esent\'es par des champs affines ont une forte tendence
\`a \^etre nilpotents. Ceci est un fait tr\`es bien connu en th\'eorie de l'homotopie rationnelle ou $p$-adique, o\`u les
groupes fondamentaux que l'on obtient \`a partir des mod\`eles alg\'ebriques sont les compl\'etions de Mal'cev
(relativement au corps de base) des groupes fondamentaux usuels. De m\^eme, pour un espace $X$
qui est non-nilpotent, le champ $(X\otimes k)^{uni}$ ne verra pas la partie r\'eductive du type d'homotopie de
$X$ (e.g. le sch\'ema en groupes $\pi_{1}((X\otimes k)^{uni},x)$ est un sch\'ema en groupes affine pro-unipotent).
On est donc naturellement amen\'es
\`a introduire une notion plus g\'en\'erale que celle de champs affines qui soit moins restrictive au niveau des groupes
fondamentaux. C'est pour cette notion plus g\'en\'erale que nous r\'eserverons l'expression \textit{types d'homotopie
sch\'ematiques}.

Tout d'abord nous introduirons une notion d'\textit{$\infty$-gerbe affine}
qui g\'en\'eralise la notion de gerbe affine utilis\'ee dans le formalisme Tannakien (voir \cite{sa,d2}), et de
fa\c{c}on \`a ce que les \textit{$\infty$-gerbes affines
soient aux gerbes affines ce que les champs affines sont aux sch\'emas affines}. En clair, nous d\'efinirons une
$\infty$-gerbe affine (point\'ee) sur un corps $k$ comme \'etant un champ point\'e et connexe $F$ sur $Spec\, k$, dont le champ
des lacets $\Omega_{*}F$ est un champ affine sur $k$ (voir Def. \ref{d10}). Les types d'homotopie sch\'ematiques point\'es
sur  $k$ seront alors les $\infty$-gerbes affines point\'ees poss\'edant une certaine propri\'et\'e de
localit\'e cohomologique (voir Def. \ref{d10}). Nous conjecturons cependant que cette propri\'et\'e est toujours
satisfaite, et que toute $\infty$-gerbe affine point\'ee est donc un type d'homotopie sch\'ematique (voir \ref{conj}).
Quoiqu'il en soit, la cat\'egorie des types d'homotopie sch\'ematiques point\'es ainsi d\'efinie est
une sous-cat\'egorie pleine de la cat\'egorie homotopique des pr\'efaisceaux simpliciaux point\'es
sur le site des sch\'emas affines sur le corps $k$.
Les $\infty$-gerbes affines et les types d'homotopie sch\'ematiques
sur un anneau de base plus g\'en\'eral ne seront pas d\'efinis dans ce travail. \\

Le premier r\'esultat que nous d\'emontrerons est le crit\`ere de
reconnaissance suivant, qui permet de donner de nombreux exemples
int\'eressants de types d'ho\-mo\-to\-pie sch\'ematiques point\'es.

\begin{thm}{(Voir Thm. \ref{t6} et Thm. \ref{t6'})}\label{ti2}
Les types d'homotopie sch\'e\-ma\-ti\-ques sont exactement les pr\'efaisceaux simpliciaux point\'es et connexes
tels que pour tout $i>0$, le faisceau $\pi_{i}(F,*)$
soit repr\'esentable par un sch\'ema en groupes affine, qui est de plus unipotent pour $i>1$.

De plus, si $F$ est une $\infty$-gerbe affine point\'ee, alors pour tout $i>1$ le faisceau
$\pi_{i}(F,*)$ est repr\'esentable par un sch\'ema en groupes affine et unipotent. Le faisceau
$\pi_{1}(F,*)$ est un sous-faisceau d'un faisceau repr\'esentable par un sch\'ema en groupes affine.
\end{thm}

Le th\'eor\`eme
pr\'ec\'edent r\'epond presque \`a la conjecture \ref{conj}. Pour donner une preuve compl\`ete \`a cette conjecture
il resterait \`a montrer que pour toute $\infty$-gerbe affine point\'ee $F$, le faisceau $\pi_{1}(F,*)$
soit repr\'esentable par un sch\'ema affine. On peut en r\'ealit\'e montrer que le faisceau $\pi_{1}(F,*)$
poss\`ede un espace de modules grossier affine $\widetilde{\pi_{1}(F,*)}$, et que le morphisme naturel
$\pi_{1}(F,*) \longrightarrow \widetilde{\pi_{1}(F,*)}$ est un monomorphisme. \\

Les types d'homotopie sch\'ematiques forment une cat\'egorie de champs
donnant une solution du probl\`eme de la sch\'ematisation (au sens de l'appendice A).
De fa\c{c}on plus pr\'ecise,
pour un corps $k$, on peut consid\'erer le foncteur d\'eriv\'e des sections globales
$\mathbb{R}\Gamma$, de la cat\'egorie des types d'homotopie sch\'ematiques point\'es sur $k$ vers la cat\'egorie homotopique
des ensembles simpliciaux point\'es. Notre second th\'eor\`eme d'existence est le suivant.

\begin{thm}{(voir Thm. \ref{t7})}
Le foncteur $\mathbb{R}\Gamma$ restreint \`a la cat\'egorie des types d'homotopie sch\'ematiques point\'es poss\`ede
un adjoint \`a gauche $X \mapsto (X\otimes k)^{sch}$. Le champ $(X\otimes k)^{sch}$ sera appel\'e une
\textit{sch\'ematisation de $X$ sur $k$}.
\end{thm}

On d\'eduit imm\'ediatement de la propri\'et\'e universelle des sch\'ematisations que le sch\'ema en groupes
affine $\pi_{1}((X\otimes k)^{sch},x)$ est le compl\'et\'e affine du groupe discret $\pi_{1}(X,x)$. Il s'agit donc
du \textit{groupe fondamental alg\'ebrique de l'espace $X$} (d\'ecrit par exemple dans \cite{d1}
et \cite[VI \S 1]{sa}).
On entrevoit ici la diff\'erence
fondamentale entre les foncteurs d'affination et de sch\'ematisation. Il existe en effet toujours un morphisme naturel
$(X\otimes k)^{sch} \longrightarrow (X\otimes k)^{uni}$, qui au niveau des groupes fondamentaux induit la projection
de $\pi_{1}((X\otimes k)^{sch},x)$ sur son quotient unipotent maximal. Ce morphisme au niveau des champs
peut aussi \^etre raisonnablement pens\'e comme un \textit{quotient unipotent maximal}.

Tout comme pour le cas des affinations, le morphisme d'adjonction $X \longrightarrow (X\otimes k)^{sch}$ poss\`ede une
certaine propri\'et\'e de pr\'eservation de la cohomologie. Plus pr\'ecis\`ement, le groupe $\pi_{1}(X,x)$ et
le faisceau en groupes $\pi_{1}((X\otimes k)^{sch},x)$ poss\`edent les m\^emes re\-pr\'e\-sen\-ta\-tions lin\'eaires. De plus,
pour $V$ une telle re\-pr\'e\-sen\-ta\-tion lin\'eaire on a
$H^{*}(X,V)\simeq H^{*}((X\otimes k)^{sch},V\otimes \mathbb{G}_{a})$. Ces deux conditions
peuvent aussi se combiner en une notion de \textit{$P$-\'equivalence} (voir Def. \ref{d12}), ce qui permet d'affirmer que
le morphisme $X \longrightarrow (X\otimes k)^{sch}$ est aussi un morphisme de localisation au sens de A.K. Bousfield.

Enfin, tout comme pour le cas du foncteur d'affination, le
foncteur $X \mapsto (X\otimes k)^{sch}$ permet de modeler une certaine partie de la th\'eorie de l'homotopie (une partie
que l'on pourra appeler \textit{sch\'ematique}). Cependant, m\^eme pour des types d'homotopie de type fini $X$,
il n'est pas raisonnable de penser pouvoir d\'ecrire explicitement le champ $(X\otimes k)^{sch}$. Ainsi, ce que mod\'ele
exactement les types d'homotopie sch\'ematiques reste un peu myst\'erieux. Quoiqu'il en soit, des exemples simples
montrent que le champ $(X\otimes k)^{sch}$ est g\'en\'eralement beaucoup plus gros que le champ $(X\otimes k)^{uni}$ (voir
$\S$ $3.4$). Enfin, lorsque $X$ est simplement connexe on a $(X\otimes k)^{sch}\simeq (X\otimes k)^{uni}$.
L'objet $(X\otimes k)^{sch}$ nous semble donc un substitut bien adapt\'e \`a la th\'eorie de
l'homotopie rationnelle et $p$-adique pour les types d'homotopie non-simplement connexes. A titre d'exemple
d'application, la construction $X \mapsto (X\otimes\mathbb{C})^{sch}$ sera utilis\'ee dans un travail
ult\'erieur pour \'etudier les types d'homotopie des vari\'et\'es complexes lisses et projectives \`a groupes
fondamentaux arbitraires. L'\'etude par la th\'eorie de Hodge de ce nouvel
invariant d'homotopie permettra alors d'obtenir de nouvelles restrictions sur les types d'homotopie
des vari\'et\'es projectives, qui semblent hors d'atteinte par une approche utilisant la th\'eorie
de l'homotopie rationnelle (voir \cite{kt} ainsi que \S 3.5.1). \\

Signalons aussi que les alg\`ebres co-simpliciales sont des mod\`eles alg\'ebriques pour les champs affines. De m\^eme,
il existe une notion de \textit{$H_{\infty}$-alg\`ebres de Hopf} qui permet de donner des mod\`eles alg\'ebriques
des types d'homotopie sch\'ematiques (voir Def. \ref{hopf}). L'\'etude de ces mod\`eles alg\'ebriques ne sera pas le but de
ce travail, mais nous esp\'erons pouvoir revenir sur le sujet par la suite
(voir \cite{small}). A titre indicatif, on peut r\'esumer
le double apport de cet article par le diagramme symbolique suivant.
$$
\xymatrix{
  \{ \text{Champs  affines} \} \ar[r] & \{\text{Types d'homotopie sch\'ematiques} \} \\
 \{ \text{Alg\`ebres co-simpliciales} \} \ar[r] \ar[u] & \ar[u] \{H_{\infty}-\text{alg\`ebres
  de  Hopf}\} } $$
Dans ce diagramme, les fl\`eches verticales sont les proc\'ed\'es de g\'eom\'etrisation induit par le foncteur
$\mathbb{R}Spec$ (i.e. le 
passage de l'alg\`ebre \`a la g\'eom\'etrie), et les fl\'eches horizontales symbolisent le passage aux \textit{objets en groupes}. \\

\begin{center} \textit{Types d'homotopie des vari\'et\'es alg\'ebriques} \end{center}

La th\'eorie g\'en\'erale des types d'homotopie sch\'ematiques et du foncteur
de sch\'ematisation poss\`ede plusieurs applications
dans le cadre de la g\'eom\'etrie alg\'e\-bri\-que. Une premi\`ere application est l'\'etude
des types d'homotopie des vari\'et\'es alg\'ebriques complexes (voir \S 3.5.1). En effet, pour
une vari\'et\'e complexe lisse et projective $X$ et $x\in X(\mathbb{C})$, on peut consid\'erer son espace
topologique sous-jacent des points complexes (muni de la topologie analytique) $X^{top}$, ainsi
que sa sch\'ematisation $(X^{top}\otimes \mathbb{C})^{sch}$. Dans le travail
\cite{kt} on construit une \textit{d\'ecomposition de Hodge} sur le champ
$(X^{top}\otimes \mathbb{C})^{sch}$, qui englobe les diff\'erentes structures de Hodge
connues sur la cohomologie, les groupes d'homotopie et le groupe fondamental. L'existence de
cette structure additionelle sur le champ $(X^{top}\otimes \mathbb{C})^{sch}$ poss\`ede
d'int\'eressantes cons\'equences sur le type d'homotopie de la vari\'et\'e $X$
et impose de s\'erieuses restrictions. On peut ainsi construire des nouveaux exemples
de types d'homotopie qui ne sont pas r\'ealisables par des vari\'et\'es projectives, et dont
l'obstruction se trouve dans des invariants d'homotopie sup\'erieurs (pr\'ecisemment
l'action du groupe fondamental sur les groupes d'homotopie).

De fa\c{c}on plus g\'en\'erale pour
une vari\'et\'e alg\'ebrique $X$ sur un corps de base $k$ on peut utiliser la notion
de types d'homotopie sch\'ematiques afin de construire des types d'homotopie
associ\'es aux th\'eories cohomologiques usuelles (Betti, de Rham, $l$-adique, crystalline, \dots, voir
\S 3.5.2, 3.5.3). Toutes ces constructions montrent que la notion de type d'homotopie sch\'ematique
est ad\'equate pour l'\'etude homotopique des vari\'et\'es alg\'ebriques, tout comme
il \'etait d\'ej\`a bien connu que le bon cadre pour une th\'eorie des groupes fondamentaux de vari\'et\'es
alg\'ebriques est en r\'ealit\'e celui des sch\'emas en groupes affines et non celui des groupes discrets
(voir par exemple \cite{d1,sa}). Cependant, ces constructions ne seront pas d\'ecrites en d\'etail
et nous nous contenterons de donner les propri\'et\'es fondamentales (et parfois
caract\'eristiques) de ces types d'homotopie. \\

\begin{center} \textit{Champs $\infty$-g\'eom\'etriques} \end{center}

Une autre application de la th\'eorie des champs affines est la notion
de \textit{champ $\infty$-g\'eom\'etrique} (voir \S 4). Par d\'efinition un champ alg\'ebrique
(par exemple au sens d'Artin) est obtenu en recollant des sch\'emas affines
\`a l'aide de l'action d'un groupoide lisse. En rempla\c{c}ant formellement la notion
de sch\'ema affine par celle plus g\'en\'erale de champ affine on donne une d\'efinition
de \textit{champ $\infty$-g\'eom\'etrique}. Il se trouve que les types
d'homotopie sch\'ematiques fournissent les premiers exemples de champs $\infty$-g\'eom\'etriques.
Un deuxi\`eme exemple fondamental est le champ des structures multiplicatives
sur un complex parfait $V$, g\'en\'eralisant le champ des structures d'alg\`ebres
sur un espace vectoriel donn\'e au cas des alg\`ebres diff\'erentielles gradu\'ees (voir
\S 4.2.1).
On peut aussi donner des exemples plus sophistiqu\'es de champs $\infty$-g\'eom\'etriques,
comme le champ des filtrations sur une alg\`ebre diff\'erentielle gradu\'ee commutative, qui
pourra \^etre utilis\'e pour construire un \textit{domaine des p\'eriodes non-ab\'eliennes} ainsi
qu'une application de p\'eriodes correspondantes qui controlera la variation
de la filtration de Hodge sur le type d'homotopie rationnel d'une famille de vari\'et\'es lisses et
projectives (voir \S 4.2.2). \\

\begin{center} \textit{Relations avec d'autres travaux} \end{center}

La notion de champs affine me semble nouvelle. Il en existe
cependant des avatars ou versions voisines dans les articles \cite[\S 6]{s5} et \cite{hin2}.
En caract\'eristique nulle et pour le cas des
types d'homotopie $1$-connexes et de type fini le th\'eor\`eme \ref{ti} apparait
dans \cite[\S 6]{s5}.

Il y a aussi les quelques six cents pages de \cite{gr}, dont il est extr\`emement difficile de faire une comparaison
avec le point de vue adopt\'e ici. Par exemple,
la notion de champs affines me semble assez proche de celle de \textit{complexes of unipotent bundles}
discut\'es dans \cite[p. $446$]{gr}. Cependant,
l'objet que A. Grothendieck appelle \textit{sch\'ematisation} d'un espace $X$
semble plus proche de notre affination $(X\otimes k)^{uni}$ que de notre sch\'ematisation $(X\otimes k)^{sch}$. A vrai dire
la question de l'existence du champ $(X\otimes k)^{sch}$, bien que fortement inspir\'ee des consid\'erations sur la th\'eorie de
l'homotopie que l'on  trouve dans les lettres de A. Grothendieck \`a L. Breen, ne semble pas avoir \'et\'e consid\'er\'ee dans
\cite{gr}.  Cependant, la motivation principale de ce travail est n\'ee d'un d\'esir de  comprendre certains passages de \cite{gr},
et de ce fait ce manuscrit a eu une
influence d\'eterminante sur les d\'efinitions et les \'enonc\'es qui apparaissent dans cet article.
Ce travail s'ins\`ere donc naturellement dans le vaste programme propos\'e par C. Simpson
et intitul\'e \textit{la poursuite de la poursuite des champs} (ou encore comme l'a sugg\'er\'e D. Husem\"oller lors d'une
conversation sur le sujet \textit{la $2$-poursuite des champs}).

Il existe aussi d'autres approches au probl\`eme de la sch\'ematisation des types
d'homotopie tel qu'il est abord\'e dans \cite{gr}. La premi\`ere, qui est expos\'ee dans \cite{be}, donne une
construction de la \textit{sph\`ere sch\'ematique stable} en caract\'eristique positive. Il n'est pas tout \`a fait clair que cette
construction soit comparable \`a celle de cet article, m\^eme si cette sph\`ere sch\'ematique stable ressemble de tr\`es pr\`es
\`a ce que l'on pourrait appeler dans notre langage un \textit{spectre affine} (version stable des champs affines). Une
seconde
approche a \'et\'e introduite dans \cite{s3}. Ici aussi la comparaison avec notre construction n'est pas imm\'ediate, bien
que la th\'eorie des \textit{champs pr\'esentables et tr\`es pr\'esentables} soit intimement li\'ee \`a celle
des types d'homotopie sch\'ematiques (voir Thm. \ref{t6}). Enfin, dans \cite{to2} nous avons construit
un foncteur de sch\'ematisation \`a l'aide d'une notion de \textit{cat\'egorie simpliciale Tannakienne}, qui est
conjectur\'e \^etre isomorphe \`a celui d\'efini dans ce travail (voir
aussi \cite{to3}). Signalons au passage que les constructions et
r\'esultats du pr\'esent travail ont \'et\'es pour la plupart devin\'es \`a l'aide de ce formalisme Tannakien.

Comme nous l'avons fait remarquer le foncteur de sch\'ematisation que nous construisons
peut \^etre consid\'er\'e comme un substitut des th\'eories de l'homotopie rationnelle et $p$-adique
pour des espaces non-nilpotents. Il existe aussi une autre approche qui consiste \`a utiliser
une version \'equivariante de l'homotopie rationnelle et $p$-adique, et \`a appliquer celle-ci aux rev\^etements universels
(voir \cite{bs,ght}). Il me semble un probl\`eme int\'eressant de savoir comment ces deux points de vue
sont reli\'es.

Enfin, la notion de champs $\infty$-g\'eom\'etriques pr\'esent\'e dans le paragraphe
\S 4 a \'et\'e inspir\'ee par \cite{s5}. Elle est aussi
tout \`a fait analogue a la notion de \emph{$D$-champs g\'eom\'etriques}
de \cite{hagdag} (voir aussi \cite{tv2}). Pour tout dire la th\'eorie des champs affines
de cet article a fortement inspir\'e le d\'eveloppement de la
\emph{g\'eom\'etrie alg\'ebrique homotopique} de \cite{tv,tv2}. \\

\begin{center} \textit{Remerciements} \end{center}

Je voudrais tout d'abord remercier C. Simpson, qui par l'interm\'ediaire d'ar\-ti\-cl\-es, de discussions et de
correspondances a beaucoup influenc\'e ce travail. Je remercie aussi particuli\`erement H. Baues,
qui m'a sugg\'er\'e d'utiliser une
version d\'eriv\'ee du foncteur de compl\'etion affine des groupes discrets afin de construire
le foncteur de sch\'ematisation.
Le lecteur remarquera que ceci est une \'etape fondamentale dans la preuve du th\'eor\`eme \ref{t7}.

Pour de nombreuses discussions, remarques et correspondances sur le sujet je remercie aussi les personnes suivantes:
L. Breen, A. Hirschowitz, L. Katzarkov, M. Olsson,
T. Pantev, J. Sauloy, M. Spitzweck, D. Stanley, J. Tapia et G. Vezzosi.

La majeure partie de ce travail a \'et\'e effectu\'ee \`a l'institut Max Planck de Bonn durant l'ann\'ee $1999-2000$, que je tiens
\`a remercier pour son hospitalit\'e et ses conditions de travail exceptionnelles.

\newpage

\textsf{Notations et conventions:}
Pour un univers $\mathbb{U}$ nous appellerons $\mathbb{U}$-ensemble (resp. $\mathbb{U}$-ensemble simplicial, resp. $\mathbb{U}$-groupe, resp. \dots) un ensemble (resp.
un ensemble simplicial, resp. un groupe, resp. \dots) qui est un \'el\'ement de $\mathbb{U}$. La cat\'egorie des $\mathbb{U}$-ensembles
(resp. $\mathbb{U}$-ensembles simpliciaux, resp. $\mathbb{U}$-groupes, resp \dots) sera alors not\'ee $\mathbb{U}-Ens$ (resp.
$\mathbb{U}-SEns$, resp. $\mathbb{U}-Gp$, resp. \dots). Nous ferons une exception avec l'expression $\mathbb{U}$-cat\'egorie, que nous r\'eservons \`a la notion usuelle
(voir \cite[$I Def. 1.2$]{sga4}), d\'esignant une cat\'egorie $C$ telle que \'etant donn\'e deux objets $X$ et $Y$, l'ensemble des morphismes
de $X$ vers $Y$ dans $C$ soit un $\mathbb{U}$-ensemble. Pour d\'esigner une cat\'egorie appartenant \`a $\mathbb{U}$ nous parlerons de cat\'egorie
$\mathbb{U}$-petite.

Nous parlerons de $\mathbb{U}$-limites et $\mathbb{U}$-colimites (resp. $\mathbb{V}$-limites et $\mathbb{V}$-colimites)
pour d\'esigner des limites et colimites dont les cat\'egories d'indices
appartiennent \`a $\mathbb{U}$ (resp. \`a $\mathbb{V}$). De m\^eme, nous parlerons de $\mathbb{U}$-limites et
$\mathbb{U}$-colimites homotopiques (resp. $\mathbb{V}$-limites et
$\mathbb{V}$-colimites homotopiques) \`a valeurs dans une cat\'egorie de mod\`eles (voir
\cite[$\S 20$]{hi}). \\

Pour la suite nous fixerons $\mathbb{U}$ un univers contenant l'ensenbles des nombres naturels, et $\mathbb{V}$ un univers avec $\mathbb{U} \in
\mathbb{V}$.
Comme il en est l'usage, nous noterons $\Delta$ la cat\'egorie simpliciale standard. C'est la cat\'egorie $\mathbb{U}$-petite dont les objets sont les
nombres naturels $[n]$, et dont les morphismes de $[m]$ vers $[n]$ sont les applications croissantes de $\{0,\dots,m\}$ vers $\{0,\dots,n\}$.
Pour tout $\mathbb{V}$-cat\'egorie $C$, nous noterons $SC$ la $\mathbb{V}$-cat\'egorie des objets simpliciaux dans $C$ (i.e. la $\mathbb{V}$-cat\'egorie des foncteurs de $\Delta^{o}$ vers $C$), et pour $F \in SC$ nous noterons $F_{n}:=F([n])$. Nous identifierons syst\'ematiquement $C$ \`a la sous-cat\'egorie
pleine de $C^{\Delta^{o}}$ form\'ee des foncteurs constants $F : \Delta^{o} \longrightarrow C$.

Si $C$ est une cat\'egorie $\mathbb{V}$-petite, et $X$ un objet de $C$, nous noterons $h_{X}$ le $\mathbb{V}$-pr\'efaisceau simplicial sur $C$ qu'il repr\'esente. C'est donc
le foncteur d\'efini par $h_{X}(Y):=Hom(Y,X)$, o\`u l'ensemble $Hom(Y,X)$ est vu comme un ensemble simplicial constant. De m\^eme, si $X$ est
un objet simplicial de $C$, nous noterons $h_{X}$ le $\mathbb{V}$-pr\'efaisceau simplicial dont le pr\'efaisceau des simplexes de dimension $n$ est
d\'efini par la formule $(h_{X})_{n}(Y):=Hom(Y,X_{n})$. \\

Nous utiliserons les d\'efinitions de cat\'egories de mod\`eles ferm\'ees \'enonc\'ees dans \cite{ho}. Elles seront toujours $\mathbb{V}$-petites.
On supposera implicitement lorsqu'elles seront engendr\'ees par cofibrations (voir \cite[$\S 2.1$]{ho}) que les ensembles g\'en\'erateurs
des cofibrations et cofibrations triviales seront des $\mathbb{V}$-ensembles (de m\^eme, les ordinaux apparaissant dans les colimites
transfinies seront des \'el\'ements de $\mathbb{V}$). Rappelons que la cat\'egorie $\mathbb{V}-SEns$ est une cat\'egorie de mod\`eles ferm\'ee, mono\"{\i}dale sym\'etrique
pour le produit direct, et engendr\'ee par cofibrations pour la convention ci-dessus (voir \cite[$\S 3$]{ho}).
Les $Hom$ internes de $\mathbb{V}-SEns$ seront not\'es
$\underline{Hom}$, et ces $Hom$ d\'eriv\'es $\mathbb{R}\underline{Hom}$ (voir \cite[$\S 4$]{ho}).
Une cat\'egorie de mod\`eles ferm\'ee simpliciale sera alors une cat\'egorie ($\mathbb{V}$-petite) de mod\`eles ferm\'ee,
qui est un module
sur $\mathbb{V}-SEns$ (au sens de \cite[$4.2.28$]{ho}). Enfin, si $M$ est une cat\'egorie de mod\`eles ferm\'ee $Ho(M)$
d\'esignera sa cat\'eorie homotopique. C'est une cat\'egorie $\mathbb{V}$-petite qui est
naturellement \'equivalente \`a une $\mathbb{U}$-cat\'egorie. Les ensembles de morphismes
dans $Ho(M)$ seront not\'es $[-,-]_{M}$.

Enfin, pour une cat\'egorie de mod\`eles ferm\'ee $M$ nous noterons $M_{*}$ la cat\'egorie de mod\`eles de ses objets point\'es.
Lorsque
de plus $M$ est une cat\'egorie de mod\`eles simpliciale, il en est de m\^eme de $M_{*}$. Ses $Hom$ simpliciaux seront alors not\'es $\underline{Hom}_{*}$ (voir \cite[$4.2.19$]{ho}).

\newpage

\section{Rappels sur les champs}

Dans ce premier chapitre nous avons rassembl\'e un certains nombres de d\'efinitions et de r\'esultats standards de la th\'eorie
des pr\'efaisceaux simpliciaux sur un site de Grothendieck. Les id\'ees essentielles de la th\'eorie ont \'et\'e d\'evelopp\'ees
par A. Joyal dans une lettre \`a A. Grothendieck (voir \cite{jo}), o\`u l'existence d'une structure de cat\'egorie de mod\`eles
sur la cat\'egorie des faisceaux simpliciaux est d\'emontr\'ee. Nous utiliserons cependant une structure un peu diff\'erente, o\`u la classe
des cofibrations est strictement plus petite que celle des monomorphismes, et o\`u les objets fibrants sont faciles \`a caract\'eriser
(voir \ref{l1}).
Bas\'ee sur des id\'ees remontant \`a D. Quillen, A.K. Bousfield et D.M. Kan, cette structure \`a \'et\'e consid\'er\'ee
pour la premi\`ere fois dans \cite[$5$]{s1}, et a \'et\'e reprise en d\'etail dans \cite{bla}.
Pour ce qui est des preuves et des d\'etails techniques nous renvoyons le lecteur \`a
\cite{j,s1,gj,s1,bla}. Signalons aussi le travail en cours \cite{du}, o\`u il est probable que le point de vue de la
localisation de  Bousfield d\'ecrit dans \cite[$6$]{s1} soit repris en d\'etail.

Nous commencerons par rappeler la d\'efinition de la structure de cat\'egorie de mod\`eles ferm\'ee sur la cat\'egorie
des pr\'efaisceaux simpliciaux, ainsi que ces principales propri\'et\'es. L'existence de cette structure de cat\'egorie de mod\`eles
nous sera utile tout au long de ce travail. Dans les paragraphes suivants nous rappellerons comment elle permet de d\'efinir
la cohomologie d'un pr\'efaisceau simplicial
\`a valeurs dans un syst\`eme local. Nous donnerons aussi une br\`eve d\'emonstration d'un th\'eor\`eme de comparaison avec la
cohomologie d\'efinie par foncteurs d\'eriv\'es, ainsi qu'un crit\`ere pour s'assurer qu'un pr\'efaisceau simplicial poss\`ede une
d\'e\-com\-po\-si\-tion de Postnikov convenable. Nous terminerons ce chapitre par un aper\c{c}u de la th\'eorie du d\'ela\c{c}age dans le cadre
des pr\'efaisceaux simpliciaux, et un bref rappel sur les sch\'emas en groupes affines sur un corps.  \\

Pour tout ce chapitre $C$ d\'esignera une cat\'egorie $\mathbb{V}$-petite, et $SPr(C)$ la cat\'egorie des pr\'efaisceaux
en $\mathbb{V}$-ensembles
simpliciaux sur $C$. La cat\'egorie $SPr(C)$ est donc une $\mathbb{V}$-cat\'egorie.
Nous supposerons que $C$ poss\`ede des produits fibr\'es, et qu'elle est munie d'une topologie de Grothendieck
(au sens de \cite[$IV$]{sga3}). Pour fixer les id\'ees, notons que par la suite $C$ sera la cat\'egorie des $\mathbb{U}$-sch\'emas affines sur $Spec\, k$, pour $k$ un anneau commutatif, et la topologie
sera engendr\'ee par les morphismes fid\`element plats et quasi-compacts (i.e. la topologie $fpqc$ de \cite[$IV.6.3$]{sga3}).

\subsection{Rappel des d\'efinitions}

On commence par munir la cat\'egorie $SPr(C)$ d'une structure de cat\'egorie de mod\`eles que nous appellerons le structure forte
(dans la litt\'erature on trouve aussi l'expression de structure projective niveau par niveau), et qui
a \'et\'ee introduite pour la premi\`ere fois dans \cite{bk}.
Il s'agit du cas
o\`u l'on consid\`ere la topologie discr\`ete, et donc o\`u l'on voit $C$ comme une simple cat\'egorie d'indices.
On appellera alors fibration forte (resp. \'equivalence forte) tout morphisme de $SPr(C)$,
$f : F \longrightarrow F'$, tel que pour tout objet $X \in C$ le morphisme $f_{X} : F(X) \longrightarrow F(X')$ soit une fibration
(resp. une \'equivalence) dans $\mathbb{V}-SEns$. Pour ces d\'efinitions, $C$ \'etant $\mathbb{V}$-petite, la cat\'egorie $SPr(C)$ est une cat\'egorie de mod\`eles
ferm\'ee simpliciale (voir \cite[$Ex. II.6.9$]{gj}).
Les $Hom$ simpliciaux de cette derni\`ere seront not\'es $\underline{Hom}$. En clair, ceci signifie que pour
$X \in \mathbb{V}-SEns$ et $F,F' \in SPr(C)$ on peut d\'efinir fonctoriellement dans $SPr(C)$ des objets
$X\otimes F$, $\underline{Hom}(F,F')$ et $(F')^{X}$ satisfaisant \`a la formule d'adjonction
$$\underline{Hom}(X,\underline{Hom}(F,F'))\simeq \underline{Hom}(X\otimes F,F')\simeq
\underline{Hom}(F,(F')^{X}).$$
Ces $Hom$ simpliciaux \'etant compatibles avec la structure de cat\'egorie de mod\`eles on peut d\'efinir des $Hom$
simpliciaux d\'eriv\'es (voir \cite[$\S 4.3$]{ho})
$$\mathbb{R}_{f}\underline{Hom}(F,F'):=\underline{Hom}(QF,RF') \in Ho(\mathbb{V}-SEns),$$
o\`u $QF \longrightarrow F$ est un remplacement cofibrant, et $F' \longrightarrow RF'$ un remplacement fibrant.
Noter l'indice $f$ dans la notation $\mathbb{R}_{f}$, qui fait r\'ef\'erence \`a la structure forte. Ces objets sont
bien d\'efinis et fonctoriels comme objets de la cat\'egorie homotopique $Ho(\mathbb{V}-SEns)$. De plus, la formule d'adjonction pr\'ec\'edente
induit une formule d'adjonction d\'eriv\'ee
$$\mathbb{R}_{f}\underline{Hom}(X,\mathbb{R}_{f}\underline{Hom}(F,F'))\simeq \mathbb{R}_{f}\underline{Hom}(X\otimes F,F')
\simeq \mathbb{R}_{f}\underline{Hom}(F,(RF')^{X}),$$
o\`u $F' \longrightarrow RF'$ est un remplacement fibrant. \\

Pour $X \in C$, on dispose du foncteur d'\'evaluation en $X$
$$\begin{array}{cccc}
j_{X} : & SPr(C) & \longrightarrow & \mathbb{V}-SEns \\
& F & \mapsto & F(X).
\end{array}$$
Ce foncteur poss\`ede un adjoint \`a gauche not\'e $j_{X}^{*}$, et l'adjonction $(j_{X}^{*},j_{X})$ est une adjonction de Quillen.
En particulier, pour tout $X \in C$, $j_{X}^{*}(*)$ est un objet cofibrant de $SPr(C)$. Or, il est clair que
$j_{X}^{*}(*)$ est isomorphe \`a $h_{X}$. Ainsi, pour tout objet du site $X \in C$, le pr\'efaisceau simplicial qu'il repr\'esente
$h_{X}$ est un objet cofibrant.
Dans le cas o\`u $X=e$ est un objet final, nous noterons comme il en est l'usage
$$\begin{array}{cccc}
\Gamma:=j_{e} : & SPr(C) & \longrightarrow & \mathbb{V}-SEns \\
& F & \mapsto & F(e),
\end{array}$$
le foncteur des sections globales. Son adjoint \`a gauche $j_{e}^{*}$ est le foncteur qui associe \`a $X \in \mathbb{V}-SEns$ le pr\'efaisceau
constant $\underline{X}$. Remarquer que lorsque $C$ ne poss\`ede pas d'objet final les pr\'efaisceaux constants ne sont g\'en\'eralement pas
cofibrants. De plus, le foncteur des sections globales $\Gamma(-):=\underline{Hom}(*,-)$ n'est alors plus de Quillen \`a droite en g\'en\'eral. \\

Revenons au cas o\`u $C$ est munie d'une topologie de Grothendieck. A un objet $F \in SPr(C)$ on associe ses pr\'efaisceaux d'homotopie
de la fa\c{c}on suivante. Pour tout $X \in C$, et tout $0$-simplexe $s \in F(X)_{0}$, on d\'efinit le $m$-\`eme pr\'efaisceau d'homotopie
de $F$ point\'e en $s$ par
$$\begin{array}{cccc}
\pi_{m}^{pr}(F,s) : & (C/X)^{o} & \longrightarrow & \mathbb{V}-Gp \\
& (u : Y\rightarrow X) & \mapsto & \pi_{m}(F(Y),u^{*}(s)).
\end{array}$$
C'est un pr\'efaisceau en $\mathbb{V}$-groupes sur $C/X$ qui est ab\'elien lorsque $m>1$. Le pr\'efaisceau des composantes connexes de $F$
est d\'efini par
$$\begin{array}{cccc}
\pi_{0}^{pr}(F) : & C^{o} & \longrightarrow & \mathbb{V}-Ens \\
& X & \mapsto & \pi_{0}(F(X)).
\end{array}$$
On d\'efinit alors les faisceaux d'homotopie de $F$ comme \'etant les faisceaux associ\'es aux pr\'efaisceaux d'homotopie. Ils seront not\'es
$$\pi_{m}(F,s):=a(\pi_{m}^{pr}(F,s)) \qquad \pi_{0}(F):=a(\pi_{0}^{pr}(F)).$$

\begin{df}\label{d15}
\emph{Soit $f : F \longrightarrow F'$ un morphisme de pr\'efaisceaux simpliciaux sur le site $C$.}
\begin{itemize}
\item
\emph{Le morphisme $f$ est une \'equivalence si
les deux conditions suivantes sont satisfaites.}
\begin{itemize}
\item \emph{Le morphisme induit $\pi_{0}(f) : \pi_{0}(F) \longrightarrow \pi_{0}(F')$ est un isomorphisme de 
faisceaux sur $C$.}
\item \emph{Pour tout objet $X \in C$, tout $m>0$ et tout $0$-simplexe $s \in F(X)_{0}$, les morphismes induits
$$\pi_{m}(f,s) : \pi_{m}(F,s) \longrightarrow \pi_{m}(F',f(s))$$
sont des isomorphismes de faisceaux sur $C/X$.}
\end{itemize}
\item \emph{Le morphisme $f$ est une cofibration si c'est une cofibration forte}.
\end{itemize}
\end{df}

Comme il est expliqu\'e dans \cite[Thm. $5.1$]{s1}, ces d\'efinitions munissent la cat\'egorie
$SPr(C)$ d'une structure de cat\'egorie de mod\`eles simpliciale.
Remarquons au passage que toute \'equivalence forte est une \'equivalence, que toute fibration est une fibration forte,
et que tout fibration triviale est une \'equivalence forte. De m\^eme, toute \'equivalence entre objets fibrants dans
$SPr(C)$ est une \'equivalence forte.

Les $Hom$ simpliciaux d\'eriv\'es relativement \`a cette nouvelle structure seront not\'es
$\mathbb{R}\underline{Hom}$.
De m\^eme, lorsque $C$ poss\`edera un objet final $e$, le foncteur d\'eriv\'e \`a droite des sections globales sera not\'e $\mathbb{R}\Gamma$. On dispose aussi des isomorphismes d'adjonctions analogues \`a ceux d\'ecrits pr\'ec\'edemment
$$\mathbb{R}\underline{Hom}(X\times F,F')\simeq \mathbb{R}\underline{Hom}(X,\mathbb{R}\underline{Hom}(F,F')) \simeq
\mathbb{R}\underline{Hom}(F,(RF')^{X}).$$

Il existe une relation entre les deux structures de cat\'egorie de mod\`eles sur $SPr(C)$ connue sous le nom de
localisation de Bousfield \`a gauche (voir \cite{hi}). Nous n'insisterons pas trop sur cette notion et nous
nous contenterons de retenir que le fait que $SPr(C)$ soit une localisation de Bousfield \`a gauche de $SPr(C)$ munie de la structure forte implique
l'existence d'une jolie caract\'erisation des objets fibrants dans $SPr(C)$. Le lecteur remarquera qu'il s'agit
d'un analogue de nature homotopique de la  d\'efinition d'un faisceau.

\begin{lem}\label{l1}
Un objet $F \in SPr(C)$ est fibrant si et seulement si les deux assertions suivantes sont satisfaites.

\begin{enumerate}
\item Pour tout objet $X \in C$, $F(X)$ est fibrant dans $\mathbb{V}-SEns$ (i.e. $F$ est fortement fibrant).

\item Pour tout objet $X \in C$, et tout hyper-recouvrement $U_{*}$ de $X$ dans $C$, le morphisme naturel
$$F(X) \longrightarrow Holim_{[m] \in \Delta}F(U_{m})$$
est une \'equivalence dans $\mathbb{V}-SEns$.

\end{enumerate}
\end{lem}

\textit{Preuve:} Voir \cite{du}. \hfill $\Box$ \\

Un corollaire important de ce dernier lemme est qu'un pr\'efaisceau simplicial de la forme $K(M,n)$, o\`u $M$ est un
faisceau flasque sur $C$, est toujours fibrant. \\

Remarquer que s'il l'on consid\`ere des pr\'efaisceaux en groupo\"{\i}des plut\^ot que des pr\'efaisceaux simpliciaux, alors la seconde condition
du lemme pr\'ec\'edent est exactement \'equivalente \`a la condition usuelle d'\^etre un champ (voir \cite{j2} pour
plus de d\'etails sur les relations entre pr\'efaisceaux simpliciaux et champs en groupo\"{\i}des). La d\'efinition suivante
est alors naturelle. Notons qu'il s'agit aussi du point
de vue adopt\'e dans \cite{s1} pour d\'evelopper la th\'eorie des $n$-champs de Segal.

\begin{df}\label{d1}
\emph{Un champ sur $C$ est un objet de la cat\'egorie $SPr(C)$
satisfaisant \`a la condition (2) du lemme \ref{l1}. Par extension, tout objet $F$ de la cat\'egorie homotopique
$Ho(SPr(C))$ sera appel\'e un champ}.
\end{df}

Avec ce langage, le champ associ\'e \`a un pr\'efaisceau simplicial $F$ est simplement un de ses mod\`eles fibrants.

\subsection{Limites homotopiques et d\'ecomposition de Postnikov}

Soit $I$ une cat\'egorie $\mathbb{V}$-petite, et $SPr(C)^{I}$ la cat\'egorie des foncteurs de $I$ vers $SPr(C)$.
On peut munir $SPr(C)^{I}$ d'une structure de cat\'egorie de mod\`eles engendr\'ee par cofibrations
o\`u les \'equivalences (resp. les cofibrations) sont d\'efinies niveau par niveau (i.e. $F \longrightarrow F'$ est une
\'equivalence si et seulement si pour tout $i \in I$, $F(i) \longrightarrow F'(i)$ est une \'equivalence dans $SPr(C)$).
Noter que cette structure n'est pas celle que nous avons utilis\'ee pour d\'efinir $SPr(C)$, mais est une structure de type
injective analogue \`a celle d\'efinie dans \cite[Thm. $2.3$]{j}. La cat\'egorie de mod\`eles ainsi obtenue
est une cat\'egorie de mod\`eles ferm\'ee simpliciale, et ces $Hom$ simpliciaux seront not\'es $\underline{Hom}_{I}$. En particulier, on dispose du
foncteur des sections globales
$$\begin{array}{cccc}
\Gamma(I,-) : & SPr(C)^{I} & \longrightarrow & \mathbb{V}-SEns \\
& F & \mapsto & \underline{Hom}_{I}(*,F)
\end{array}$$
dont l'adjoint \`a gauche est le foncteur qui associe \`a $X \in \mathbb{V}-SEns$ le foncteur constant de valeurs $X$, $\underline{X} : I \longrightarrow SPr(C)$.
L'adjonction $(\underline{-},\Gamma(I,-))$ est alors une adjonction de Quillen, et
le foncteur d\'eriv\'e \`a droite de $\Gamma(I,-)$ sera alors not\'e
$$Holim_{I}:=\mathbb{R}\Gamma(I,-) : Ho(SPr(C)^{I}) \longrightarrow Ho(SPr(C)).$$
C'est le foncteur de limite homotopique le long de $I$. Il est naturellement isomorphe au foncteur $Holim_{I}$ construit
de fa\c{c}on standard \`a partir de la structure simpliciale sur $SPr(C)$ (voir \cite[$\S 19$]{hi}).

Il est important de remarquer que la d\'efinition de $Holim_{I}$ d\'epend de la topologie de $C$, et en g\'en\'eral, la comparaison entre
les limites homotopiques pour diff\'erentes topologies n'est pas ais\'ee. On dispose cependant du crit\`ere g\'en\'eral suivant.
Il permettra par la suite d'avoir dans certains cas des d\'ecompositions de Postnikov
raisonnables, que l'on sait ne pas exister en g\'en\'eral.

\begin{lem}\label{l2}
Soit $F : I \longrightarrow SPr(C)$ un $I$-diagramme de pr\'efaisceaux simpliciaux, tel que pour tout $i \in I$, $F_{i}$ soit un champ.
Alors, si on note $Holim_{I}^{triv}F \in Ho(SPr(C))$ la limite homotopique de $F$ lorsque $C$ est muni de la topologie triviale, le morphisme
naturel
$$Holim_{I}^{triv}F \longrightarrow Holim_{I}F$$
est un isomorphisme dans $Ho(SPr(C))$.
\end{lem}

\textit{Preuve:} C'est imm\'ediat, car toute r\'esolution fibrante $F \longrightarrow F'$ est telle que les morphismes induits
$F_{i} \longrightarrow F_{i}'$ soient des \'equivalences fortes pour tout $i \in I$.
Ainsi, $F \longrightarrow F'$ est aussi une r\'esolution fibrante lorsque
$C$ est muni de la topologie triviale. \hfill $\Box$ \\

Nous allons maintenant \'etudier deux cas particulier du foncteur pr\'ec\'edent. Il s'agit du cas o\`u
$I$ poss\`ede trois objets $0$, $1$ et $2$ et deux uniques morphismes non triviaux $1 \rightarrow 0$, $2 \rightarrow 0$,
et du cas o\`u
$I$ est la cat\'egorie associ\'ee \`a l'ensemble totalement ordonn\'e des entiers naturels. \\

Commen\c{c}ons par examiner le premier cas, et notons $0$, $1$ et $2$ les objets de $I$, et $a : 1 \rightarrow 0$, $b : 2 \rightarrow 0$
ses morphismes non triviaux. Un objet de $SPr(C)^{I}$ est alors la donn\'ee d'un diagramme de pr\'efaisceaux simpliciaux
$$\xymatrix{F_{1} \ar[r]^{a} & F_{0} & \ar[l]_{b} F_{2}.}$$
Pour tout objet $F \in Ho(SPr(C)^{I})$ nous utiliserons la notation
$$F_{1}\times^{h}_{F_{0}}F_{2}:=Holim_{I}F,$$
que nous appellerons le produit fibr\'e homotopique de $F_{1}$ et $F_{2}$ au-dessus de $F_{0}$. Un cas particuli\`erement important est celui o\`u $F_{2}=*$. L'objet $F_{1}\times^{h}_{F_{0}}*$ est alors la fibre
homotopique du morphisme $a : F_{1} \longrightarrow F_{0}$ au point $b : * \longrightarrow F_{0}$.

Il est facile de v\'erifier que les objets fibrants dans $SPr(C)^{I}$ sont exactement les diagrammes tels que $F_{i}$ soit un champ pour tout $i$,
et tels que  les morphismes $a$ et $b$ soient des fibrations. En particulier, si $F \in SPr(C)^{I}$ est fibrant, pour tout $i$ et tout $X \in C$,
$F_{i}(X)$  est un ensemble simplicial fibrant, et les morphismes $a_{X} : F_{1}(X) \rightarrow F_{0}(X)$, $b_{X} : F_{2}(X) \rightarrow F_{0}(X)$
 sont des fibrations d'ensembles simpliciaux.

Supposons maintenant que $F \in SPr(C)^{I}$ soit fibrant et muni d'un morphisme $s : * \longrightarrow F$, qui induit donc des morphismes
$s_{i} : * \longrightarrow F_{i}$, et $s : * \longrightarrow F_{1}\times^{h}_{F_{0}}F_{2}$. Alors, en utilisant les suites exactes longues d'homotopie associ\'ees aux produits fibr\'es
homotopiques d'ensembles simpliciaux, on obtient une suite exacte longue de pr\'efaisceaux en groupes sur $C$
$$\xymatrix{
\pi_{m+1}^{pr}(F_{0},s_{0}) \ar[r]^-{\partial} & \pi_{m}^{pr}(F_{1}\times^{h}_{F_{0}}F_{2},s) \ar[r]^-{a\times b} &
\pi_{m}^{pr}(F_{1},s_{1})\times \pi_{m}^{pr}(F_{2},s_{2}) \ar[r]^-{a-b} & }$$ 
$$\xymatrix{ \pi_{m}^{pr}(F_{0},s_{0}) \ar[r]^-{\partial} & \pi_{m-1}^{pr}(F_{1}\times^{h}_{F_{0}}F_{2},s) \ar[r]^-{a\times b} & \dots }$$
$$\xymatrix{\dots \ar[r] & \pi_{1}^{pr}(F_{1}\times^{h}_{F_{0}}F_{2},s) \ar[r]^-{a\times b} &
\pi_{1}^{pr}(F_{1},s_{1})\times \pi_{1}^{pr}(F_{2},s_{2}) \ar[r]^-{a-b} &}$$ 
$$\xymatrix{\pi_{1}^{pr}(F_{0},s_{0}) \ar[r]^-{\partial} &  \pi_{0}^{pr}(F_{1}\times^{h}_{F_{0}}F_{2}) \ar[r]^-{a\times b} & \pi_{0}^{pr}(F_{1})\times \pi_{0}^{pr}(F_{2}) \ar[r]^-{a-b} & \pi_{0}^{pr}(F_{0}).}$$
Par convention, une suite de morphismes d'ensembles $\xymatrix{M \ar[r]^{u} & N \ar[r]^{a-b} & P}$
est exacte lorsque l'image de $u$ est exactement l'\'equaliseur des morphismes $a$ et $b$.
De m\^eme, dire qu'une suite de morphismes d'ensembles point\'es $\xymatrix{M \ar[r]^{u} & N \ar[r]^{v} & P}$
est exacte signifie que l'image de $u$ est \'egale \`a la fibre de $v$ au point distingu\'e.

En utilisant l'exactitude du foncteur de faisceautisation, on en d\'eduit une suite exacte longue sur les faisceaux d'homotopie
$$\xymatrix{
\pi_{m+1}(F_{0},s_{0}) \ar[r]^-{\partial} & \pi_{m}(F_{1}\times^{h}_{F_{0}}F_{2},s) \ar[r]^-{a\times b} &
\pi_{m}(F_{1},s_{1})\times \pi_{m}(F_{2},s_{2}) \ar[r]^-{a-b} & }$$ 
$$\xymatrix{ \pi_{m}(F_{0},s_{0}) \ar[r]^-{\partial} & \pi_{m-1}(F_{1}\times^{h}_{F_{0}}F_{2},s) \ar[r]^-{a\times b} & \dots }$$
$$\xymatrix{\dots \ar[r] & \pi_{1}(F_{1}\times^{h}_{F_{0}}F_{2},s) \ar[r]^-{a\times b} &
\pi_{1}(F_{1},s_{1})\times \pi_{1}(F_{2},s_{2}) \ar[r]^-{a-b} &}$$ 
$$\xymatrix{\pi_{1}(F_{0},s_{0}) \ar[r]^-{\partial} &  \pi_{0}(F_{1}\times^{h}_{F_{0}}F_{2}) \ar[r]^-{a\times b} & 
\pi_{0}(F_{1})\times \pi_{0}(F_{2}) \ar[r]^-{a-b} & \pi_{0}(F_{0}).}$$

L'existence de cette suite exacte longue et un analogue du lemme des cinqs impliquent en particulier que l'objet
$F_{1}\times^{h}_{F_{0}}F_{2}$ ne d\'epend  pas de la topologie de $C$. En d'autres termes, si $F \in SPr(C)^{I}$ est un objet, et si on
note $F_{1}\times^{hf}_{F_{0}}F_{2}$ le produit fibr\'e homotopique lorsque $C$ est muni de la topologie triviale, alors le morphisme
naturel dans $Ho(SPr(C))$, $F_{1}\times^{hf}_{F_{0}}F_{2} \longrightarrow F_{1}\times^{h}_{F_{0}}F_{2}$ est un isomorphisme. Une cons\'equence
remarquable de ce fait est que la cat\'egorie de mod\`eles $SPr(C)$ est propre \`a droite (i.e. les changements de bases le long
de fibrations pr\'eservent les \'equivalences). Pour plus de r\'ef\'erences sur ce point le lecteur pourra consulter \cite[$\S 3$]{s1}. \\

Supposons maintenant que $I$ soit la cat\'egorie d\'efinie par l'ensemble ordonn\'e des entiers naturels. Ses objets
sont donc les entiers $n \in \mathbb{N}$, et il existe un unique morphisme non trivial de $n$ vers $m$ si et seulement si
$n > m$. La cat\'egorie $I$ peut donc se repr\'esenter sch\'ematiquement par
$$ \dots m \longrightarrow m-1 \longrightarrow \dots 1 \longrightarrow 0,$$
et un objet de $SPr(C)^{I}$ est donc la donn\'ee pour tout entier $n$ d'un pr\'efaisceau simplicial $F_{n} \in SPr(C)$, et de morphismes
$f_{n} : F_{n} \longrightarrow F_{n-1}$. La limite homotopique de $F \in SPr(C)^{I}$ sera alors not\'ee
$Holim_{n}F_{n}$.

Comme le foncteur de faisceautisation ne commute pas avec les limites infinies, on ne disposera pas en g\'en\'eral d'une suite exacte
de Milnor reliant les faisceaux d'homotopie de $Holim_{n}F_{n}$ avec ceux de $F$. Ceci entraine aussi qu'il n'existe pas
de d\'ecomposition de Postnikov en g\'en\'eral. Nous allons cependant utiliser le lemme \ref{l2} et la suite exacte de Milnor
de \cite[Prop. $VI.2.15$]{gj} pour donner un crit\`ere de convergence des d\'ecompositions de Postnikov. Pour cela nous supposerons connu
le cas absolu (tel qu'il est d\'ecrit dans \cite[$\S VI$]{gj} par exemple). \\

Soit $* \longrightarrow F$ un pr\'efaisceau simplicial point\'e tel que $\pi_{0}(F)\simeq *$. Quitte \`a se restreindre au sous-pr\'efaisceau
simplicial de $F$ des simplexes au-dessus du point $*$ (ce qui changera $F$ en un pr\'efaisceau simplicial \'equivalent),
on pourra m\^eme supposer que $F$ est $0$-r\'eduit, et donc que $\pi_{0}^{pr}(F)\simeq *$. De plus, quitte \`a prendre un remplacement fibrant pour
la structure forte on pourra supposer que $F(X)$ est un ensemble simplicial fibrant pour tout objet $X \in C$.
Rappelons que le pr\'efaisceau en groupes $\pi_{1}^{pr}(F,*)$ op\`ere naturellement sur $\pi_{n}^{pr}(F,*)$, et de m\^eme
le faisceau en groupes $\pi_{1}(F,*)$ op\`ere sur $\pi_{n}(F,*)$. Nous utiliserons alors les pr\'efaisceaux simpliciaux
$K(\pi_{1}^{pr}(F,*),\pi_{n}^{pr}(F,*),n+1)$  d\'efinis dans le paragraphe suivant. Ce sont des pr\'efaisceaux simpliciaux point\'es, connexes et
fibrants pour la structure forte, et avec
$$\pi_{i}^{pr}(K(\pi_{1}^{pr}(F,*),\pi_{n}^{pr}(F,*),n+1),*)\simeq * \; pour\;  i\neq 1,n+1$$
$$
\pi_{1}^{pr}(K(\pi_{1}^{pr}(F,*),\pi_{n}^{pr}(F,*),n+1),*)\simeq \pi_{1}^{pr}(F,*) $$
$$\pi_{n+1}^{pr}(K(\pi_{1}^{pr}(F,*),\pi_{n}^{pr}(F,*),n+1),*)\simeq \pi_{n}^{pr}(F,*),$$
et dont l'action du $\pi_{1}^{pr}$ sur le $\pi_{n}^{pr}$ est celle donn\'ee pr\'ec\'edemment.

Consid\'erons la tour de Postnikov de $F$, $\{\tau_{\leq n}F\}_{n}$, d\'efinie objet par objet de $C$ comme dans \cite[$VI.3$]{gj}.
Il s'agit d'un diagramme
$$\xymatrix{\dots \ar[r] & \tau_{\leq n}F \ar[r]^-{p_{n}} & \tau_{\leq n-1}F \ar[r]^-{p_{n-1}} & \dots \tau_{\leq 0}F=*},$$
o\`u la fibre homotopique du morphisme $p_{n}$ est \'equivalente pour la structure forte au pr\'efaisceau simplicial $K(\pi_{n}^{pr}(F,*),n)$.
On peut donc trouver pour $n>1$, un diagramme homotopiquement cart\'esien pour la structure forte
$$\xymatrix{
\tau_{\leq n}F \ar[r] \ar[d] & \tau_{\leq n-1}F \ar[d] \\
K(\pi_{1}^{pr}(F,*),1) \ar[r] & K(\pi_{1}^{pr}(F,*),\pi_{n}^{pr}(F,*),n+1),}$$
qui est un diagramme classifiant pour la fibration $\tau_{\leq n}F \longrightarrow \tau_{\leq n-1}F$.
Mais comme les produits fibr\'es homotopiques sont ind\'ependants
\`a \'equivalence pr\`es de la topologie sur $C$, on trouve aussi des diagrammes homotopiquement cart\'esiens
$$\xymatrix{
\tau_{\leq n}F \ar[r] \ar[d] & \tau_{\leq n-1}F \ar[d] \\
K(\pi_{1}(F,*),1) \ar[r] & K(\pi_{1}(F,*),\pi_{n}(F,*),n+1).}$$

La proposition suivante est due \`a R. Thomason dans le cas stable, et \`a
J.F. Jardine dans le cas g\'en\'eral (voir \cite[Lem. $3.4$]{j}).

\begin{prop}\label{p1}
Soit $* \longrightarrow F$ un pr\'efaisceau simplicial point\'e avec $\pi_{0}(F)\simeq *$, et
$E_{n}$ un mod\`ele fibrant pour $K(\pi_{n}(F,*),n+1)$.
Supposons qu'il existe un entier $N$, tel que pour tout entier $i\geq 0$ et tout $n>i+N$ on ait
$$\pi_{i}^{pr}(E_{n},*)\simeq 0.$$
Alors le morphisme naturel
$$F \longrightarrow Holim_{n}\tau_{\leq n}F$$
est un \'equivalence.
\end{prop}

\textit{Preuve:} Consid\'erons le diagramme commutatif
$$\xymatrix{
F \ar[rr]^-{f} \ar[rd]_{p} & & Holim_{n>0}\tau_{\leq n}F \ar[dl]^{q} \\
& \tau_{\leq 1}F & }$$
Comme $Holim_{n>0}\tau_{\leq n}F$ est naturellement \'equivalent \`a $Holim_{n}\tau_{\leq n}F$, il suffit de montrer que
le morphisme $f$ est une \'equivalence. De plus, une utilisation de la suite exacte longue en homotopie montre qu'il suffit
de d\'emontrer que le morphisme induit sur les fibres homotopiques des morphismes $p$ et $q$ est une \'equivalence.
Ceci implique que l'on peut remplacer $F$ par la fibre homotopique de $p$, et donc se ramener au cas o\`u $F$ est simplement connexe
(i.e. $\pi_{1}(F,*)\simeq \pi_{0}(F)\simeq *$).

On peut clairement supposer que $F$ est un champ, ce que nous ferons.
De plus, quitte \`a remplacer le diagramme $\{\tau_{\leq n}F_{n}\}_{n}$ par un diagramme qui lui est \'equivalent, on peut aussi supposer
que chaque $\tau_{\leq n}F$ est un champ, et que chacun des diagrammes
$$\xymatrix{
\tau_{\leq n}F \ar[r] \ar[d] & \tau_{\leq n-1}F \ar[d] \\
\bullet \ar[r] & E_{n},}$$
est homotopiquement cart\'esien pour la structure forte. D'apr\`es
\cite[Prop. $VI.2.15$]{gj}, il existe pour tout $i\geq 0$ d'une suite exacte courte de pr\'efaisceaux
$$\xymatrix{
Lim^{1}_{n}\pi_{i+1}^{pr}(\tau_{\leq n}F) \ar[r] &
\pi_{i}^{pr}(Holim^{triv}_{n}\tau_{\leq n}F) \ar[r] & Lim_{n}\pi_{i}^{pr}(\tau_{\leq n}F) \ar[r] & 0},$$
o\`u $Holim^{triv}$ d\'esigne la limite homotopique pour la topologie triviale (et o\`u nous ne
mentionons pas le point de base).
Cependant, l'hypoth\`ese faite sur les pr\'efaisceaux d'homotopie de $E_{n}$
entraine que le syst\`eme projectif $\{\pi_{i}^{pr}(\tau_{\leq n}F,*)\}_{n}$ est stationnaire.
La suite exacte pr\'ec\'edente se r\'eduit donc un isomorphisme de pr\'efaisceaux
$$\pi_{i}^{pr}(Holim^{triv}_{n}\tau_{\leq n}F,*)\simeq \pi_{i}^{pr}(\tau_{\leq n}F,*),$$
pour $n>i+N$. En passant aux faisceaux d'homotopie on en d\'eduit que le morphisme naturel
$Holim^{triv}_{n}\tau_{\leq n}F \longrightarrow \tau_{\leq n}F$ induit des isomorphismes de faisceaux
$$\pi_{i}(Holim^{triv}_{n}\tau_{\leq n}F)\simeq \pi_{i}(\tau_{\leq n}F,*),$$
pour $n>i+N$. Or, comme on a $\pi_{i}(F,*)\simeq \pi_{i}(\tau_{\leq n}F,*)$ d\`es que $i\leq n$,
ceci implique clairement que le morphisme $F \longrightarrow Holim_{n}^{triv}\tau_{\leq n}F$
est une \'equivalence. Le lemme \ref{l2} implique alors qu'il en est de m\^eme du morphisme naturel
$F \longrightarrow Holim_{n}\tau_{\leq n}F$.
\hfill $\Box$ \\

\textit{Remarque:} Remarquer que la condition $\pi_{i}^{pr}(E_{n})=0$ pour $n>i+N$ est \'equivalent \`a la condition
que $H^{n+1-i}(X,\pi_{n}(F,*))=0$,  pour tout $n>i+N$ et tout objet du site $X \in C$. \\

En corollaire de la d\'emonstration on trouve le r\'esultat suivant.
\begin{cor}\label{p1'}
Soit 
$$\xymatrix{F \ar[r] & \dots \ar[r] & F_{n} \ar[r] & F_{n-1} \ar[r] & \dots \ar[r] & \dots \ar[r] & F_{1} \ar[r] & F_{0}=*}$$
une
tour de champs point\'es et connexes, telle que pour tout $n>0$,
la fibre homotopique de $F_{n} \longrightarrow F_{n-1}$ soit \'equivalente \`a un pr\'efaisceau simplicial de la forme
$K(\pi_{n},n)$. Supposons de plus qu'il existe un entier $N>0$ tel que pour tout objet $X \in C$ et tout $n>0$ on ait
$H^{n+i-1}(X,\pi_{n})=0$ d\`es que $n>i+N$. Alors, on a
$$\pi_{0}(Holim_{n}F_{n})=* \qquad \pi_{i}(Holim_{n}F_{n},*)\simeq \pi_{i} \; pour \; i>0.$$
\end{cor}

\subsection{Cohomologie des pr\'efaisceaux simpliciaux}

Nous consid\`ererons essentiellement le cas particulier des pr\'efaisceaux simpliciaux point\'es $s : * \longrightarrow F$ tels que
le morphisme induit sur le faisceau des composantes connexes $* \longrightarrow \pi_{0}(F)$ soit un isomorphisme. Un tel
$F$ sera appel\'e un pr\'efaisceau simplicial point\'e et connexe.

\begin{df}\label{d2}
\emph{Un syst\`eme local $M$ sur un pr\'efaisceau simplicial point\'e et 
connexe $F$, est la donn\'ee d'un faisceau en groupes
ab\'eliens $M$, et d'une action de $\pi_{1}(F,s)$ sur $M$. 
Un morphisme entre deux tels syst\`emes locaux $M$ et $M'$, est la donn\'ee
d'un morphisme de faisceaux de groupes $M \longrightarrow M'$ qui 
commute avec les actions de $\pi_{1}(F,s)$.}

\emph{La cat\'egorie ab\'elienne des syst\`emes locaux sur $F$ ainsi d\'efinie sera not\'ee $SysLoc(F)$.}
\end{df}

Rappelons que pour un $\mathbb{V}$-groupe $G$ on peut construire son classifiant $K(G,1)$, qui est un $\mathbb{V}$-ensemble simplicial d\'efini de la fa\c{c}on suivante.
On commence par d\'efinir $E(G,1)$, qui est l'ensemble simplicial nerf du morphisme naturel $G \longrightarrow *$. Ainsi,
$E(G,1)_{m}:=G^{m+1}$, et les faces et d\'eg\'en\'erescences sont donn\'ees par les projections et les diagonales. Le groupe
$G$ op\`ere sur $E(G,1)$ en op\`erant sur lui-m\^eme \`a gauche. On d\'efinit alors $K(G,1):=E(G,1)/G$.
Cette construction s'\'etend naturellement au cas o\`u $G$ est un $\mathbb{V}$-ensemble simplicial muni d'une loi de groupe. En effet, on forme
alors l'ensemble bisimplicial $K(G,1)_{p,q}:=K(G_{p},1)_{q}$, et on d\'efinit $K(G,1)$ comme  son ensemble simplicial diagonal. Ainsi, on a
$K(G,1)_{n}=K(G_{n},1)_{n}$. On peut aussi commencer par former l'ensemble bisimplicial nerf du morphisme $G \longrightarrow *$, et
appeler sa diagonale $E(G,1)$. Le groupe $G$ op\`ere alors sur $E(G,1)$ et on a $K(G,1)=E(G,1)/G$. Rappelons aussi que
l'ensemble simplicial $K(G,1)$ est naturellement point\'e (par l'image de l'identit\'e de $E(G,1)$), et qu'il existe des isomorphismes
naturels $\pi_{m}(K(G,1),*)\simeq \pi_{m-1}(G,e)$.

Lorsque $G$ est un groupe ab\'elien simplicial l'ensemble simplicial $K(G,1)$ h\'erite d'une loi de groupe ab\'elien induite
par celle de $G$. On peut alors it\'erer la construction pr\'ec\'edente et poser
$$E(G,m):=E(K(G,m-1),1) \qquad K(G,m):=K(K(G,m-1),1).$$
La fonctorialit\'e des constructions pr\'ec\'edentes permet de d\'efinir pour tout pr\'e\-fais\-ceau en groupes simpliciaux $G$ sur $C$, des
pr\'efaisceaux simpliciaux $E(G,1)$ et $K(G,1)$. De plus, lorsque $G$ est ab\'elien on dispose aussi des pr\'efaisceaux simpliciaux
$E(G,m)$ et $K(G,m)$. Le groupe simplicial $K(G,m-1)$ op\`ere sur $E(G,m)$ et on a $K(G,m)=E(G,m)/K(G,m-1)$. Nous d\'efinirons par convention
$K(G,0)$ comme \'etant le faisceau de groupes simpliciaux $G$. \\

Soit $G$ un pr\'efaisceau en groupes op\`erant sur un pr\'efaisceau en groupes ab\'eliens $M$ \`a travers
un morphisme de pr\'efaisceaux $\rho : G \longrightarrow \underline{Aut}_{C}(M)$. Alors le groupe $G$ op\`ere encore
de fa\c{c}on naturelle sur les pr\'efaisceaux en groupes ab\'eliens simpliciaux $K(M,m-1)$, pour $m>0$. On peut donc former le produit semi-direct
de $G$ par $K(M,m-1)$, qui est un pr\'efaisceau en groupes simpliciaux not\'e $G\times_{\rho} K(M,m-1)$. On notera par la suite
$K(G,M,m):=K(G\times_{\rho} K(M,m-1),1)$. On dispose de plus d'une suite exacte
de pr\'efaisceaux en groupes simpliciaux
$$K(M,m-1) \longrightarrow G\times_{\rho}K(M,m) \longrightarrow G.$$
En en prenant l'image par le foncteur $K(-,1)$ on en d\'eduit deux nouveaux morphismes
$$\xymatrix{ K(M,m) \ar[r] & K(G,M,m) \ar[r]^-{p} & K(G,1).}$$
Pour $m=0$, on d\'efinit par convention $K(G,M,0)$ comme \'etant le pr\'efaisceau en groupes $G\times_{\rho}M$, produit semi-direct de $G$ par $M$.
La suite exacte longue en homotopie implique alors que $K(M,m)$ est naturellement \'equivalent dans $SPr(C)$ \`a la fibre homotopique
du morphisme $p$, et ceci pour tout $m\geq 0$. Le pr\'efaisceau simplicial $K(G,M,m)$ est donc caract\'eris\'e \`a \'equivalence forte
pr\`es par le fait que
$$\pi_{i}^{pr}(K(G,M,m),*)=* \; pour \; i\neq 1,m$$
$$\pi_{1}^{pr}(K(G,M,m),*)\simeq G \; 
\pi_{m}^{pr}(K(G,M,m),*)\simeq M,$$
et par l'action de $G$ sur $M$. \\

Soit maintenant $F$ un pr\'efaisceau simplicial connexe, et $\rho : G=\pi_{1}(F,s) \longrightarrow \underline{Aut}_{C}(M)$ un syst\`eme local.
On dispose d'un morphisme d'ensembles simpliciaux
$$p^{*} : \mathbb{R}\underline{Hom}(F,K(G,M,m)) \longrightarrow \mathbb{R}\underline{Hom}(F,K(G,1)).$$
De plus, il existe un morphisme naturel dans $Ho(SPr(C))$, $F \longrightarrow \tau_{\leq 1}F\simeq K(G,1)$, qui d\'efinit donc un point
$* \in \mathbb{R}\underline{Hom}(F,K(G,1))$. Nous d\'efinissons alors l'ensemble simplicial
$\mathbb{R}\underline{Hom}_{\rho}(F,K(M,m))$ comme \'etant la fibre homotopique du morphisme $p^{*}$ au point $*$. Remarquer
que si $\rho$ est le morphisme trivial alors l'ensemble simplicial $\mathbb{R}\underline{Hom}_{\rho}(F,K(M,m))$ est naturellement \'equivalent \`a $\mathbb{R}\underline{Hom}(F,K(M,m))$.

\begin{df}\label{d3}
\emph{Le $m$-\`eme groupe de cohomologie du pr\'efaisceau simplicial point\'e et connexe $F$, \`a coefficients dans le
syst\`eme local $\rho : \pi_{1}(F,s) \longrightarrow \underline{Aut}_{C}(M)$ est d\'efini par}
$$H^{m}(F,M):=\pi_{0}(\mathbb{R}\underline{Hom}_{\rho}(F,K(M,m))).$$

\emph{Plus g\'en\'eralement, si $F \in SPr(C)$ est un pr\'efaisceau simplicial (non point\'e) et $M$ un faisceau de groupes ab\'eliens sur $C$, nous
noterons
$$H^{m}(F,M):=\pi_{0}\mathbb{R}\underline{Hom}(F,K(M,m)).$$
Lorsque $F=h_{X}$, pour $X$ un objet du site $C$, ces groupes seront not\'es
$H^{m}(X,M)$, et $H^{m}(C,M)$ si $X$ est l'objet final de $C$.}
\end{df}

En utilisant les suites exactes longues en homotopie associ\'ees \`a des fibrations, il est facile de v\'erifier les deux faits suivants.

\begin{itemize}
\item Le groupe $H^{0}(F,M)$ s'identifie naturellement au sous-groupe de $\Gamma(M)$ invariant par
l'action de $\Gamma(\pi_{1}(F,s))$.
\item Pour toute suite exacte  $\xymatrix{0 \ar[r] & M \ar[r] & M' \ar[r] & M'' \ar[r] & 0}$ de syst\`emes locaux, il existe une
suite exacte longue fonctorielle
$$\xymatrix{0 \ar[r] & H^{0}(F,M) \ar[r] & H^{0}(F,M') \ar[r] & H^{0}(F,M'') \ar[r] & H^{1}(F,M) \ar[r]
& \dots }$$
$$\xymatrix{H^{m}(F,M) \ar[r] & H^{m}(F,M') \ar[r] & H^{m}(F,M'') \ar[r] &
H^{m+1}(F,M) \ar[r] & \dots }$$
\end{itemize}

Lorsque $(F,s)=(K(G,1),*)$, pour $G$ un pr\'efaisceau en groupes sur $C$, les foncteurs $H^{m}(F,M)$ sont en r\'ealit\'e des
foncteurs d\'eriv\'es. Comme nous ne connaissons pas de r\'ef\'erences g\'en\'erales sur le sujet nous donnerons une esquisse de preuve
du th\'eor\`eme suivant. Il s'agit de toute \'evidence d'un r\'esultat faisant partie du folklore de la th\'eorie.

\begin{thm}\label{t1}
Soit $G$ un pr\'efaisceau en groupes sur $C$ et $(F,s)=(K(G,1),*)$. Alors, pour tout entier $m$, le foncteur
$$\begin{array}{cccc}
H^{m}(F,-) : & SysLoc(F) & \longrightarrow & \mathbb{V}-Ab \\
& (\rho : G \rightarrow M) & \mapsto & H^{m}(F,M)
\end{array}$$
est isomorphe au $m$-\`eme foncteur d\'eriv\'e du foncteur $H^{0}(F,-)$.
\end{thm}

\textit{Esquisse de preuve:} Commen\c{c}ons par traiter le cas o\`u $G=\{1\}$ est trivial. Il faut alors montrer que pour tout faisceau de groupes
ab\'eliens $M$ sur $C$, $H^{m}(C,M)$ est isomorphe \`a la cohomologie de $C$, calcul\'ee par foncteur d\'eriv\'e, \`a coefficients dans $M$
(que nous noterons $H_{der}^{m}(C,-)$).
Pour cela, notons $C^{-}(C,Ab)$ la cat\'egorie des complexes de pr\'efaisceaux en $\mathbb{V}$-groupes ab\'eliens sur $C$ qui sont concentr\'es en degr\'es
n\'egatifs ou nuls (avec une diff\'erentielle de degr\'e $+1$). Pour chaque objet $E$ de cette cat\'egorie, on dispose de ses pr\'efaisceaux
de cohomologie $\underline{H}_{pr}^{n}(E)$, dont les faisceaux associ\'es seront not\'es $\underline{H}^{n}(E)$. Un quasi-isomorphisme
est par d\'efinition un morphisme de $C^{-}(C,Ab)$, $f : E \longrightarrow E'$,
tel que pour tout $m$, le morphisme induit $H^{m}(f) : \underline{H}^{m}(E) \longrightarrow \underline{H}^{m}(E')$ soit un isomorphisme
de faisceaux. On notera alors $D^{-}(C,Ab)$ la cat\'egorie obtenue \`a partir de $C^{-}(C,Ab)$ en inversant formellement
les quasi-isomorphismes. Il est alors bien connu que les foncteurs de cohomologie par foncteurs d\'eriv\'es se calculent par
la formule suivante
$$H^{m}_{der}(C,M)\simeq [\underline{\mathbb{Z}},M[m]]_{D^{-}(C,Ab)},$$
o\`u $\underline{\mathbb{Z}}$ est le pr\'efaisceau constant de fibre $\mathbb{Z}$, et $M[m]$ est le complexe concentr\'e en degr\'e
$-m$ et de valeurs $M$.

Introduisons $SAb(C)$, la cat\'egorie des pr\'efaisceaux en $\mathbb{V}$-groupes ab\'eliens simpliciaux sur $C$. D'apr\`es \cite[Cor. $III.2.3$]{gj}
il existe une \'equivalence de
cat\'egories $\Gamma : C^{-}(C,Ab) \longrightarrow SAb(C)$ (o\`u $\Gamma$ est d\'efini terme \`a terme au-dessus de $C$).
De plus, en rempla\c{c}ant les ensembles simpliciaux par des groupes ab\'eliens simpliciaux dans les paragraphes pr\'ec\'edents, on
munit $SAb(C)$ d'une structure de cat\'egorie de mod\`eles o\`u les \'equivalences et les fibrations sont detect\'ees dans
$SPr(C)$ (i.e. en oubliant la structure de groupes). A l'aide de \cite[Cor. $III.2.7$]{gj} on voit que
le foncteur $\Gamma$ est compatible avec la formation des faisceaux de cohomologie et d'homotopie, et qu'il d\'efinit une
bijection de l'ensemble des
quasi-isomorphismes de $C^{-}(C,Ab)$ vers l'ensemble des \'equivalences de $SAb(C)$.
Il induit donc une \'equivalence sur les cat\'egories homotopiques associ\'ees
$$\Gamma : D^{-}(C,Ab) \simeq Ho(SAb(C)).$$
Ainsi, pour tout pr\'efaisceau en groupes ab\'eliens $M$  on a
$$[\underline{\mathbb{Z}},M[m]]_{D^{-}(C,Ab)}\simeq [\Gamma(\underline{\mathbb{Z}}),\Gamma(M[m])]_{Ho(SAb(C))}.$$
Cependant $\Gamma(\underline{\mathbb{Z}})$ est \'equivalent dans $SAb(C)$ au pr\'efaisceau constant de fibre
$\mathbb{Z}$, et $\Gamma(M[m])$ est \'equivalent \`a $K(M,m)$. On a donc
$$H^{m}_{der}(C,M)\simeq [\underline{\mathbb{Z}},M[m]]_{Ho(C^{-}Ab(C))}\simeq [\underline{\mathbb{Z}},K(M,m)]_{Ho(SAb(C))}.$$
Pour achever la d\'emonstration dans ce premier cas, on fait appel \`a l'adjonction de Quillen
$$-\otimes \mathbb{Z} :  SPr(C)  \longrightarrow SAb(C)  \qquad  j :  SAb(C)  \longrightarrow  SPr(C),$$
o\`u $-\otimes \mathbb{Z}$ est le foncteur d'ab\'elianisation, et $j$ est le foncteur d'oublie de la structure de groupes.
Cette adjonction descend en une adjonction sur les cat\'egories homotopiques, et donne
$$
H^{m}_{der}(C,M)\simeq [\underline{\mathbb{Z}},K(M,m)]_{Ho(SAb(C))}\simeq [*,j(K(M,m))]_{Ho(SPr(C))}$$
$$\simeq \pi_{0}\mathbb{R}\underline{Hom}(*,K(M,m))\simeq H^{m}(C,M).$$

Revenons au cas g\'en\'eral o\`u $F=K(G,1)$, pour $G$ un pr\'efaisceau en groupes.
Rappelons que l'on peut d\'efinir un site
$G$-\'equivariant $C/BG$. Ses objets sont ceux de $C$, et ils seront not\'es symboliquement $X\longrightarrow BG$ pour
$X \in C$. La donn\'ee d'un morphisme dans $C/BG$
$$\xymatrix{
X \ar[rd] \ar[rr]^{u} &  & Y \ar[dl] \\
& BG & }$$
est la donn\'ee d'un couple $(f,x)$, o\`u $f : X \rightarrow Y$ est un morphisme de $C$, et $x \in G(X)$. Les morphismes se composent naturellement
par la formule $(g,y)\circ (f,x)=(g\circ f,f^{*}(y).x)$, pour $\xymatrix{X\ar[r]^{f} & Y \ar[r]^{g} & Z}$ et $x \in G(X)$, $y\in G(Y)$. La topologie
sur $C/BG$ est d\'efinie en d\'ecr\'etant que les familles couvrantes dans $C/BG$ sont celles qui sont couvrantes dans $C$.
Il est clair que $C/BG$ poss\`ede des produits fibr\'es, mais il ne poss\`ede pas de produits finis en g\'en\'eral (et encore moins d'objet final).
De plus, la donn\'ee d'un pr\'efaisceau (resp. d'un faisceau) sur $C/BG$ est naturellement \'equivalente \`a
celle d'un pr\'efaisceau (resp. d'un faisceau) sur  $C$ muni d'une action de $G$. Il existe donc une \'equivalence de cat\'egorie $SPr(C/BG)\simeq G-SPr(C)$, o\`u $G-SPr(C)$ est la cat\'egorie des objets de $SPr(C)$ munis d'une action de $G$. Ceci montre aussi que la cat\'egorie
ab\'elienne $SysLoc(K(G,1))$ est \'equivalente \`a la cat\'egorie des faisceaux en groupes ab\'eliens sur $C/BG$, et
donc poss\`ede suffisemment d'injectifs.

Rappelons aussi que pour tout objet $F \in SPr(C)$, on peut consid\'erer la cat\'egorie des objets sur $F$, $SPr(C)/F$, qui est naturellement munie
d'une structure de cat\'egorie de mod\`eles ferm\'ee simpliciale.
De plus, pour tout morphisme $f : F \longrightarrow F'$ on dispose d'une adjonction
de Quillen
$$\begin{array}{cccc}
f_{!} : & SPr(C)/F & \longrightarrow & SPr(C)/F' \\
& (G\rightarrow F) & \mapsto & (G\rightarrow F\rightarrow F') 
\end{array}$$
$$\begin{array}{cccc}
f^{*} : & SPr(C)/F' & \longrightarrow & SPr(C)/F \\
& (G\rightarrow F') & \mapsto & (F\times_{G}F' \rightarrow F),
\end{array}$$
o\`u $f_{!}$ est l'adjoint \`a gauche.
Comme $SPr(C)$ est une cat\'egorie de mod\`eles propre \`a droite, l'adjonction pr\'ec\'edente est en r\'ealit\'e
une \'equivalence de Quillen d\`es que $f$ est une \'equivalence.

Notons alors $\underline{Hom}_{F}$ les $Hom$ simpliciaux de $SPr(C)/F$. On v\'erifie \`a l'aide des d\'efinitions et de l'\'equivalence de Quillen
ci-dessus, que pour tout syst\`eme local
$\rho : G \longrightarrow \underline{Aut}_{C}(M)$ sur $K(G,1)$, on dispose d'isomorphismes naturels
$$H^{m}(K(G,1),M)\simeq \pi_{0}\mathbb{R}\underline{Hom}_{K(G,1)}(K(G,1),K(G,M,m)).$$
En utilisant un argument tout \`a fait similaire \`a celui utilis\'e dans \cite[$\S 2.3$]{to1} on peut d\'emontrer qu'il existe
une \'equivalence de Quillen de la cat\'egorie de mod\`eles $SPr(C)/K(G,1)$ vers la cat\'egorie de mod\`eles $SPr(C/BG)$
(on renvoit \`a \cite{kt} pour plus de d\'etails sur cette \'equivalence. Lorsque $G$ n'est
pas un pr\'efaisceau simplicial cofibrant on en prendra tout d'abord un remplacement cofibrant
en tant que pr\'efaisceau en groupes simpliciaux). Ainsi,
on en d\'eduit un isomorphisme naturel
$$H^{m}(K(G,1),M)\simeq \pi_{0}\mathbb{R}\underline{Hom}_{BG}(*,K(M,m)_{BG}),$$
o\`u $\underline{Hom}_{BG}$ d\'esigne le $Hom$ simplicial de la cat\'egorie de mod\`eles simpliciale $SPr(C/BG)$, et $K(M,m)_{BG}$ est le pr\'efaisceau simplicial sur $C/BG$ correspondant \`a travers l'\'equivalence $SPr(C/BG)\simeq G-SPr(C)$ au pr\'efaisceau simplicial $K(M,m)$
muni de son action naturelle. En faisant appel au cas particulier trait\'e pr\'ec\'edemment pour le site $C/BG$, on voit
que $H^{m}_{\rho}(K(G,1),M)$ est isomorphe \`a $H^{m}_{der}(C/BG,M)$. En particulier, si $M$ correspond \`a un objet injectif
dans la cat\'egorie des syst\`emes locaux sur $K(G,1)$, on a $H^{m}(K(G,1),M)=0$ pour $m>0$. On en d\'eduit que la collection de foncteurs
$\{H^{m}(K(G,1),-)\}_{m}$ forme un foncteur cohomologique effa\c{c}able, et qu'il est donc isomorphe au foncteur d\'eriv\'e
de $H^{0}(K(G,1),-)$. \hfill $\Box$ \\

Remarquer que le th\'eor\`eme pr\'ec\'edent appliqu\'e au site $C/X$, des objets sur un objet $X$, montre que le foncteur
$H^{m}(X,-)$ est le $m$-\`eme foncteur d\'eriv\'e du foncteur des sections globales $\Gamma : Ab(C/X) \longrightarrow Ab$. Ainsi, pour les pr\'efaisceaux
simpliciaux repr\'esentables par des objets du site, la cohomologie que nous avons d\'efinie dans la d\'efinition \ref{d3}
coincide avec  la cohomologie d\'efinie par foncteurs d\'eriv\'es.

\subsection{$H_{\infty}$-Champs}

Dans ce paragraphe nous avons rassembl\'e quelques r\'esultats de la th\'eorie du d\'ela\c{c}age dans le cadre des
pr\'efaisceaux simpliciaux. La r\'ef\'erence pour le cas topologique est \cite{se2}. \\

Soit $s : * \longrightarrow F$ un pr\'efaisceau simplicial point\'e. On d\'efinit le pr\'efaisceau
des lacets de base $s$ par la formule
$$\begin{array}{cccc}
\Omega_{s}F : & C^{o} & \longrightarrow & \mathbb{V}-SEns \\
& X & \mapsto & \Omega_{s(X)}F(X).
\end{array}$$
Ici, pour un ensemble simplicial point\'e $(X,x)$, $\Omega_{x}X$ est d\'efini comme \'etant l'ens\-em\-ble simplicial des morphismes
de $\Delta^{1}$ vers $X$ qui envoient les deux sommets de $\Delta^{1}$ sur $x$.
La d\'efinition pr\'ec\'edente fournit un foncteur $\Omega : SPr(C)_{*} \longrightarrow SPr(C)_{*}$ de la cat\'egorie des pr\'efaisceaux simpliciaux point\'es vers elle m\^eme. On peut v\'erifier que ce foncteur est en r\'ealit\'e de Quillen \`a droite. Il se d\'erive donc
en un foncteur sur les cat\'egories homotopiques
$$\mathbb{R}\Omega : Ho(SPr(C)_{*}) \longrightarrow Ho(SPr(C)_{*}).$$
A l'aide de la suite exacte longue en homotopie on v\'erifie ais\'ement que les faisceaux d'homotopie de $\mathbb{R}\Omega_{s}F$ sont donn\'es par
$$\pi_{m}(\mathbb{R}\Omega_{s}F,*)\simeq \pi_{m+1}(F,s) \qquad \pi_{0}(\mathbb{R}\Omega_{s}F)\simeq \pi_{1}(F,s).$$
En particulier, ceci montre que l'on peut calculer $\mathbb{R}\Omega_{*}F$ en supposant que $C$ est muni de la topologie triviale.
Pour \^etre plus pr\'ecis, ceci veut dire que si l'on note $\mathbb{R}_{f}\Omega_{*}$ le foncteur d\'eriv\'e \`a droite de $\Omega_{*}$ lorsque
$C$ est muni de la topologie triviale, alors le morphisme naturel dans $Ho(SPr(C))$
$$\mathbb{R}_{f}\Omega_{*}(F) \longrightarrow \mathbb{R}\Omega_{*}(F)$$
est un isomorphisme. Ceci va nous permettre de g\'en\'eraliser sans trop de probl\`emes les r\'esultats standards de \cite{se2} au cas
des pr\'efaisceaux simpliciaux. \\

\begin{df}\label{d4}
\emph{Un pr\'e-$\Delta^{o}$-pr\'efaisceau simplicial est un foncteur}
$$\begin{array}{cccc}
F : & \Delta^{o} & \longrightarrow & SPr(C)\\
& [m] & \mapsto & F_{m}
\end{array}$$
\emph{tel que $F_{0}=*$. La cat\'egorie des pr\'e-$\Delta^{o}$-pr\'efaisceaux simpliciaux sera not\'ee $Pr\Delta^{o}-SPr(C)$.}
\end{df}

Cette d\'efinition et les r\'esultats qui vont suivre valent en particulier pour le cas o\`u $C=*$ est le site trivial. Nous noterons
dans ce cas $Pr-\Delta^{o}-SEns:=Pr-\Delta^{o}-SPr(*)$, que nous appellerons la cat\'egorie des pr\'e-$\Delta^{o}$-ensembles
simpliciaux. \\

Comme $[0]$ est initial dans $\Delta^{o}$, on peut identifier la cat\'egorie $Pr\Delta^{o}-SPr(C)$ avec la cat\'egorie des pr\'efaisceaux
sur $\Delta-[0]$ \`a valeurs dans $SPr(C)_{*}$, la cat\'egorie des pr\'efaisceaux simpliciaux point\'es.
Comme $SPr(C)_{*}$ est une cat\'egorie de mod\`el\-es engendr\'ee par cofibrations on peut alors appliquer
\cite[$\S II.6.9$]{gj} et munir $Pr\Delta^{o}-SPr(C)$ d'une structure de cat\'egorie de mod\`eles.
Rappelons que les fibrations (resp. les \'equivalences)
pour cette structure sont les morphismes $f : G \longrightarrow G$ tels que
pour tout $[m] \in \Delta$, le morphisme induit $f_{m} : F_{m} \longrightarrow G_{m}$ soit une fibration (resp. une
\'equivalence) dans $SPr(C)$. \\

Soit $F : \Delta^{o} \longrightarrow D$ un objet simplicial d'une cat\'egorie $D$ poss\'edant des produits finis, et supposons que
$F_{0}$ soit l'objet final dans $D$. Rappelons que l'on peut alors d\'efinir pour tout $n\geq 1$ le $n$-\`eme morphisme de Segal
$F_{n} \longrightarrow F_{1}^{n}$.
Il est induit par les $n$-morphismes dans $\Delta$, $p_{i} : [1] \longrightarrow [n]$, o\`u chaque $p_{i}$ est d\'efini par
$p_{i}(0)=i$ et  $p_{i}(1)=i+1$.

\begin{df}\label{d5}
\emph{Un objet $F$ de $Pr\Delta^{o}-SPr(C)$ est appel\'e un $H_{\infty}$-champ 
s'il v\'erifie les trois conditions suivantes.}
\begin{itemize}
\item \emph{Pour tout $n$, $F_{n}$ est un champ.}
\item \emph{Pour tout $n$, le morphisme de Segal
$$F_{n} \longrightarrow F_{1}^{n}$$
est une \'equivalence dans $SPr(C)$.}
\item \emph{La loi de mono\"{\i}de naturellement induite sur le faisceau $\pi_{0}(F_{1})$ est une loi de groupe.}
\end{itemize}

\emph{Par extensions, nous appellerons encore $H_{\infty}$-champ tout objet de la cat\'egorie homotopique $Ho(Pr-\Delta^{o}-SPr(C))$ qui est
isomorphe \`a un $H_{\infty}$-champ.}
\end{df}

Pour un objet $F \in Pr\Delta^{o}-SPr(C)$ on d\'efinit son pr\'efaisceau simplicial classifiant $BF$, d\'efini par
$BF_{n}:=(F_{n})_{n}$. Si $F$ est vu comme un pr\'efaisceau bi-simplicial sur $C$, alors $BF$ en est le pr\'efaisceau simplicial diagonal.
De plus, comme $F_{0}=*$, $BF$ est un pr\'efaisceau simplicial r\'eduit (i.e. $F_{0}=*$), et en particulier
peut-\^etre vu comme un pr\'efaisceau simplicial point\'e, $BF \in SPr(C)_{*}$. Ceci d\'efinit \'evidemment un foncteur
$$B : Pr\Delta^{o}-SPr(C) \longrightarrow SPr(C)_{*}.$$
Ce foncteur poss\`ede un adjoint \`a droite, not\'e usuellement $\Omega_{*}$.
Pour m\'emoire rappelons que pour $F \in SPr(C)_{*}$, $\Omega_{*}(F)$ est d\'efini par
$$\begin{array}{cccc}
\Omega_{*}(F) : & \Delta^{o} & \longrightarrow & SPr(C) \\
& [m] & \mapsto & F^{\Delta_{*}^{n}}
\end{array}$$
o\`u $\Delta^{m}_{*}$ est  l'ensemble simplicial point\'e obtenu \`a partir de $\Delta^{m}$ en identifiant tous ses sommets en un point unique.
En d'autres termes, $\Omega_{*}F_{m}$ est le pr\'efaisceau simplicial des morphismes $\Delta^{m} \longrightarrow F$ qui envoient
tous les sommets de $\Delta^{m}$ sur le point distingu\'e de $F$. \\

Le th\'eor\`eme suivant est une g\'en\'eralisation du th\'eor\`eme principal de la th\'eorie du d\'ela\c{c}age de G. Segal.
Comme nous ne connaissons pas de r\'ef\'erences pr\'ecises nous donnerons quelques indications pour montrer qu'il se d\'eduit
du th\'eor\`eme original.

\begin{thm}\label{t2}
Le foncteur $\Omega_{*} : SPr(C)_{*} \longrightarrow Pr-\Delta^{o}-SPr(C)$ est de Quillen \`a droite, et son adjoint \`a gauche
$B$ pr\'eserve les \'equivalences. De plus les deux assertions suivantes sont vraies.

\begin{enumerate}
\item
Pour tout objet $F \in SPr(C)_{*}$, $\mathbb{R}\Omega_{*}(F)$ est un $H_{\infty}$-espace, et le morphisme d'adjonction
$$B\mathbb{R}\Omega_{*}F \longrightarrow F$$
est un isomorphisme dans $Ho(SPr(C)_{*})$ lorsque $F$ est connexe (i.e. $\pi_{0}(F)\simeq *$).

\item
De m\^eme, si $F$ est un $H_{\infty}$-champ, alors le morphisme d'adjonction
$$F \longrightarrow \mathbb{R}\Omega_{*}(BF)$$
est un isomorphisme dans la cat\'egorie homotopique $Ho(Pr-\Delta^{o}-SPr(k))$.
\end{enumerate}
\end{thm}

\textit{Esquisse de preuve:} Le fait que $(B,\Omega_{*})$ soit une adjonction de Quillen pour la structure forte est
imm\'ediat. De plus,
on v\'erifie sans probl\`emes que $\Omega_{*}$ envoie tout objet fibrant sur un objet fibrant. Le fait que le couple de foncteurs
$(B,\Omega_{*})$ est une adjonction de Quillen
est alors une propri\'et\'e g\'en\'erale des localisations de Bousfield \`a gauche (On pourra par exemple
appliquer \cite[Thm. $3.4.20$]{hi}).
Comme $B$ pr\'eserve les \'equivalences fortes et que les fibrations triviales
sont aussi des \'equivalences fortes, il est facile de voir que $B$ pr\'eserve aussi les \'equivalences. \\

$(1)$ Le fait que $\Omega_{*}(F)$ soit un $H_{\infty}$-champ lorsque $F$ est fibrant provient imm\'ediatement de sa
d\'efinition et du fait que $\pi_{0}(\Omega_{*}(F))\simeq \pi_{1}(F,*)$. Pour voir que le morphisme d'adjonction est un isomorphisme,
on peut par exemple utiliser le fait que le foncteur d\'eriv\'e de $\Omega_{*}$ est ind\'ependant de la topologie sur $C$. Quitte \`a se restreindre
au sous-pr\'efaisceau simplicial de $F$ des simplexes au-dessus du point $*$ on peut ainsi se
ramener au cas o\`u la topologie est  triviale. Dans ce cas c'est une application de \cite{se2} objet par objet de $C$. \\

$(2)$ Soit $F$ un $H_{\infty}$-champ, que l'on peut supposer fibrant dans $Pr-\Delta^{o}-SPr(C)$.
Comme les morphismes de Segal $F_{n} \longrightarrow F_{1}^{n}$ sont des \'equivalences entre
objets fibrants, ce sont des \'equivalences fortes. De plus, dire que $\pi_{0}(F)$ est un groupe est \'equivalent \`a dire que le morphisme
$F_{2} \longrightarrow F_{1}^{2}$, donn\'e par la "multiplication" et "la premi\`ere projection" est une \'equivalence. De m\^eme, c'est
une \'equivalence forte. On voit ainsi que $F$ est un $H_{\infty}$-champ lorsque $C$ est muni de la topologie triviale.
De plus, on a vu que le foncteur d\'eriv\'e $\mathbb{R}\Omega_{*}$ pouvait se calculer en utilisant la topologie triviale. On en d\'eduit donc
que l'on peut supposer que $C$ est muni de la topologie triviale. Le r\'esultat d\'ecoule alors d'une application de \cite{se2} objet par
objet de $C$. \hfill $\Box$ \\

Tout comme dans \cite{to1}, on d\'eduit du th\'eor\`eme \ref{t2}
une \'equivalence entre la th\'eorie homotopique des objets en
groupes dans $SPr(C)$ et celle des $H_{\infty}$-champs sur $C$.
Nous utiliserons en particulier que tout $H_{\infty}$-champ est
naturellement \'equivalent \`a un pr\'efaisceau en groupes simpliciaux.

\subsection{Sch\'emas en groupes affines}

Pour tout ce paragraphe $k$ sera un corps de caract\'eristique quelconque,
et nous noterons par $(Aff/k)_{fpqc}$ le site des $\mathbb{U}$-sch\'emas affines sur $Spec \, k$ muni
de la topologie fid\`element plate et quasi-compacte (voir \cite[$VI.6.3$]{sga3}). Nous noterons alors
$SPr(k)$ la cat\'egorie des pr\'efaisceaux en $\mathbb{V}$-ensembles simpliciaux sur le site 
de Grothendieck $(Aff/k)_{fpqc}$, munie de sa structure de cat\'egorie de mod\`eles
d\'ecrite dans le paragraphe $1.1$. \\

Notons $GpAff/k$ la cat\'egorie des groupes dans $Aff/k$, dont nous appellerons les objets des sch\'emas en groupes affines
(ainsi un sch\'ema en groupes affines sera toujours un \'el\'ement de $\mathbb{U}$). Tout
objet $G \in GpAff/k$ repr\'esente un faisceau en groupes $h_{G}$. En notant $Gp(k)$ la cat\'egorie des faisceaux en $\mathbb{V}$-groupes
sur $(Aff/k)_{fpqc}$,  ceci induit un foncteur
$$\begin{array}{cccc}
h : & GpAff/k & \longrightarrow & Gp(k) \\
& G & \mapsto & h_{G},
\end{array}$$
qui par le lemme de Yoneda est pleinement fid\`ele. Les objets dans l'image essentielle de ce foncteur seront aussi appel\'es
des sch\'emas en groupes affines.

Rappelons alors que le fait que $k$ soit un corps entraine que le foncteur $h_{G}$ est exact (voir \cite[$III, \S 3, n^{o} 7$]{dg}). Ceci signifie
que pour toute suite exacte dans $GpAff/k$
$$\xymatrix{1 \ar[r] & G_{1} \ar[r] & G_{2} \ar[r] & G_{3} \ar[r] & 1}$$
la suite induite
$$\xymatrix{1 \ar[r] & h_{G_{1}} \ar[r] & h_{G_{2}} \ar[r] & h_{G_{3}} \ar[r] & 1}$$
est une suite exacte de faisceaux sur $(Aff/k)_{fpqc}$. En particulier, la sous-cat\'egorie de $Gp(k)$ form\'ee des sch\'emas en groupes affines est stable
par noyaux, co-noyaux et images.
De plus, d'apr\`es \cite[$III$ Thm. $4.3$]{mi} on sait que si $G$ est un sch\'ema en groupe affine,
$X$ un sch\'ema affine et $P$ un $h_{G}$-torseur de faisceau quotient $h_{X}$, alors $P$ est repr\'esentable par un sch\'ema affine.
Ainsi, la sous-cat\'egorie
pleine de $Gp(k)$ form\'ee des sch\'emas en groupes affines est aussi stable par extensions. Enfin, comme le foncteur
d'oublie de la structure de groupe commute avec les limites, et qu'une $\mathbb{U}$-limite projective
de sch\'emas affines est un sch\'ema affine, on voit que la sous-cat\'egorie pleine de $Gp(k)$ des sch\'emas en groupes affines
est aussi stable par $\mathbb{U}$-limites projectives. Notons aussi que le faisceau image d'un morphisme de sch\'emas en groupes affines
$f : G \longrightarrow H$ est toujours repr\'esent\'e par un sous-sch\'ema ferm\'e dans $H$ (\cite[$III, \S 3, Thm. 7.2$]{dg}).

Rappelons enfin qu'un sch\'ema en groupes affine $G$ est la limite projective de ses quotients de types finis sur $k$. Ceci permet
alors d'identifier la cat\'egorie des sch\'emas en groupes affines avec celle des Pro-objets (dans $\mathbb{U}$) dans la cat\'egorie des
sch\'emas en groupes affines de type fini sur $k$ (voir \cite[$III, \S 3, Cor. 7.5$]{dg}).

Ces propri\'et\'es nous permettent d'identifier sans risque un sch\'ema en groupes affine $G$ au faisceau qu'il repr\'esente $h_{G}$. Ainsi, par la
suite nous noterons $G$ aussi bien pour d\'esigner l'un ou l'autre. \\

Nous appellerons repr\'esentation lin\'eaire d'un sch\'ema en groupes affine $G$, la donn\'ee d'un $G-\mathcal{O}$-module dont
le $\mathcal{O}$-module sous-jacent est quasi-coh\'erent et appartient \`a $\mathbb{U}$ (voir \cite[$II, \S 2$]{dg}). C'est donc la donn\'ee d'un $k$-espace vectoriel
$V$ (\'eventuellement de dimension infini) de $\mathbb{U}$, et pour tout sch\'ema affine $Spec\, A \in Aff/k$ d'un morphisme de groupes
$G(A) \longrightarrow Aut_{A}(V\otimes A)$, fonctoriel en $A$. La cat\'egorie des repr\'esentations lin\'eaires
de $G$ sera not\'ee $Rep_{k}(G)$. Pour $V$ une telle repr\'esentation, nous noterons encore par $V$ le faisceau de $G-\mathcal{O}$-modules
sous-jacent.
Notons aussi que toute repr\'esentation lin\'eaire est limite inductive de ses sous-repr\'esentations
de dimension finie sur $k$. Ainsi, la cat\'egorie $Rep_{k}(G)$ s'identifie \`a la cat\'egorie des Ind-objets de la cat\'egorie
des repr\'esentations lin\'eaires sur $G$ de dimension finie. Remarquons enfin qu'une repr\'esentation lin\'eaire $V$ de dimension finie $d$,
d'un sch\'ema en groupes affine $G$, d\'etermine une action de $G$ sur le sch\'ema en groupes affine $V\simeq \mathbb{G}_{a}^{d}$.

L'oublie de la structure de $\mathcal{O}$-module permet de voir une repr\'esentation lin\'eaire de $G$ comme un syst\`eme local sur
$K(G,1)$ (voir Def. \ref{d2}). Ceci fournit un foncteur $Rep_{k}(G) \longrightarrow SysLoc(K(G,1))$, qui est exact et fid\`ele, mais qui
n'est pas pleinement fid\`ele en g\'en\'eral. Ce foncteur nous permet de parler de cohomologie d'un sch\'ema en groupes affine $G$ \`a valeurs
dans une repr\'esentation lin\'eaire $\mathbb{V}$. Par d\'efinition il s'agira de la cohomologie de $K(G,1)$ \`a valeurs dans
le syst\`eme local correspondant. Ces groupes de cohomologie seront simplement not\'es $H^{i}(K(G,1),V)$, et sont naturellement
munis d'une structure de $k$-espaces vectoriels.
Remarquer que $H^{0}(K(G,1),V)$ s'identifie naturellement
au sous-$k$-espace vectoriel de $V$ des points fixes par $G$, qui sera not\'e par la suite $V^{G}$.

\begin{lem}\label{l4}
Pour tout entier $i$, le foncteur
$$H^{i}(K(G,1),-) : Rep_{k}(G) \longrightarrow k-Vect$$
est le $i$-\`eme foncteur d\'eriv\'e du foncteur $H^{0}(K(G,1),-)$.
\end{lem}

\textit{Esquisse de preuve:}   Soit $V \in Rep_{k}(G)$ une repr\'esentation lin\'eaire de $G$.
Le faisceau sous-jacent au syst\`eme local $V\in SysLoc(K(G,1))$ \'etant quasi-coh\'erent il est
acyclique sur $(Aff/k)_{fpqc}$.
Le pr\'efaisceau simplicial $K(V,i)$ est donc
fibrant en tant qu'objet de la cat\'egorie $G-SPr(k)$, des pr\'efaisceaux simpliciaux
munis d'une action de $G$ (ou encore des pr\'efaisceaux simpliciaux sur le site \'equivariant $(Aff/k)/BG$). La cohomologie
de $K(G,1)$ \`a valeurs dans $V$ peut donc se calculer \`a l'aide de la topologie triviale sur $Aff/k$. Mais alors, on sait
d'apr\`es le th\'eor\`eme \ref{t1} que le foncteur $H^{i}(K(G,1),-)$ est le $i$-\`eme foncteur d\'eriv\'e du foncteur des sections globales
sur la cat\'egorie
des pr\'efaisceaux en groupes ab\'eliens sur $Aff/k$, munis d'une action de $G$. Ainsi, d'apr\`es
la proposition \cite[$I$ Prop. $5.2.1$]{sga3} appliqu\'ee au pr\'efaisceau d'anneaux constant $\mathbb{Z}$,
on sait alors que $H^{i}(K(G,1),V)$ est la cohomologie de Hochschild de $G$ \`a coefficients dans $\mathbb{V}$. Une application
du th\'eor\`eme  \cite[$I$ Thm. $5.3.1$]{sga3} montre alors
que $H^{i}(K(G,1),V)$ se calcule aussi par foncteurs d\'eriv\'es dans la cat\'egorie $Rep_{k}(G)$. Ce qu'il fallait d\'emontrer.
\hfill $\Box$ \\

\begin{cor}\label{c'}
Pour toute repr\'esentation lin\'eaire $V$ d'un sch\'ema en groupes affine $G$, on a pour tout $i\geq 0$
$$H^{i}(K(G,1),V)\simeq Colim_{j}H^{i}(K(G,1),V_{j}),$$
o\`u $V_{j}$ parcourt les sous-repr\'esentations lin\'eaires de $V$ qui sont de dimension finie sur $k$.
\end{cor}

\textit{Preuve:} D'apr\`es le lemme pr\'ec\'edent et \cite[$I$ Prop. $5.2.1$]{sga3} on sait que les groupes de cohomologie
$H^{i}(K(G,1),V)$ sont isomorphes aux groupes de cohomologie de Hochschild de $G$ \`a coefficients dans $\mathbb{V}$. Mais il est
facile de voir que le complexe de Hochschild de $V$ est la limite inductive des complexes de Hochschild
de ses sous-repr\'esentations de dimension finie. \hfill $\Box$ \\

Rappelons qu'un sch\'ema en groupes affine $G$ est unipotent, si pour toute repr\'esentation lin\'eaire $V$ non-nulle on a
$V^{G}\neq 0$. D'apr\`es \cite[$IV, \S 2, Prop. 2.5$]{dg} un sch\'ema en groupe affine est unipotent si et seulement si pour tout
sous-sch\'ema en groupe ferm\'e $K \hookrightarrow G$ il existe un morphisme non nul de sch\'emas en groupes
$K \longrightarrow \mathbb{G}_{a}$. De plus, la sous-cat\'egorie pleine de $Gp(k)$ form\'ee des sch\'emas en groupes affines  unipotents
est stable par noyaux, co-noyaux, images, extensions et $\mathbb{U}$-limites projectives (voir \cite[$VI, \S 2, Prop. 2.3$]{dg}). Comme pour le cas des sch\'emas en
groupes affines, tout sch\'ema en groupes affine unipotent est la limite projective de ses quotients de type fini sur $k$. En particulier,
la cat\'egorie des sch\'emas en groupes affines unipotents est donc \'equivalente \`a la cat\'egorie des Pro-objets
dans la cat\'egorie des sch\'emas en groupes affines unipotents et de type fini sur $k$.  \\

Notons que tout sch\'ema en groupes affine $G$ poss\`ede un quotient unipotent maximal $G \longrightarrow G^{uni}$. Par d\'efinition,
$G^{uni}$ est la limite projective de tous les quotients de $G$ qui sont des sch\'emas en groupes affines unipotents et de
type fini sur $k$. On v\'erifie alors que le morphisme $G \longrightarrow G^{uni}$ est universel pour les morphismes vers
des sch\'emas en groupes affines et unipotents. Le foncteur induit sur les cat\'egories de repr\'esentations
$Rep_{k}(G^{uni}) \longrightarrow Rep_{k}(G)$ est exact et pleinement fid\`ele, et identifie $Rep_{k}(G^{uni})$ \`a la sous-cat\'egorie
pleine de $Rep_{k}(G)$ form\'ee des rep\'esentations unipotentes
(i.e. colimites filtrantes d'extensions succ\'essives de repr\'esentations triviales).

\begin{df}\label{d6}
\emph{Soit $G$ un $\mathbb{U}$-groupe, consid\'er\'e comme un pr\'efaisceau constant sur $(Aff/k)_{fpqc}$. Une compl\'etion affine (resp. unipotente) de $G$
est un sch\'ema en groupes affine $\mathcal{A}(G)$ (resp. un sch\'ema en groupes affine et unipotent $\mathcal{U}(G)$), muni
d'un morphisme de faisceaux en groupes $G \longrightarrow \mathcal{A}(G)$
(resp. $G \longrightarrow \mathcal{U}(G)$) dans $Gp(k)$, qui est universel
pour les morphismes
vers les sch\'emas en groupes affines (resp. les sch\'emas en groupes affines unipotents).}
\end{df}

Le lemme suivant est un fait tr\`es bien connu (voir \cite{hm}).

\begin{lem}\label{l5}
Tout $\mathbb{U}$-groupe $G$ admet une compl\'etion affine (resp. une compl\'etion unipotente).
\end{lem}

\textit{Preuve:} Traitons le cas des compl\'etions affines, celui des compl\'etions unipotentes et en tout point similaire.
Pour les besoins de la d\'emonstration, rappelons que pour
tout objet $X\in C$ dans une $\mathbb{V}$-cat\'egorie et toute sous-cat\'egorie pleine $\mathbb{V}$-petite $D\subset C$, on peut
former la cat\'egorie $\mathbb{V}$-petite $X/D$, des objets de $D$ sous $X$.
Ses objets sont les paires $(Y,u)$, o\`u $Y$ est un objet de $D$ et $u : X \longrightarrow Y$ est un morphisme dans $C$.
Un morphisme dans $X/D$ de $(Y,u)$ vers $(Z,v)$ est la donn\'ee d'un morphisme $f : Y \longrightarrow Z$ dans $D$ tel que
$f\circ u=v$.

Notons $GpAlg$ la sous-cat\'egorie pleine de $Gp(k)$ form\'ee des sch\'emas en grou\-pes affines et de type fini sur $k$.
Posons
$$G\longrightarrow \mathcal{A}(G):=Lim_{H\in G/GpAlg}H$$
o\`u la limite est prise dans la cat\'egorie des faisceaux en $\mathbb{V}$-groupes $Gp(k)$. Comme la cat\'egorie $GpAlg$
consiste en des sch\'emas de types finis sur $k$, elle poss\`ede une
sous-cat\'egorie pleine $\mathbb{U}$-petite qui lui est \'equivalente. Ainsi, $G$ \'etant un \'el\'ement de $\mathbb{U}$, la cat\'egorie $G/GpAlg$ est
en r\'ealit\'e \'equivalente \`a une cat\'egorie $\mathbb{U}$-petite.
Ceci implique que $\mathcal{A}(G)$ est isomorphe \`a une
$\mathbb{U}$-limite projective de sch\'emas en groupes affines et est donc lui-m\^eme un sch\'ema un groupes affine.
De plus, il est clair par construction que tous les morphismes vers un sch\'ema en groupes affine de type fini
$G \longrightarrow H$ se factorisent par le morphisme naturel
$G \longrightarrow \mathcal{A}(G)$.
Mais comme tout sch\'ema en groupes affine est une $\mathbb{U}$-limite projective de sch\'emas en groupes affines et de type fini, on
voit que $G \longrightarrow \mathcal{A}(G)$ poss\`ede la propri\'et\'e universelle requise. \hfill $\Box$ \\

Remarquons que d'apr\`es les propri\'et\'es universelles, pour tout $\mathbb{U}$-groupe $\Gamma$ on a un isomorphisme
$\mathcal{U}(\Gamma)\simeq \mathcal{A}(\Gamma)^{uni}$. \\

\section{Champs affines}

Nous pr\'esentons dans ce chapitre la premi\`ere des deux notions fondamentales de cet article, celle de champ affine.
Il s'agit d'une version de
nature homotopique de la notion de sch\'ema affine, qui est d\'efinie directement \`a partir d'une structure de
cat\'egorie de mod\`eles simpliciale
sur la cat\'egorie des $k$-alg\`ebres co-simpliciales. Cette cat\'egorie de mod\`eles sera alors utilis\'ee pour
d\'efinir le foncteur d\'eriv\'e du foncteur qui \`a une $k$-alg\`ebre $A$ associe le sch\'ema affine $Spec\, A$.
Nous montrerons que ce foncteur d\'eriv\'e induit un foncteur pleinement fid\`ele de la cat\'egorie
homotopique des $k$-alg\`ebres co-simpliciales dans la cat\'egorie homotopique des pr\'efaisceaux simpliciaux
sur $(Aff/k)_{fpqc}$. La cat\'egorie des champs affines sera alors d\'efinie comme la sous-cat\'egorie
image essentielle de ce foncteur. Ainsi, tout comme la cat\'egorie des sch\'emas affines est \'equivalente
\`a la cat\'egorie (oppos\'ee) des $k$-alg\`ebres, la cat\'egorie des champs affines sera \'equivalente \`a la
cat\'egorie homotopique (oppos\'ee) des $k$-alg\`ebres co-simpliciales.

Par la suite, nous d\'efinirons la notion d'affination d'un pr\'efaisceau simplicial, et nous en donnerons un
crit\`ere d'existence. En particulier, nous montrerons que tout ensemble simplicial appartenant \`a $\mathbb{U}$
poss\`ede une affination, dont nous donnerons une description conjecturale des faisceaux d'homotopie pour
des ensembles simpliciaux simplement connexes et de type fini.

Enfin, lorsque l'anneau de base est un corps, nous d\'emontrerons un crit\`ere pour caract\'eriser les champs
affines et connexes \`a l'aide de leurs faisceaux d'ho\-mo\-to\-pie. Nous en d\'eduirons
les analogues des th\'eor\`emes standards de l'homotopie rationnelle et $p$-adique, montrant ainsi que
la th\'eorie des champs affines est un cadre ad\'equat pour l'\'etude des mod\`eles alg\'ebriques
des types d'homotopie. \\

Pour tout ce chapitre $k$ sera un anneau commutatif unitaire appartenant
\`a $\mathbb{U}$, et $(Aff/k)_{fpqc}$ le site des $\mathbb{U}$-sch\'emas affines sur $Spec \, k$ muni
de la topologie fid\`element plate et quasi-compacte (voir \cite[$VI.6.3$]{sga3}). Nous noterons
$SPr(k)$ la cat\'egorie des pr\'efaisceaux en $\mathbb{V}$-ensembles simpliciaux sur $(Aff/k)_{fpqc}$, munie de sa structure de cat\'egorie de mod\`eles
d\'ecrite dans le chapitre pr\'ec\'edent, et $Ho(SPr(k))$ sa cat\'egorie homotopique. Nous verrons syst\'ematiquement
les $\mathbb{V}$-ensembles simpliciaux comme des pr\'efaisceaux simpliciaux constants sur le site $(Aff/k)_{fpqc}$. Ceci permet
de voir $\mathbb{V}-SEns$ somme une sous-cat\'egorie pleine de $SPr(k)$. Il faut cependant se m\'efier du fait que ce plongement
n'induit pas un plongement au niveau des cat\'egories homotopiques. \\

La cat\'egorie des $k$-alg\`ebres commutatives et unitaires de $\mathbb{V}$ sera not\'ee $k-Alg$. A tout objet $A \in k-Alg$ nous associerons
le pr\'efaisceau d'ensembles sur $Aff/k$, qui \`a un sch\'ema affine $S=Spec\, B$ fait correspondre
l'ensemble $Hom_{k-Alg}(A,B)$. Ce pr\'efaisceau sera not\'e $Spec\, A \in SPr(k)$. Remarquer que lorsque
$A$ est un \'el\'ement de $\mathbb{U}$,
alors $Spec\, A$ est repr\'esentable par le sch\'ema affine $Spec\, A \in Aff/k$.

Par convention, nous r\'eserverons l'expression \textit{sch\'ema affine} pour d\'esigner les objets de $Aff/k$,
ou encore les pr\'efaisceaux
qu'ils repr\'esentent. Ces seront donc toujours des $\mathbb{U}$-sch\'emas. Nous parlerons
aussi de \emph{champ repr\'esentable} pour d\'esigner tout objet de $Ho(SPr(k))$ isomorphe
\`a un sch\'ema affine (et donc de la forme $Spec\, A$ avec $A$ appartenant \`a $\mathbb{U}$).

\subsection{Rappel sur les alg\`ebres co-simpliciales}

Rappelons tout d'abord que la cat\'egorie $k-Alg$ des $\mathbb{V}$-alg\`ebres commutatives et unitaires sur $k$ poss\`ede
tout type de limites et de colimites param\'etr\'ees par des ensembles de $\mathbb{V}$. Les limites se calculent dans la cat\'egorie
des $k$-modules et les sommes et coproduits finis sont donn\'es par les produits tensoriels. Comme il est expliqu\'e dans
\cite[$II$ Ex. $2.8$]{gj}, la cat\'egorie $k-Alg^{\Delta}$, des objets co-simpliciaux dans $k-Alg$, est alors munie d'une structure de cat\'egorie
simpliciale. Rappelons que pour $X$ un $\mathbb{V}$-ensemble simplicial et $A \in k-Alg^{\Delta}$ l'objet $A^{X}$ est d\'efini
par
$$\begin{array}{cccc}
A^{X} : & \Delta & \longrightarrow & k-Alg \\
& [n] & \mapsto & \prod_{X_{n}}A_{n}.
\end{array}$$
On d\'efinit alors l'ensemble simplicial des morphismes entre deux objets $A,B \in k-Alg^{\Delta}$ par la formule usuelle
$$\begin{array}{cccc}
\underline{Hom}_{k-alg}(A,B) : & \Delta^{o} & \longrightarrow & \mathbb{V}-SEns \\
& [n] & \mapsto & Hom_{k-Alg^{\Delta}}(A,B^{\Delta^{n}}).
\end{array}$$
Enfin, pour $X \in \mathbb{V}-SEns$ et $A \in k-Alg^{\Delta}$ nous noterons $X\otimes A$ le produit externe de $A$ par $X$. On dispose des eternels isomorphismes d'adjonctions
$$\underline{Hom}(X,\underline{Hom}_{k-alg}(A,B))\simeq \underline{Hom}_{k-alg}(X\otimes A,B)\simeq
\underline{Hom}_{k-alg}(A,B^{X}),$$
pour $X \in \mathbb{V}-SEns$ et $A,B \in k-Alg^{\Delta}$. \\

Nous identifierons la cat\'egorie $k-Alg$ \`a la sous-cat\'egorie pleine des objets constants dans $k-Alg^{\Delta}$. Remarquons que
pour tout $\mathbb{V}$-ensemble simplicial $X$, $k^{X}$ est la $k$-alg\`ebre co-simpliciale de cohomologie de $k$. Explicitement on a
$$\begin{array}{cccc}
k^{X} : & \Delta & \longrightarrow & k-Alg \\
& [n] & \mapsto & Hom(X_{n},k).
\end{array}$$
Remarquer aussi que le foncteur $k^{(-)} : \mathbb{V}-SEns \longrightarrow (k-Alg^{\Delta})^{o}$ est l'adjoint \`a gauche du foncteur
$\underline{Hom}_{k-alg}(-,k) : k-Alg^{\Delta} \longrightarrow \mathbb{V}-SEns$. \\

Notons $C_{+}(k)$ la cat\'egorie des complexes de co-chaines concentr\'es en degr\'es positifs et appartenant \`a $\mathbb{V}$ ,
et $k-Mod^{\Delta}$ la cat\'egorie
des $k$-modules co-simpliciaux de $\mathbb{V}$. Il existe alors deux foncteurs inverses l'un de l'autres (c'est la correspondance duale
de Dold-Puppe, voir \cite[$1.4$]{k})
$$N : k-Mod^{\Delta} \longrightarrow C_{+} \qquad \Gamma : C_{+} \longrightarrow k-Mod^{\Delta},$$
o\`u pour $M \in k-Mod^{\Delta}$, $N(M)$ est le complexe des co-chaines normalis\'ees de $M$. Par d\'efinition, les $k$-modules
de cohomologie d'un objet $M \in k-Mod^{\Delta}$ seront d\'efinis par
$$H^{i}(M):=H^{i}(N(M)).$$
En transportant la structure de mod\`ele de $C_{+}$ \`a $k-Mod^{\Delta}$ on d\'eduit une
structure de cat\'egorie de mod\`eles sur les $k$-modules co-simpliciaux
pour laquelle les \'equivalences sont les quasi-isomorphismes (i.e. morphismes induisant
des isomorphismes sur tous les $H^{i}$), et les fibrations sont les \'epimorphismes.

Le produit tensoriel de $k$-modules s'\'etend de fa\c{c}on naturelle en un produit tensoriel sur $k-Mod^{\Delta}$. De plus, la cat\'egorie
$k-Mod^{\Delta}$ est aussi munie d'une structure simpliciale, qui fait de $k-Mod^{\Delta}$ une cat\'egorie en alg\`ebre
ferm\'ee sur la cat\'egorie mono\"{\i}dale ferm\'ee $\mathbb{V}-SEns$
(voir \cite[Def.$4.1.9$, Def. $4.1.13$]{ho}). Pour deux objets $M$ et
$N$ dans $k-Mod^{\Delta}$, on dispose de deux morphismes naturels, qui sont des quasi-isomorphismes inverses l'un de l'autre
$$N(M\otimes N) \longrightarrow N(M)\otimes N(N) \longrightarrow N(M\otimes N),$$
et qui sont de plus associatifs et unitaires (c'est le  th\'eor\`eme d'Eilenberg-Zilber, voir \cite[$\S IV$ Thm. $2.4$]{gj}). Ceci permet en particulier de montrer que les
foncteurs $N$ et $\Gamma$ pr\'eservent les homotopies. Ainsi, deux morphismes dans $k-Mod^{\Delta}$ sont homotopes
si et seulement si leur image dans $C_{+}$ le sont (c'est l'\'enonc\'e dual de \cite[Thm. $2.6$]{do}).\\

L'oubli de la loi d'alg\`ebre d\'efinit un foncteur $k-Alg^{\Delta} \longrightarrow k-Mod^{\Delta}$, qui \`a
une $k$-alg\`ebre co-simpliciale fait correspondre son $k$-module co-simplicial sous-jacent. Par souci de l\'eg\`eret\'e nous noterons
encore $A$ le $k$-module co-simplicial sous-jacent \`a une $k$-alg\`ebre co-simpliciale $A$. En particulier nous
utiliserons les notations $N(A)$, $H^{*}(A)$ \dots pour $A \in k-Alg^{\Delta}$. Notons au passage que $H^{*}(A)$ est toujours muni
d'une structure de $k$-alg\`ebre commutative gradu\'ee et unitaire induite par la structure de $k$-alg\`ebre sur $A$.

\begin{df}\label{d7}
\emph{Soit $f : A \longrightarrow B$ un morphisme de $k$-alg\`ebres co-simpliciales.}
\begin{itemize}
\item \emph{Nous dirons que $f$ est une \'equivalence si pour tout entier $i$ le morphisme induit
$$H^{i}(f) : H^{i}(A) \longrightarrow H^{i}(B)$$
est un isomorphisme de $k$-modules.}
\item \emph{Nous dirons que $f$ est une fibration si pour tout entier $n$ le morphisme induit
$$f_{n} : A_{n} \longrightarrow B_{n}$$
est un morphisme surjectif.}
\item \emph{Nous dirons que $f$ est une cofibration s'il poss\`ede la propri\'et\'e de rel\`evement \`a gauche par rapport \`a toutes les
fibrations triviales (voir \cite[Def. $1.1.2$]{ho}).}

\end{itemize}
\end{df}

\begin{thm}\label{t3}
Avec les d\'efinitions pr\'ec\'edentes, la cat\'egorie $k-Alg^{\Delta}$ est une cat\'egorie de mod\`eles ferm\'ee
simpliciale.
\end{thm}

\textit{Esquisse de preuve:} C'est un fait standard qui se d\'emontre \`a l'aide de l'ar\-gu\-ment du petit objet. La preuve que nous
allons esquisser ici est une adaptation imm\'ediate de la preuve du th\'eor\`eme \cite[Thm. $4.3$]{bg}. Nous aurons besoin pour cela d'introduire les cofibrations et cofibrations triviales \'el\'ementaires suivantes.

Pour tout entier $n\geq 0$, notons $k[n]$ le complexe de co-chaine de $k$-modules \'egal \`a $k$ en degr\'e $n$ et $0$ partout
ailleurs. Nous noterons aussi $k<n-1,n>$ le complexe de co-chaine \'egal \`a $k$ en degr\'es $n-1$ et $n$, $0$ partout ailleurs
et avec comme diff\'erentielle $d^{n-1}=id : k \rightarrow k$. Par convention $k<-1,0>=0$.
Il existe des morphismes naturels de complexes
$k[n] \longrightarrow k<n-1,n>$. Nous noterons pour tout $n$
$$S^{n}_{k}:=\Gamma(k[n]) \qquad D_{k}^{n}:=\Gamma(k<n-1,n>).$$
Ce sont des $k$-modules co-simpliciaux munis de morphismes $S^{n}_{k} \longrightarrow D_{k}^{n}$. Les $k$-alg\`ebres
commutatives unitaires co-simpliciales libres engendr\'ees par ces $k$-modules co-simpliciaux seront not\'es
$$S(n):=S^{*}(S^{n}_{k}) \in k-Alg^{\Delta} \qquad D(n):=S^{*}(D_{k}^{n})\in k-Alg^{\Delta}.$$
Remarquer que $S(0)\simeq k[T]$, et $D(0)\simeq k$.

Nous d\'efinissons $I$ comme \'etant l'ensemble des morphismes naturels $S(n) \longrightarrow D(n)$ pour $n\geq 0$, et des morphismes
(correspondant \`a l'unit\'e) $k \longrightarrow S(n)$ pour $n\geq 0$. De m\^eme, l'ensemble $J$ est l'ensemble des morphismes
$k \longrightarrow D(n)$ pour $n\geq 0$. \\

On peut v\'erifier que pour tout $A \in k-Alg^{\Delta}$, la donn\'ee d'un morphisme $x : S(n) \longrightarrow A$ est
\'equivalente \`a la donn\'ee d'un cocycle dans le complexe normalis\'e associ\'e $x \in Z^{n}(N(A))$. De m\^eme,
la donn\'ee d'un morphisme $y : D(n) \longrightarrow A$ est \'equivalent \`a la donn\'ee d'un
\'el\'ement $y \in N(A)_{n-1}$. Ainsi, comme un morphisme de $k-Alg^{\Delta}$ est surjectif si et seulement si
le morphisme correspondant sur les complexes normalis\'es est surjectif, on v\'erifie \`a l'aide des propri\'et\'es universelles
pr\'ec\'edentes que les morphismes de $I$ sont des cofibrations, et que ceux de $J$ sont des cofibrations triviales. \\

Pour d\'emontrer que $k-Alg^{\Delta}$ est une cat\'egorie de mod\`eles ferm\'ee on utilise l'argument du petit objet relativement
aux ensembles $I$ et $J$ (voir \cite[$2.1.17$]{ho}). Il est alors bien connu qu'il suffit de montrer les trois assertions suivantes.

\begin{enumerate}
\item Un morphisme $f : A \longrightarrow B$ de $k-Alg^{\Delta}$ est une fibration (resp. une fibration triviale) si et seulement s'il
poss\`ede la propri\'et\'e de rel\`evement par rapport \`a  tous les morphismes de $J$ (resp. de $I$).

\item Une colimite filtrante d'\'equivalences dans $k-Alg^{\Delta}$ est une \'equivalence.

\item Pour tout morphisme $\tau_{n} : k \longrightarrow D(n)$ dans $J$, et tout objet $A \in k-Alg^{\Delta}$, le morphisme
naturel $A \longrightarrow A\otimes D(n)$ est une \'equivalence.

\end{enumerate}

La propri\'et\'e $(1)$ se d\'eduit imm\'ediatement des propri\'et\'es universelles des objets $S(n)$ et $D(n)$ et des d\'efinitions
des fibrations et des \'equivalences. La propri\'et\'e $(2)$ est vraie car les foncteurs $H^{i} : k-Alg^{\Delta} \longrightarrow k-Mod$
commutent avec les colimites filtrantes.

Pour d\'emontrer $(3)$ on commence par remarquer que les morphismes $0\longrightarrow k<n-1,n>$ sont des \'equivalences d'homotopie.
Comme le foncteur $\Gamma$ pr\'eserve les homotopies on trouve que $0 \longrightarrow D_{k}^{n}$ est aussi une \'equivalence par
homotopie dans $k-Mod^{\Delta}$. De plus, le foncteur qui \`a $M \in k-Mod^{\Delta}$ associe $S^{*}(M) \in k-Alg^{\Delta}$
est obtenu par prolongation du foncteur $S^{*} : k-Mod \longrightarrow k-Alg$, et d'apr\`es
l'\'enonc\'e dual de \cite[Thm. $5.6$]{do} il pr\'eserve donc les homotopies.
Ainsi, on trouve que les morphismes $k \longrightarrow D(n)$ sont des \'equivalences d'homotopie. Soit
$h : D(n) \longrightarrow D(n)^{\Delta^{1}}$ une contraction de $D(n)$ sur $k$. En tensorisant par $A$ on trouve un nouveau morphisme
$$h : A\otimes D(n) \longrightarrow A\otimes (D(n)^{\Delta^{1}})\simeq (A\otimes D(n))^{\Delta^{1}},$$
qui est un contraction de $A\otimes D(n)$ sur $A$. Ainsi, le morphisme $A \longrightarrow A\otimes D(n)$ est une \'equivalence
d'homotopie et donc une \'equivalence. \\

Il nous reste \`a d\'emontrer l'axiome $\mathbf{SM7}$ de \cite[$\S II$ Def. $3.1$]{gj}. Pour cela, il faut montrer que pour toute cofibration d'ensembles simpliciaux
$u : X \longrightarrow Y$, et toute cofibration $v : A \longrightarrow B$ dans $k-Alg$, le morphisme naturel
$$A^{X} \longrightarrow B^{X}\prod_{B^{Y}}A^{Y}$$
est une fibration, qui est une \'equivalence d\`es que $u$ ou $v$ en est une. Le fait que c'est une fibration est clair.
De plus, il est imm\'ediat de voir que
le morphisme $B^{X} \longrightarrow B^{Y}$ est une fibration, et donc que le morphisme de complexes $N(B^{X}) \longrightarrow
N(B^{Y})$ est surjectif. Comme $N$ commute avec les produits fibr\'es, on dispose ainsi d'une suite exacte longue de Mayer-Vietoris
$$\xymatrix{ 0 \ar[r] & H^{0}(B^{X}\prod_{B^{Y}}A^{Y}) \ar[r] &  H^{0}(B^{X})\prod H^{0}(A^{Y}) \ar[r] & H^{0}(B^{Y}) \ar[r] & \dots}$$
$$\xymatrix{ \dots \ar[r] & H^{i}(B^{X}\prod_{B^{Y}}A^{Y}) \ar[r] & H^{i}(B^{X})\prod H^{i}(A^{Y}) \ar[r] &
H^{i}(B^{Y}) \ar[r] & }$$
$$\xymatrix{ H^{i+1}(B^{X}\prod_{B^{Y}}A^{Y}) \ar[r] & \dots.}$$
On voit donc qu'il suffit de montrer que le morphisme $A^{X} \longrightarrow B^{X}$ (resp. $A^{X} \longrightarrow A^{Y}$)
est une \'equivalence lorsque $u$ (resp. $v$) en est une. Mais ceci se ram\`ene \`a l'\'enonc\'e connu que le
foncteur qui \`a un bi-complexe positif associe son complexe total pr\'eserve les quasi-isomorphismes. \hfill $\Box$ \\

\textit{Remarque:} Remarquer que tous les objets de $k-Alg^{\Delta}$ sont fibrants. De plus, si $A$ est une $k$-alg\`ebre co-simpliciale
constante, alors $A$ est un objet cofibrant dans $k-Alg^{\Delta}$. En effet, la donn\'ee d'un morphisme $A \longrightarrow B$ dans
$k-Alg^{\Delta}$ est \'equivalente \`a la donn\'ee d'un morphisme de $k$-alg\`ebres $A \longrightarrow H^{0}(B)$. Le morphisme
unit\'e $k \longrightarrow A$ poss\`ede donc la propri\'et\'e de rel\`evement par rapport \`a toutes les fibrations triviales. \\

Comme la cat\'egorie de mod\`eles $k-Alg^{\Delta}$ est simpliciale, on peut d\'efinir des $Hom$ simpliciaux en un sens d\'eriv\'es
(voir \cite[Thm. $4.3.2$]{ho}).  Ils seront not\'es $\mathbb{R}\underline{Hom}_{k-alg}$. Ainsi, pour $A$ et $B$ deux objets de $k-Alg^{\Delta}$,
$\mathbb{R}\underline{Hom}_{k-alg}(A,B)$ est un objet bien d\'efini dans $Ho(\mathbb{V}-SEns)$ et fonctoriel en $A$ et $B$. Cela munit
la cat\'egorie homotopique $Ho(k-Alg^{\Delta})$ d'une structure de cat\'egorie enrichie dans $Ho(\mathbb{V}-SEns)$.

Nous noterons aussi comme il en est l'usage $[-,-]_{k-alg}$ les ensembles de morphismes dans la cat\'egorie
homotopique $Ho(k-Alg^{\Delta})$. Remarquons au passage qu'il existe des isomorphismes naturels pour tout $A,B \in k-Alg^{\Delta}$
$$[A,B]_{k-alg}\simeq \pi_{0}\mathbb{R}\underline{Hom}_{k-alg}(A,B).$$

La structure simpliciale sur la cat\'egorie de mod\`eles $k-Alg^{\Delta}$ nous permettra par la suite de parler de limites et colimites
homotopiques, qui sont d\'efinies de fa\c{c}on standard (voir \cite[$\S 19$]{hi}). Nous utiliserons souvent que toute $k$-alg\`ebre co-simpliciale
$A$ est \'equivalente \`a la limite homotopique $Holim_{[n]\in \Delta}A_{n}$, o\`u chaque $A_{n}$ est vu comme une $k$-alg\`ebre co-simpliciale
constante (ce qui se v\'erifie imm\'e\-dia\-te\-ment en utilisant que tout ensemble simplicial $X$ est \'equivalent \`a la colimite
homotopique $Hocolim_{[n] \in \Delta^{o}}X_{n}$). \\

Pour cl\^ore ce paragraphe signalons que nous aurons a utiliser la sous-cat\'egorie pleine de $k-Alg^{\Delta}$ form\'ee des $k$-alg\`ebres
co-simpliciales appartenant \`a $\mathbb{U}$. Cette sous-cat\'egorie poss\`ede de tr\`es nombreuses propri\'et\'es de stabilit\'e. Nous utiliserons
en particulier que pour tout morphisme $f : A \longrightarrow B$ de $k$-alg\`ebres co-simpliciales de $\mathbb{U}$, les factorisations
fonctorielles d\'efinies \`a l'aide des ensembles g\'en\'erateurs $I$ et $J$ des cofibrations et des cofibrations triviales, appartiennent encore
\`a $\mathbb{U}$. En particulier, toute $k$-alg\`ebre co-simpliciale de $\mathbb{U}$ poss\`ede un mod\`ele cofibrant qui appartient encore \`a $\mathbb{U}$.

\subsection{Champs affines}

Pour $A \in k-Alg^{\Delta}$ une $k$-alg\`ebre co-simpliciale, on peut consid\'erer le pr\'efaisceau simplicial $Spec\, A \in SPr(k)$
d\'efini par
$$\begin{array}{cccc}
Spec\, A : & (Aff/k)^{o} & \longrightarrow & \mathbb{V}-SEns \\
& Spec\, B & \mapsto & \underline{Hom}_{k-alg}(A,B).
\end{array}$$
Le pr\'efaisceau des $n$-simplexes de $Spec\, A$ est donc isomorphe au pr\'efaisceau 
re\-pr\-\'es\-en\-t\'e par le sch\'ema affine $Spec\, A_{n}$. Ceci d\'efinit \'evidemment un
foncteur
$$Spec : (k-Alg^{\Delta})^{o} \longrightarrow SPr(k).$$
Ce foncteur poss\`ede un adjoint \`a gauche
$$\mathcal{O} : SPr(k) \longrightarrow (k-Alg^{\Delta})^{o}$$
qui est d\'efini de la fa\c{c}on suivante. Pour tout pr\'efaisceau simplicial $F \in SPr(k)$, on d\'efinit une $k$-alg\`ebre co-simpliciale
par
$$\begin{array}{cccc}
\mathcal{O}(F) : & \Delta & \longrightarrow & k-Alg \\
& [n] & \mapsto & Hom(F_{n},\mathcal{O}).
\end{array}$$
Dans cette d\'efinition $\mathcal{O}$ d\'esigne aussi le pr\'efaisceau en $k$-alg\`ebres tautologique sur $Aff/k$, qui au sch\'ema
affine $Spec\, A$ associe $A$. Le loi de $k$-alg\`ebre sur $Hom(F_{n},\mathcal{O})$ est alors induite naturellement par la loi
de $k$-alg\`ebre sur le pr\'efaisceau $\mathcal{O}$. Remarquer que $\mathcal{O}(F)$ n'est g\'en\'eralement pas un \'el\'ement de $\mathbb{U}$.

A l'aide de ces d\'efinitions on v\'erifie qu'il existe des isomorphismes naturels, d\'ecrivant le
comportement des foncteurs $Spec$ et $\mathcal{O}$ relativement \`a la structure simpliciale,
$$Spec\, (X\otimes A)\simeq Spec\, (A)^{X} \qquad \mathcal{O}(X\otimes F)\simeq \mathcal{O}(F)^{X},$$
pout tout $X \in \mathbb{V}-SEns$, $A \in k-Alg^{\Delta}$ et $F \in SPr(k)$. En particulier, l'adjonction
entre $Spec$ et $\mathcal{O}$ s'\'etend en une adjonction simpliciale
$$\underline{Hom}_{k-alg}(A,\mathcal{O}(F))\simeq \underline{Hom}(F,Spec\, A).$$

Le lemme suivant est une remarque cl\'e que nous utiliserons implicitement par la suite.

\begin{lem}\label{l6}
Si $A \in k-Alg^{\Delta}$ appartient \`a $\mathbb{U}$ alors le morphisme d'adjonction
$A \longrightarrow \mathcal{O}Spec\, A$ est un isomorphisme.
\end{lem}

\textit{Preuve:}  Le pr\'efaisceau $\mathcal{O}$ \'etant repr\'esent\'e par le sch\'ema affine $Spec\, k[T]$, c'est une application
du lemme de Yoneda. \hfill $\Box$ \\

Pour la proposition suivante, nous rappelons que $SPr(k)$ est munie da la structure de cat\'egorie de mod\`eles d\'ecrite dans le chapitre
pr\'ec\'edent, et que ses \'equivalences sont d\'efinies relativement \`a la topologie fid\`element plate et quasi-compacte. De m\^eme,
$Ho(SPr(k))$ d\'esignera la cat\'egorie homotopique des pr\'e\-fai\-sc\-eaux simpliciaux sur le site $(Aff/k)_{fpqc}$.

\begin{prop}\label{p2}
Le foncteur $Spec : (k-Alg^{\Delta})^{o} \longrightarrow SPr(k)$
est de Quillen \`a droite.
\end{prop}

\textit{Preuve:} L'axiome $\mathbf{SM7}$ des cat\'egories de mod\`eles simpliciales implique que le
foncteur $Spec$ est de Quillen \`a droite lorsque $SPr(k)$ est munie de sa structure forte. Par la propri\'et\'e universelle
des localisations de Bousfied \`a gauche il faut donc montrer que pour tout objet cofibrant $A \in k-Alg^{\Delta}$, le
pr\'efaisceau simplicial $Spec\, A$ est fibrant
(voir \cite[$3.4.20$]{hi}). Pour cela nous allons utiliser le crit\`ere \ref{l1}.

Soit $A \in k-Alg^{\Delta}$ un objet cofibrant, et $Spec\, B_{*} \longrightarrow Spec\, C$ un
hyper-recouvrement dans $(Aff/k)_{fpqc}$. Ce morphisme est donc donn\'e par un morphisme
d'alg\`ebres co-simpliciales
$$C \longrightarrow B_{*}.$$
Il s'agit alors
de montrer que le morphisme naturel
$$\underline{Hom}_{k-alg}(A,C) \longrightarrow Holim_{[n] \in \Delta}\underline{Hom}_{k-alg}(A,B_{n})$$
est une \'equivalence. Or, comme tous les objets de $k-Alg^{\Delta}$ sont fibrants et que $A$ est cofibrant, le morphisme
ci-dessus est isomorphe dans la cat\'egorie homotopique $Ho(\mathbb{V}-SEns)$ au morphisme
$$\mathbb{R}\underline{Hom}_{k-alg}(A,C) \longrightarrow Holim_{[n] \in \Delta}\mathbb{R}\underline{Hom}_{k-alg}(A,B_{n})
\simeq \mathbb{R}\underline{Hom}_{k-alg}(A,B_{*}).$$
En effet, en tant qu'objet de $k-Alg^{\Delta}$, $Holim_{[n] \in \Delta}B_{n}$ est naturellement \'equivalent \`a
$B_{*}$. Il suffit ainsi de montrer que le morphisme naturel
$$C \longrightarrow B_{*}$$
est une \'equivalence dans $k-Alg^{\Delta}$. Mais en passant aux complexes normalis\'es associ\'es on voit facilement que
cette derni\`ere assertion se d\'emontre en appliquant
la descente cohomologique fild\`element plate (voir \cite[Ex. $V^bis$]{sga42})
\`a l'hyper-recouvrement $Spec\, B_{*} \longrightarrow
Spec\, C$, et au complexe de faisceaux quasi-coh\'erents concentr\'e en degr\'e z\'ero $\mathcal{O}_{Spec\, C}$.
\hfill $\Box$ \\

\begin{cor}\label{c1}
Pour toute $k$-alg\`ebre co-simpliciale $A$ de $\mathbb{U}$, le morphisme d'ad\-jon\-ct\-ion
$$A\longrightarrow \mathbb{L}\mathcal{O}\mathbb{R}Spec\, A$$
est un isomorphisme dans $Ho(k-Alg^{\Delta})$.

En particulier, le foncteur d\'eriv\'e
$$\mathbb{R}Spec : Ho(k-Alg^{\Delta}) \longrightarrow Ho(SPr(k))$$
restreint aux $k$-alg\`ebres co-simpliciales isomorphes dans $Ho(k-Alg^{\Delta})$ \`a des objets de $\mathbb{U}$, est un foncteur
pleinement fid\`ele.
\end{cor}

\textit{Preuve:} Il suffit de montrer que pour tout $k$-alg\`ebre co-simpliciale $A$ de $\mathbb{U}$, qui est un objet cofibrant dans $k-Alg^{\Delta}$, le morphisme naturel $A\longrightarrow \mathbb{L}\mathcal{O}Spec\, A$ est un isomorphisme.
Consid\`erons l'isomorphisme naturel dans $Ho(SPr(k))$
$$Spec\, A \simeq Hocolim_{[n] \in \Delta^{o}}Spec\, A_{n}.$$
Remarquons que comme toute les $k$-alg\`ebres $A_{n}$ sont des \'el\'ements de $\mathbb{U}$, le pr\'e\-fai\-sc\-eau simplicial
repr\'esent\'e par $Spec\, A_{n}$ est un objet cofibrant dans $SPr(k)$. Ainsi, comme $\mathcal{O}$ est de Quillen
\`a gauche et que $A_{n}$ est un objet cofibrant dans $k-Alg^{\Delta}$, le lemme \ref{l6} implique que l'on a
$A_{n}\simeq \mathbb{L}\mathcal{O}(\mathbb{R}Spec\, A_{n})$.
Le morphisme d'adjonction s'\'ecrit donc
$$A \simeq Holim_{[n]\in \Delta}A_{n} \simeq Holim_{[n] \in \Delta} \mathbb{L}\mathcal{O}(\mathbb{R}Spec\, A_{n})$$
$$\simeq \mathbb{L}\mathcal{O}(Hocolim_{[n]\in \Delta}\mathbb{R}Spec\, A_{n})\simeq
\mathbb{L}\mathcal{O}(Spec\, A).$$
Il s'agit donc bien d'un isomorphisme. \hfill $\Box$ \\

\begin{df}\label{d8}
\emph{Un champ affine (sur $k$) est un champ $F \in SPr(k)$, qui vu comme un objet de $Ho(SPr(k))$ est isomorphe
\`a $\mathbb{R}Spec\, A$, pour $A$ une $k$-alg\`ebre co-simpliciale appartenant \`a $\mathbb{U}$.
Les morphismes de champs affines sont simplement les morphismes de champs.}

\emph{Par abus de langage, nous dirons qu'un pr\'efaisceau simplicial $F$ est un champ affine, si un mod\`ele fibrant de $F$ est
un champ affine.}
\end{df}

\textit{Remarques:} \begin{itemize}

\item Il est important de noter que $\mathbb{R}Spec\, A$ est par d\'efinition le pr\'efaisceau simplicial $Spec\, QA$, o\`u
$QA \longrightarrow A$ est un remplacement cofibrant de $A$ dans $k-Alg^{\Delta}$. Ainsi, un pr\'efaisceau de la forme
$Spec\, B$, o\`u $B$ est une $k$-alg\`ebre co-simpliciale de $\mathbb{U}$ quelconque n'est \`a prioris pas un champ affine. Comme le montre
le th\'eor\`eme \ref{t4} et le contre-exemple qui suit le corollaire \ref{c4} la notion de champ
affine est en r\'ealit\'e beaucoup plus forte que celle d'\^etre \'equivalent \`a un sch\'ema simplicial affine de $\mathbb{U}$.

\item
Noter qu'un pr\'efaisceau simplicial de la forme
$Spec\, A$, o\`u $A$ est une $k$-alg\`ebre co-simpliciale de $\mathbb{V}$ qui n'est pas \'equivalente dans $k-Alg^{\Delta}$
\`a un \'el\'ement de $\mathbb{U}$, n'est pas
consid\'er\'e comme un champ affine. En contre partie tous les pr\'efaisceaux repr\'esentables par des objets de $Aff/k$ le sont.
En effet, pour toute $k$-alg\`ebre $A$ de $\mathbb{U}$, la $k$-alg\`ebre co-simplciale constante de valeurs $A$ appartient \`a
$\mathbb{U}$ et
est un objet cofibrant. Ainsi, $\mathbb{R}Spec\, A$ est isomorphe dans $Ho(SPr(k))$ au pr\'efaisceau repr\'esent\'e par
le sch\'ema affine $Spec A$.

\item D'un point de vue terminologique il est important
de ne pas confondre \emph{champs affines} et \emph{sch\'emas
affines}. Tout sch\'ema affine est un champ affine, mais il existe
des champs affines et $0$-tronqu\'es (i.e. des faisceaux qui sont
affines en tant que champs) qui ne sont pas repr\'esentables par
des sch\'emas affines (un tel exemple est donn\'e \`a la suite du
corollaire \ref{c3}). Il est donc parfois pr\'ef\'erable
d'utiliser l'expression \emph{champs repr\'esentables} pour
signifier \emph{sch\'emas affines} et \'eviter ainsi des confusions.

\end{itemize}

Les premiers exemples de champs affines sont les champs classifiants $K(\mathbb{G}_{a},i)$, o\`u $\mathbb{G}_{a}$ est le
groupe additif. Ce sont en r\'ealit\'e les pi\`eces \'el\'ementaires de la th\'eorie, et nous verrons que
tout champ affine
est une certaine limite homotopique construite \`a partir de tels champs.

\begin{lem}\label{l7}
Notons $\mathbb{G}_{a}$ le pr\'efaisceau en groupes additifs sous-jacent au pr\'e\-fai\-sc\-eaux en $k$-alg\`ebres $\mathcal{O}$, et
$S(i)$ la $k$-alg\`ebre
co-simpliciale d\'efinie lors de la preuve du th\'eor\`eme \ref{t3}.
Alors, pour tout entier $i$, le pr\'efaisceau simplicial $K(\mathbb{G}_{a},i)\in SPr(k)$ est un champ affine, \'equivalent
pour la structure forte \`a $Spec\, S(i)$.
\end{lem}

\textit{Preuve:} A l'aide de la propri\'et\'e
universelle satisfaite par $S(i)$ on voit facilement que pour tout $Spec\, B \in Aff/k$,
$\pi_{m}^{pr}(Spec\, (S(i))(B))=0$ pour $m\neq i$, et $\pi_{i}^{pr}(Spec\, (S(i))(B))\simeq B$. De plus, le pr\'efaisceau des $0$-simplexes
de $Spec\, S(i)$ est isomorphe au pr\'efaisceau constant $*$. C'est alors un fait standard que tout pr\'efaisceau simplicial $F$, tel que
$F_{0}=*$, et $\pi_{m}^{pr}(F,*)=0$ pour $m\neq i$ est \'equivalent pour la structure forte sur $SPr(k)$ au
pr\'efaisceau $K(\pi^{pr}_{i}(F,*),i)$. On a donc une \'equivalence forte $Spec\, S(i)\simeq K(\mathbb{G}_{a},i)$.
\hfill $\Box$ \\

\begin{cor}\label{c2}
Pour tout pr\'efaisceau simplicial $F \in SPr(k)$ on a des isomorphismes fonctoriels
$$H^{i}(F,\mathbb{G}_{a})\simeq H^{i}(\mathbb{L}\mathcal{O}(F)).$$

En particulier, pour tout $k$-alg\`ebre co-simpliciale $A$ de $\mathbb{U}$, il existe des isomorphismes fonctoriels
$$H^{i}(A)\simeq H^{i}(\mathbb{R}Spec\, A,\mathbb{G}_{a}).$$
\end{cor}

\textit{Preuve:} Le membre de droite est par d\'efinition isomorphe \`a l'ensemble des morphismes
$[F,K(\mathbb{G}_{a},i)]_{SPr(k)}$. Ainsi, en utilisant l'adjonction
$(\mathcal{O},Spec)$ et la propri\'et\'e universelle satifaite par $S(i)$,  on trouve des isomorphismes naturels
$$H^{i}(F,\mathbb{G}_{a})\simeq [F,K(\mathbb{G}_{a},i)]_{SPr(k)}\simeq [S(i),\mathbb{L}\mathcal{O}(F)]_{k-alg}
\simeq H^{i}(\mathbb{L}\mathcal{O}(F)).$$
La seconde assertion se d\'eduit de la premi\`ere et du corollaire \ref{c1}. \hfill $\Box$ \\

\begin{prop}\label{p2'}
Les champs affines sont stables par $\mathbb{U}$-limites homotopiques.
\end{prop}

\textit{Preuve:} Comme toute $\mathbb{U}$-limite homotopique se d\'ecompose en une composition de produits homotopiques,
de produits fibr\'es homotopiques et de limites homotopiques suivant la cat\'egorie $\mathbb{N} : \dots \rightarrow  n \rightarrow n-1 \rightarrow \dots
\rightarrow 1$, il suffit de v\'erifier que les champs affines sont stables par ces trois types de limites homotopiques.

Si $I \in \mathbb{U}$, et $F_{i}\in Ho(SPr(k))$ sont des champs affines o\`u $i$ parcourt $I$, \
on peut \'ecrire $F_{i}\simeq Spec\, A_{i}$, avec $A_{i}$ une $k$-alg\`ebre co-simpliciale cofibrante appartenant \`a $\mathbb{U}$.
Dans ce cas, comme $Spec$ est un foncteur de Quillen \`a droite, on a
$$\mathbb{R}\prod_{i\in I}F_{i}\simeq Spec\, (\coprod_{i\in I}A_{i}),$$
ce qui montre que $\mathbb{R}\prod_{i\in I}F_{i}$ est bien un champ affine.

Soit $\xymatrix{F_{1} \ar[r]^-{u} & F_{0} & F_{2} \ar[l]_-{v}}$ un diagramme de champs affines dans $SPr(k)$. A \'equivalence
pr\`es, on peut d'apr\`es le corollaire \ref{c1} repr\'esenter ce diagramme comme l'image d'un diagramme de $k$-alg\`ebres
co-simpliciales cofibrantes et appartenant \`a $\mathbb{U}$,
$\xymatrix{A_{1} & \ar[l]_-{u}  A_{0} \ar[r]_-{v} & A_{2}}$, avec $u$ et $v$ des cofibrations.
On a alors
$$F_{1}\times^{h}_{F_{0}}F_{2}\simeq Spec\, (A_{1}\coprod_{A_{0}}A_{2}),$$
ce qui implique que le produit fibr\'e homotopique $F_{1}\times^{h}_{F_{0}}F_{2}$ est un champ affine.

Enfin, soit $\xymatrix{\dots F_{n} \ar[r]^-{u_{n}} & F_{n-1} \ar[r]^-{u_{n-1}} & \dots \ar[r]^-{u_{2}} & F_{1}}$ un
foncteur de $\mathbb{N}$ vers la cat\'egorie  $SPr(k)$ tel que chaque $F_{n}$ soit un champ affine. Par un argument de r\'ecurence
sur $n$, on voit que le peut remplacer ce diagramme par un diagramme \'equivalent qui est image par $Spec$ d'une
suite de cofibrations de $k$-alg\`ebres co-simpliciales cofibrantes et appartenant \`a $\mathbb{U}$
$$\xymatrix{\dots A_{n} & \ar[l]_-{u_{n}}  A_{n-1} & \ar[l]_-{u_{n-1}} \dots & \ar[l]_-{u_{2}}  A_{1}}.$$
Ainsi, $Spec$ \'etant un foncteur de Quillen \`a droite on trouve
$$Holim_{n}F_{n}\simeq Holim_{n}Spec\, A_{n} \simeq Spec\, (Colim_{n}A_{n}).$$
Ceci montre que $Holim_{n}F_{n}$ est un champ affine. \hfill $\Box$ \\

Nous sommes maintenant en mesure de donner une caract\'erisation des ch\-amps affines, qui \`a l'aide de la proposition p\'ec\'edente
permettra de construire de nombreux exemples. Pour cela
nous aurons besoin des d\'efinitions suivantes. Le lecteur pourra trouver des d\'etails sur la notion
d'objet locaux et de localisation dans \cite{hi}.

\begin{df}\label{d9}
\begin{itemize}

\item
\emph{Un pr\'efaisceau simplicial $F \in SPr(k)$ est sous-affine s'il existe un objet simplicial $X$ de $Aff/k$, tel que $F$ soit isomorphe
dans la cat\'egorie $Ho(SPr(k))$ \`a $h_{X}$, le pr\'efaisceau simplicial repr\'esent\'e par $X$.}

\item \emph{Un morphisme $f : G \longrightarrow H$ dans $SPr(k)$ est une $\mathcal{O}$-\'equivalence si pour tout $i\geq 0$ le
morphisme induit
$$f^{*} : H^{i}(H,\mathbb{G}_{a}) \longrightarrow H^{i}(G,\mathbb{G}_{a})$$
est un isomorphisme.}

\item \emph{Un pr\'efaisceau simplicial $F \in SPr(k)$ est $\mathcal{O}$-local si pour toute $\mathcal{O}$-\'equivalence $f : G
\longrightarrow H$, le morphisme induit
$$f^{*} : \mathbb{R}\underline{Hom}(H,F) \longrightarrow \mathbb{R}\underline{Hom}(G,F)$$
est un isomorphisme dans $Ho(\mathbb{V}-SEns)$.}

\end{itemize}
\end{df}

Par la suite, nous utiliserons tr\`es souvent le fait \'el\'ementaire qu'une $\mathbb{V}$-limite homotopique d'objets $\mathcal{O}$-locaux est
encore un objet $\mathcal{O}$-local.

\begin{thm}\label{t4}
Un champ $F$ est affine si et seulement s'il est sous-affine et $\mathcal{O}$-local.
\end{thm}

\textit{Preuve:}
Commen\c{c}ons par la necessit\'e.  Si $F$ est un champ affine on peut trouver une $k$-alg\`ebre co-simpliciale $A$ dans $\mathbb{U}$,
qui est de plus cofibrante dans $k-Alg^{\Delta}$, et tel que $F$ soit \'equivalent \`a $Spec\, A$. Or, comme chaque
$A_{n}$ est une $k$-alg\`ebre de $\mathbb{U}$, $Spec\, A$ est clairement isomorphe au pr\'efaisceau simplicial repr\'esent\'e par le
sch\'ema affine simplicial $[n] \mapsto Spec\, A_{n}$. Ceci montre que $F$ est sous-affine.

Montrons que $F$ est $\mathcal{O}$-local. Quitte \`a remplacer $A$ par un objet \'equivalent, on peut l'\'ecrire comme une composition
$$A^{(n)} \longrightarrow A^{(n+1)} \longrightarrow \dots Colim_{n}A^{(n)}=A,$$
o\`u pour chaque entier $n$ il existe un carr\'e cocart\'esien
$$\xymatrix{A^{(n)} \ar[r] & A^{(n+1)} \\
X=\coprod_{i\in I} X_{i} \ar[r] \ar[u] & Y=\coprod_{i \in I}Y_{i},\ar[u] }$$
tel que chaque $X_{i} \longrightarrow Y_{i}$ soit l'une des cofibrations \'el\'ementaires d\'ecrite dans la preuve du th\'eor\`eme \ref{t3}
(voir \cite[$2.1.14$]{ho}).
Ainsi, comme $Spec\, A$ est alors \'equivalent \`a la limite homotopique des $Spec\, A^{(n)}$ on se ram\`ene au cas o\`u $A=A^{(n)}$. De plus,
comme $Spec\, A^{(n)}$ est \'equivalent au produit fibr\'e homotopique de $Spec\, A^{(n-1)}$ et $Spec\, Y$ au-dessus
de $Spec\, X$
on peut supposer par r\'ecurrence que $A=X$ ou $A=Y$. Enfin, comme $Spec\, X$ (resp. $Spec\, Y$) est le produit homotopique
des $Spec\, X_{i}$ (resp. des $Spec\, Y_{i}$), on peut supposer que $A$ est l'un des domaines ou l'un des codomaines
des cofibrations \'el\'ementaires. Il nous suffit donc d'examiner le cas $A=S(p)$ et le cas $A=D(p)$, pour $p$ un entier positif.

Nous savons d'apr\`es le lemme \ref{l7} que $Spec\, S(p)$ est \'equivalent \`a $K(\mathbb{G}_{a},p)$, qui est
\'evidemment $\mathcal{O}$-local par d\'efinition. De m\^eme, on montre facilement \`a l'aide de la propri\'et\'e universelle
satisfaite par $D(p)$ que $\pi_{m}^{pr}(Spec\, D(p))=*$ pour tout entier $m$, et donc que $Spec \, D(p)$ est
\'equivalent au pr\'efaisceau simplicial $*$, qui est bien entendu $\mathcal{O}$-local. \\

Pour d\'emontrer la suffisance, soit $F$ un pr\'efaisceau simplicial $\mathcal{O}$-local, et consid\'erons le morphisme
d'adjonction dans $Ho(SPr(k))$
$$u : F \longrightarrow \mathbb{R}Spec\, \mathbb{L}\mathcal{O}(F).$$
Comme $F$ est sous-affine, l'objet $\mathbb{L}\mathcal{O}(F) \in Ho(k-Alg^{\Delta})$ est en r\'ealit\'e isomorphe \`a un objet
de $\mathbb{U}$. En effet, si $X$ est un sch\'ema affine simplicial repr\'esentant $F$, la $k$-alg\`ebre co-simpliciale
$\mathbb{L}\mathcal{O}(F)$ est \'equivalente \`a la $k$-alg\`ebre co-simpliciale $[n] \mapsto \mathcal{O}(X_{n})$. Ainsi, une application du corollaire \ref{c2} montre que le morphisme d'adjonction induit
des isomorphismes pour tout $i$
$$u^{*} : H^{i}(\mathbb{R}Spec\, \mathbb{L}\mathcal{O}(F),\mathbb{G}_{a})\simeq
H^{i}(\mathbb{L}\mathcal{O}(F))\simeq H^{i}(F,\mathbb{G}_{a}).$$
Ceci montre que $u$ est une $\mathcal{O}$-\'equivalence. Mais par hypoth\`ese $F$ est $\mathcal{O}$-local. De m\^eme, nous avons d\'ej\`a
d\'emontr\'e que $\mathbb{R}Spec\, \mathbb{L}\mathcal{O}(F)$ \'etait $\mathcal{O}$-local. Or, il est facile de v\'erifier \`a l'aide du lemme
de Yoneda qu'une $\mathcal{O}$-\'equivalence entre deux objets $\mathcal{O}$-locaux est une \'equivalence. Ainsi, $u$ est
une \'equivalence, et donc $F$ est un champ affine. \hfill $\Box$ \\

\textit{Remarques:} \begin{itemize}
\item La preuve du th\'eor\`eme pr\'ec\'edent montre en particulier que tout champ affine est construit par
limite homotopique \`a partir de champs de la forme $K(\mathbb{G}_{a},n)$. Ceci et la proposition \ref{p2'} montrent
que la cat\'egorie des champs affines est la plus petite sous-cat\'egorie pleine de $Ho(SPr(k))$ stables
par $\mathbb{U}$-limites homotopiques et contenant les champs $K(\mathbb{G}_{a},n)$. C'est la raison pour laquelle
les champs affines nous semblent tr\`es proches des \textit{complexes of unipotent bundles} discut\'es dans
\cite[p. $446$]{gr}.

\item
Le th\'eor\`eme \ref{t4} montre en particulier que la notion de champ affine est intrins\`eque \`a la th\'eorie
homotopique des pr\'efaisceaux simpliciaux. En particulier il s'agit d'une notion qui est ind\'ependante de la cat\'egorie
de mod\`eles $k-Alg^{\Delta}$ (qui n'a comme r\^ole que de donner un mod\`ele alg\'ebrique des champs affines). On pourrait tout
aussi bien utiliser d'autres cat\'egorie d'alg\`ebres (comme par exemple la cat\'egorie des $E_{\infty}$-alg\`ebres, ou encore
celle des alg\`ebres diff\'erentielles gradu\'ees commutatives lorsque $k$ est de caract\'eristique nulle), la th\'eorie
des champs affines en resterait inchang\'ee. Remarquons cependant que le corollaire \ref{c1} ne sera plus vrai pour
des $E_{\infty}$-alg\`ebre (sauf lorsque $k$ est de caract\'eristique nulle).
\end{itemize}

Le corollaire suivant est une extension du corollaire \ref{c1}. Il montre en particulier que tout pr\'efaisceau simplicial
de taille raisonable poss\`ede
une $\mathcal{O}$-localisation (i.e. un morphisme vers un objet $\mathcal{O}$-local qui est une $\mathcal{O}$-\'equivalence).
Nous reviendrons sur cette construction lors du pargraphe suivant (voir \ref{paff}).

\begin{cor}\label{c3}
Soit $F$ un pr\'efaisceau simplicial tel que $\mathbb{L}\mathcal{O}(F)$ soit isomorphe dans
$Ho(k-Alg^{\Delta})$ \`a un objet de $\mathbb{U}$. Alors,
le morphisme naturel $F \longrightarrow \mathbb{R}Spec\, \mathbb{L}\mathcal{O}(F)$
est universel dans la cat\'egorie homotopique $Ho(SPr(k))$ pour
les morphismes de $F$ vers des objets $\mathcal{O}$-locaux (et en particulier
vers des champs affines).

En particulier, si $F$ est un champ affine alors le morphisme d'adjonction
$F \longrightarrow \mathbb{R}Spec\, \mathbb{L}\mathcal{O}(F)$
est un isomorphisme dans $Ho(SPr(k))$.
\end{cor}

\textit{Preuve:} Le th\'eor\`eme \ref{t4} montre que $\mathbb{R}Spec\, \mathbb{L}\mathcal{O}(F)$ est $\mathcal{O}$-local, et
une application du corollaire \ref{c2} montre que le morphisme naturel $F \longrightarrow \mathbb{R}Spec\, \mathbb{L}\mathcal{O}(F)$ est une
$\mathcal{O}$-\'equivalence. La propri\'et\'e universelle des objets $\mathcal{O}$-locaux implique alors le r\'esultat.
\hfill $\Box$ \\

\textit{Exemples:}
\begin{itemize}
\item
Voici le contre-exemple typique d'un champ sous-affine $F$ qui n'est pas affine. Soit $D$ un $\mathbb{U}$-sch\'ema en groupes
diagonalisable sur $k$ (voir \cite[$Ex. I$ $4.4$]{sga3}), et $F=K(D,1)$, qui est clairement sous-affine. Alors, d'apr\`es le lemme \ref{l4} et
\cite[$Ex. I$ $5.3.3$]{sga3} on voit que
$$H^{i}(F,\mathbb{G}_{a})\simeq H^{i}_{H}(D,\mathbb{G}_{a})=0 \; pour \; i>0,$$
o\`u $H_{H}$ d\'esigne la cohomologie de Hochschild des sch\'emas en groupes affines.
En particulier, on a $\mathbb{R}Spec\, \mathbb{L}\mathcal{O}(F)\simeq *$, et donc d'apr\`es le corollaire \ref{c3} $F$ n'est certainement pas
un champ affine.

Cet exemple est assez repr\'esentatif de ce qui peut se passer, et en g\'en\'eral il est bon de retenir que les faisceaux
d'homotopie d'un champ affine auront tendence \`a \^etre des sch\'emas
en groupes affines et unipotents (voir Thm. \ref{t5'}).

\item Il existe des faisceaux d'ensembles qui sont des champs affines sans \^etre repr\'esentables
par des sch\'emas affines. Par exemple, le plan \'epoint\'e $\mathbb{A}^{2}-\{0\}$
est un champ affine. En effet, on peut \'ecrire $\mathbb{A}^{2}-\{0\}$ comme le champ
quotient $[Sl_{2}/\mathbb{G}_{a}]$, ou en d'autres termes il existe une suite exacte
de fibrations
$$\xymatrix{Sl_{2} \ar[r] & \mathbb{A}^{2}-\{0\} \ar[r] & K(\mathbb{G}_{a},1).}$$
On d\'eduit ais\'ement de ce diagramme et du lemme des cinqs que le morphisme
induit $\mathbb{A}^{2}-\{0\} \longrightarrow \mathbb{L}\mathcal{O}(\mathbb{A}^{2}-\{0\})$
est un isomorphisme de champs. Le corollaire \ref{c3} montre que
$\mathbb{A}^{2}-\{0\}$ est un champ affine.
\end{itemize}

\begin{cor}\label{c4}
Un pr\'efaisceau simplicial sous-affine $F$, qui est une $\mathbb{V}$-limite homotopique de
champs affines, est un champ affine.
\end{cor}

\textit{Preuve:} On invoque le th\'eor\`eme \ref{t4} et le fait qu'une
$\mathbb{V}$-limite homotopique de pr\'efaisceaux simpliciaux $\mathcal{O}$-locaux est encore $\mathcal{O}$-local.
\hfill $\Box$ \\

\subsection{Affination des types d'homotopie}

Dans le paragraphe pr\'ec\'edent nous avons d\'efini la notion de champs affines. Nous allons maintenant d\'efinir
la notion d'affination. Tout comme nous proposons les champs affines comme mod\`eles
aux \textit{complexes of unipotent bundles}, nous proposons la notion d'affination comme r\'eponse au
\textit{probl\`eme de la sch\'ematisation} de \cite{gr}. Nous r\'es\`erverons l'expression de \textit{sch\'ematisation} pour
la notion mieux adapt\'ee au cas non-simplement connexe discut\'ee dans le chapitre suivant.

\begin{df}\label{daff}
\emph{Soit $F \in Ho(SPr(k))$ un champ. Une affination de $F$, est la donn\'ee d'un champ affine $(F\otimes k)^{uni}$, et
d'un morphism dans $Ho(SPr(k))$, $u : F \longrightarrow (F\otimes k)^{uni}$ qui soit universel pour les morphismes vers des
champs affines.}
\end{df}

\textit{Remarque:} L'exposant \textit{uni} dans la notation $(F\otimes k)^{uni}$ fait r\'ef\'erence au mot
\textit{unipotent}. \\

Commen\c{c}ons par remarquer que lorsqu'une affination existe elle est unique. De plus, comme
les $K(\mathbb{G}_{a},n)$ sont des champs affines, le morphisme $u$ dans la d\'efinition pr\'ec\'edente
est toujours une $\mathcal{O}$-\'equivalence. Ainsi, le th\'eor\`eme \ref{t4} montre
qu'une affination est aussi une $\mathcal{O}$-localisation (i.e. une $\mathcal{O}$-equivalence vers un objet
$\mathcal{O}$-local), et le morphisme $u$ dans la d\'efinition pr\'ec\'edente est donc aussi universel pour
les morphismes vers des objects $\mathcal{O}$-locaux.

\begin{prop}\label{paff}
Soit $F \in Ho(SPr(k))$ un champ. Pour qu'une affination de $F$ existe il faut et il suffit
que $\mathbb{L}\mathcal{O}(F)$ soit isomorphe dans $Ho(k-Alg^{\Delta})$ \`a un objet de $\mathbb{U}$.

De plus, pour un tel pr\'efaisceau simplicial on a
$$(F\otimes k)^{uni}\simeq \mathbb{R}Spec\, (\mathbb{L}\mathcal{O}(F)).$$
\end{prop}

\textit{Preuve:} Il suffit de faire appel au corollaire \ref{c3}. \hfill $\Box$ \\

\begin{cor}\label{caff}
Tout ensemble simplicial $X$ appartenant \`a $\mathbb{U}$ (et vu comme pr\'e\-fai\-sc\-eau simplicial
constant sur $(Aff/k)_{fpqc}$) poss\`ede une affination $(X\otimes k)^{uni}$. De plus, pour un tel ensemble
simplicial $X$ on a
$$(X\otimes k)^{uni}\simeq \mathbb{R}Spec(k^{X}),$$
o\`u $k^{X}$ est la $k$-alg\`ebre co-simpliciale de cohomologie de $X$ \`a valeurs dans $k$.
\end{cor}

\textit{Preuve:} Le pr\'efaisceau simplicial constant associ\'e \`a un $\mathbb{U}$-ensemble simplicial $X$ \'etant un objet
cofibrant dans $SPr(k)$, on a
$$\mathbb{L}\mathcal{O}(X)\simeq \mathcal{O}(X)\simeq k^{X}.$$
Ceci montre que $\mathbb{L}\mathcal{O}(X)$ est bien isomorphe \`a un objet de $\mathbb{U}$, et termine
la preuve. \hfill $\Box$ \\

Notons $Ho(CAff/k)$  la sous-cat\'egorie pleine de $Ho(SPr(k))$ form\'ee
des champs affines sur $k$. En utilisant le corollaire \ref{c1}, on voit que le foncteur d\'eriv\'e des sections globales
$$\mathbb{R}\Gamma : Ho(SPr(k)) \longrightarrow Ho(\mathbb{V}-SEns)$$
envoie $Ho(CAff/k)$ dans la sous-cat\'egorie de $Ho(\mathbb{V}-SEns)$ form\'ee des objets isomorphes \`a des
$\mathbb{U}$-ensembles simpliciaux. Il se factorise donc en un foncteur
$$\mathbb{R}\Gamma : Ho(CAff/k) \longrightarrow Ho(\mathbb{U}-SEns).$$

\begin{cor}
Le foncteur d\'eriv\'e du foncteur des sections globales
$$\mathbb{R}\Gamma : Ho(CAff/k) \longrightarrow Ho(\mathbb{U}-SEns)$$
poss\`ede un adjoint \`a gauche
$$(-\otimes k)^{uni} : Ho(\mathbb{U}-SEns) \longrightarrow Ho(CAff/k).$$
\end{cor}

\textit{Preuve:} C'est une autre fa\c{c}on d'\'enoncer le corollaire \ref{caff}. \hfill $\Box$ \\

Soit $u : k \longrightarrow k'$ un morphisme d'anneaux dans $\mathbb{U}$. On peut alors d\'efinir
un foncteur de restriction $u^{*} : SPr(k) \longrightarrow SPr(k')$. Il est d\'efini par la formule suivante
$$u^{*}(F)(Spec\, A):=F(Spec\, A),$$
o\`u $F \in SPr(k)$, et $A \in Aff/k'$ qui est aussi consid\'er\'ee comme objet de $Aff/k$ via le morphisme
$\xymatrix{k \ar[r]^{u} & k' \ar[r] & A}$. Ce foncteur est clairement un foncteur de Quillen \`a droite qui de plus
pr\'eserve les \'equivalences. Il induit donc un foncteur de changement de bases
$$u^{*} : Ho(SPr(k)) \longrightarrow Ho(SPr(k')).$$
Nous utiliserons aussi la notation
$$F\otimes_{k}k':=u^{*}(F).$$

\begin{cor}\label{c12}
Soit $u : k \longrightarrow k'$ un morphisme de $\mathbb{U}$-anneaux, et $X$ un $\mathbb{U}$-ensemble simplicial
fini.
Alors, on a
$$(X\otimes k)^{uni}\otimes_{k}k'\simeq (X\otimes k')^{uni}.$$
\end{cor}

\textit{Preuve:} D'apr\`es la proposition \ref{paff} on a $(X\otimes k)^{uni}\simeq \mathbb{R}Spec\, (k^{X})$. Comme
$k^{X}$ est une $k$-alg\`ebre plate sur $k$, on a de plus
$$k^{X}\coprod^{\mathbb{L}}_{k}k'\simeq k^{X}\otimes^{\mathbb{L}}_{k}k'\simeq (k')^{X},$$
le second isomorphisme utilisant que $X$ est fini.
Ainsi, on a pour tout $A \in Aff/k'$
$$
(X\otimes k)^{uni}\otimes_{k}k'\simeq \mathbb{R}\underline{Hom}_{k'-Alg^{\Delta}}(k^{X},A)\simeq
\mathbb{R}\underline{Hom}_{k'-Alg^{\Delta}}(k^{X}\otimes^{\mathbb{L}}_{k}k',A)$$
$$\simeq
\mathbb{R}Spec\, ((k')^{X})(A)\simeq (X\otimes k')^{uni}(A). \hfill \Box$$

Le corollaire pr\'ec\'edent montre que l'affination d'un $\mathbb{U}$-ensemble simplcial fini
$X$ est stable par
changement de bases, et donc que le champ $(X\otimes \mathbb{Z})^{uni}$ permet de d\'eterminer les
champs $(X\otimes k)^{uni}$ pour tout anneau $k$, par la formule suivante
$$(X\otimes k)^{uni}\simeq (X\otimes \mathbb{Z})^{uni}\otimes_{\mathbb{Z}}k.$$
En g\'en\'eral, il est impossible de d\'ecrire explicitement le
champ $(X\otimes \mathbb{Z})^{uni}$, ni
m\^eme ses faisceaux d'homotopie. Par contre, lorsque $X$ est simplement connexe et de type fini (i.e. tous ses
groupes d'homotopie sont de type finis), on peut utiliser le th\'eor\`eme d'Eilenberg-Moore afin de
calculer les faisceaux $\pi_{i}(X\otimes \mathbb{Z})$. Le proc\'ed\'e standard permet alors de se ramener au cas
\'el\'ementaire o\`u $X=K(\mathbb{Z},n)$ (voir la preuve du th\'eor\`eme \ref{t8}).
Nous ne d\'emontrerons rien \`a ce sujet dans cet article, et nous
nous contenterons de donner une conjecture. \\

Pour toute $\mathbb{Z}$-alg\`ebre $A$ dans $\mathbb{U}$, notons $\Lambda(A)$ le groupe additif des vecteurs
de Witt universels \`a coefficients dans $A$ (voir \cite[$V$, \S $4$, $5.1.1$]{dg}). En d'autres termes
$\Lambda(A)$ est le groupe multiplicatif des s\'eries formelles \`a coefficients dans $A$ et de la forme
$1+a_{1}X+a_{2}X^{2}+ \dots + a_{n}X^{n}+\dots$. Le foncteur $Spec\, A \mapsto \Lambda(A)$ est clairement repr\'esentable
par un sch\'ema en groupes affine sur $Spec\, \mathbb{Z}$, qui sera encore not\'e $\Lambda$.

On d\'efinit un sous-foncteur $\mathbb{H}$ de $\Lambda$ de la fa\c{c}on suivante. Pour tout $A \in Aff/\mathbb{Z}$,
une s\'erie formelle $P(X) \in 1+XA[[X]]=\Lambda(A)$ appartient \`a $\mathbb{H}(A)$ si et seulement si elle
satisfait l'\'equation fonctionnelle $P(X).P(Y)=P(X+Y+XY)$. Le foncteur $\mathbb{H}$ est un sous-sch\'ema en groupes
ferm\'e dans $\Lambda$ et est donc aussi repr\'esentable par un sch\'ema an groupes affine
sur $Spec\, \mathbb{Z}$. Le sch\'ema en groupes $\mathbb{H}$ a \'et\'e introduit dans \cite[App. $C$]{js}
pour les besoins du calcul du compl\'et\'e pro-alg\'ebrique du groupe $\mathbb{Z}$, o\`u il est implicitement identifi\'e
au \textit{compl\'et\'e unipotent de $\mathbb{Z}$ sur $Spec\, \mathbb{Z}$}. On peut justifier
ce point de vue en observant que $\mathbb{H}$ est aussi le sch\'ema en groupes des automorphismes mono\"{\i}idaux du foncteur
d'oubli de la cat\'egorie des $\mathbb{Z}$-modules libres de rang fini munis d'un automorphisme unipotent.

Pour tout entier $n \in \mathbb{Z}$, on peut associ\'er la s\'erie formelle $(1+X)^{n} \in \mathbb{H}(\mathbb{Z})$. Ceci d\'efinit
un morphisme de pr\'efaisceaux en groupes $\mathbb{Z} \longrightarrow \mathbb{H}$. On en d\'eduit donc un morphisme de
pr\'efaisceaux simpliciaux $K(\mathbb{Z},n) \longrightarrow K(\mathbb{H},n)$. \\

\begin{conj}\label{conj2}
Le champ $K(\mathbb{H},n)$ est affine et le morphisme pr\'ec\'edemment d\'efini 
$K(\mathbb{Z},n) \longrightarrow K(\mathbb{H},n)$ est une affination.
\end{conj}

Une r\'eponse positive \`a la question pr\'ec\'edente impliquerait que pour tout $\mathbb{U}$-ensemble simplicial
$X$, qui est simplement connexe et de type fini, on ait
$$\pi_{i}((X\otimes \mathbb{Z})^{uni})\simeq \pi_{i}(X)\otimes_{\mathbb{Z}}\mathbb{H}.$$
En d'autres termes, les faisceaux d'homotopie de $X$ sont
les \textit{compl\'et\'es unipotents sur $Spec\, \mathbb{Z}$} des
groupes d'homotopie de $X$.

\subsection{Champs affines sur un corps}

Dans ce paragraphe nous allons \'etudier les champs affines et connexes lorsque $k$ est un corps. Il est remarquable que
l'on poss\`ede alors une caract\'erisation des champs affines \`a l'aide de leurs faisceaux d'homotopie. Ce crit\`ere
permet une grande lib\'ert\'e de manipulation, et il permettra en particulier
de d\'emontrer un crit\`ere analogue pour les types d'homotopie sch\'ematiques qui seront introduits
dans le chapitre suivant. Ce crit\`ere jouera par ailleurs un r\^ole essentiel dans la preuve
de l'existence du foncteur de sch\'ematisation. Enfin, les r\'esultats de ce paragraphe confortent aussi
le point de vue que les champs affines sont un mod\`ele pour les \textit{complexes of unipotent bundles}
de \cite{gr}. \\

Pour tout ce paragraphe $k$ sera un corps de caract\'eristique quelconque. \\

Nous commencerons par le th\'eor\`eme suivant. Pour la notion de sch\'emas en groupes unipotents on
renvoie le lecteur \`a \cite[$IV, \S 2$]{dg}, ou encore au bref r\'esum\'e se trouvant au paragraphe $1.5$.

\begin{thm}\label{t5}
Soit $* \longrightarrow F$ un pr\'efaisceau simplicial point\'e et connexe tel que pour tout entier $i>0$ le
faisceau en groupes $\pi_{i}(F,*)$ soit repr\'esentable par un sch\'ema en groupes affine et unipotent. Alors
$F$ est un champ affine.
\end{thm}

\textit{Preuve:} Pour r\'eduire le probl\`eme nous allons syst\'ematiquement appliquer la proposition \ref{p2'}. Le preuve
du th\'eor\`eme consiste donc \`a montrer que le champ $F$ s'\'ecrit comme une $\mathbb{U}$-limite homotopique
de champs de la forme $K(\mathbb{G}_{a},n)$.
Pour cela,
commen\c{c}ons par montrer que l'on est sous les hypoth\`eses de la proposition \ref{p1}
(et que l'on peut prendre $N=0$). Plus pr\'ecis\'emment nous allons montrer que pour tout entier $i$, et tout
sch\'ema affine $X$, on a $H^{p}(X,\pi_{i}(F,*))=0$ pour $p>1$ (il s'agit ici de cohomologie pour la topologie
$fpqc$).

Pour tout entier $n$, notons $\pi_{i}(F,*)_{n}$ le sous-groupe ferm\'e de $\pi_{i}(F,*)$ image du $n$-\`eme morphisme
de d\'ecalage (voir \cite[$IV$, $\S 3$ No. $4$]{dg})
(si $k$ est de caract\'eristique nulle on prendra $\pi_{i}(F,*)_{n}=0$ pour tout $n>0$). Comme
$\pi_{i}(F,*)$ est un sch\'ema en groupes unipotent l'intersection des $\pi_{i}(F,*)_{n}$ est nulle, et il existe
un isomorphisme de sch\'emas en groupes affines
$$\pi_{i}(F,*) \simeq Lim_{n}M_{n},$$
o\`u $M_{n}:=\pi_{i}(F,*)/\pi_{i}(F,*)_{n}$. On dispose alors d'une suite
exacte de Milnor
$$\xymatrix{0 \ar[r] & Lim_{n}^{1}H^{p-1}(X,M_{n}) \ar[r] & H^{p}(X,\pi_{i}(F)) \ar[r] & Lim_{n}H^{p}(X,M_{n}) \ar[r] & 0.}$$
Cette suite exacte montre en particulier qu'il suffit de montrer que pour tout $p>1$ et tout entier $n$ on a
$H^{p}(X,M_{n})=0$. En effet, dans ce cas les morphismes de transitions du syst\`eme projectif $\{H^{p-1}(X,M_{n})\}_{n}$
seront tous surjectifs pour $p>1$, et ainsi le terme $Lim_{n}^{1}H^{p-1}(X,M_{n})$ sera nul pour $p>1$.

On se ram\`ene ainsi au cas o\`u $\pi_{i}(F,*)_{n}=0$ pour un certain $n$. De plus, \`a l'aide de la suite exacte longue en cohomologie on
peut m\^eme
supposer que $\pi_{i}(F,*)_{1}=0$, ou en d'autres termes que $\pi_{i}(F,*)$ est annul\'e par le d\'ecalage. Dans ce cas, on sait
qu'il existe une suite exacte courte de sch\'emas en groupes (voir \cite[$IV$ $\S 3$ Thm. $6,6$]{dg} et \cite[$IV$ $\S 3$ Cor. $6.8$]{dg})
$$\xymatrix{0 \ar[r] & \pi_{i}(F,*) \ar[r] & \mathbb{G}_{a}^{I} \ar[r] & \mathbb{G}_{a}^{J} \ar[r] & 0.}$$
Or, comme chaque $\mathbb{G}_{a}$ est le faisceau en groupes sous-jacent \`a un $\mathcal{O}$-module quasi-coh\'erent, on a
pour $p>0$
$$H^{p}(X,\mathbb{G}_{a}^{I})\simeq H^{p}(X,\mathbb{G}_{a})^{I}\simeq 0,$$
et on voit \`a l'aide  de la suite exacte longue en cohomologie que ceci implique que $H^{p}(X,\pi_{i}(F,*))\simeq 0$ pour
$p>1$. \\

On vient de voir que les hypoth\`eses de \ref{p1} \'etaient satisfaites, et donc que $F$ est \'equivalent \`a la limite
homotopique de sa tour de Postnikov. Ainsi, $F$ est \'equivalent \`a la limite projective homotopique d'un diagramme
$\{F_{n}\}_{n}$, avec $F_{0}=*$, $F_{1}=K(\pi_{1}(F,*),1)$, et o\`u il existe des diagrammes homotopiquement cart\'esiens dans
$SPr(k)$
$$\xymatrix{
F_{n} \ar[r] \ar[d] & F_{n-1} \ar[d] \\
\bullet \ar[r] & K(\pi_{1}(F,*),\pi_{n}(F,*),n+1). }$$

Comme une $\mathbb{U}$-limite homotopique
de champs affines est un champ affine, la discussion pr\'ec\'edente montre que
l'on peut supposer que $F=K(G,M,n)$, pour $G$ un sch\'ema en groupes affine unipotent op\'erant sur un
sch\'ema en groupes affine unipotent ab\'elien $M$. Notons alors $M_{j}$ le sous-groupe de $M$ image du $j$-\`eme morphisme
de d\'ecalage. Comme les sous-groupes $M_{j}$ sont caract\'eristiques, on en d\'eduit une filtration
$G$-\'equivariante
$$\dots M_{j} \subset M_{j-1} \subset \dots M_{1} \subset M_{0}=M.$$
De plus, comme $M$ est un sch\'ema en groupes unipotent, le morphisme naturel $M \longrightarrow Lim_{j}M/M_{j}$ est un isomorphisme
$G$-\'equivariant.

\begin{lem}\label{lt1}
Le morphisme naturel
$$K(G,M,n) \longrightarrow Holim_{j}K(G,M/M_{j},n)$$
est une \'equivalence.
\end{lem}

\textit{Preuve:} Le diagramme commutatif suivant
$$\xymatrix{K(G,M,n) \ar[rr] \ar[dr]_{a} & & Holim_{j}K(G,M/M_{j},n) \ar[dl]^{b} \\
& K(G,1) & }$$
montre qu'il suffit de d\'emontrer que le morphisme induit sur les fibres homotopiques de $a$ et $b$ est
une \'equivalence. En d'autres termes on peut supposer que $G$ est trivial. Soit $E$ (resp. $E_{j}$) un mod\`ele
fibrant de $K(M,n)$ (resp. de $K(M/M_{j},n)$). Comme les groupes $M$ et $M/M_{j}$ sont des sch\'emas en groupes unipotents
on a d\'ej\`a vu que $\pi_{i}^{pr}(E)=\pi_{i}^{pr}(E_{j})=0$ pour $i\neq n,n-1$. On utilise alors la suite exacte courte de Milnor
$$\xymatrix{0 \ar[r] & Lim_{j}^{1}\pi_{i+1}^{pr}(E_{j}) \ar[r] & \pi_{i}^{pr}(Holim_{j}E_{j}) \ar[r] &
Lim_{j}\pi_{i}^{pr}(E_{j}) \ar[r] & 0.}$$
Cette suite exacte montre en particulier que $\pi_{i}(Holim_{j}E_{j})=0$ pour $i<n-2$. De plus, pour $i=n-1$, on trouve
$\pi_{n-2}^{pr}(Holim_{j}E_{j})\simeq Lim^{1}_{j}\pi_{n-1}^{pr}(E_{j})$. En d'autres termes, pour tout $X \in C$ on a
$$\pi_{n-2}(Holim_{j}E_{j}(X))\simeq Lim_{j}^{1}\pi_{n-1}(E_{j}(X))\simeq Lim_{j}^{1}H^{1}(X,M/M_{j}).$$
Or, comme $M/M_{j} \longrightarrow M/M_{j-1}$ est surjectif, et que $H^{2}(X,M/M_{j})=0$ pour tout $j$, les morphismes
de transitions $H^{1}(X,M/M_{j}) \longrightarrow H^{1}(X,M/M_{j-1})$ sont surjectifs. Ainsi, on a
$Lim_{j}^{1}H^{1}(X,M/M_{j})=0$. On en d\'eduit que $\pi_{n-2}(Holim_{j}E_{j})$ est trivial, et donc
que $Holim_{j}E_{j}$ est $(n-2)$-connexe. On peut donc invoquer le th\'eor\`eme \ref{t2} pour se ramener au cas o\`u
$n=1$. On conclut alors \`a l'aide du lemme \ref{l''} (on v\'erifiera qu'il n'y a pas d'argument circulaire).
\hfill $\Box$ \\

Le lemme pr\'ec\'edent permet de se ramener au cas o\`u $M_{j_{0}}$ est nul pour un certain $j_{0}$. On dispose dans ce cas de diagrammes
homotopiquement cart\'esiens dans $SPr(k)$
$$\xymatrix{K(G,M_{j-1},n) \ar[r] \ar[d] & K(G,M_{j-1}/M_{j},n) \ar[d] \\
K(G,1) \ar[r] & K(G,M_{j},n+1),}$$
qui montrent que l'on peut r\'eduire le probl\`eme par r\'ecurrence au cas o\`u $M$ est annul\'e par le morphisme de d\'ecalage.

Soit $M^{*}:=\underline{Hom}_{gp}(M,\mathbb{G}_{a})$ le $\mathcal{O}$-module quasi-coh\'erent des morphismes de sch\'emas en groupes de
$M$ vers $\mathbb{G}_{a}$. Le groupe $G$ op\`ere lin\'eairement sur $M^{*}$, ainsi que sur son $\mathcal{O}$-module
dual $V(M):=\underline{Hom}_{\mathcal{O}}(M^{*},\mathcal{O})$. De plus, comme $M$ est annul\'e par le d\'ecalage, il peut-\^etre
r\'ealis\'e comme un sous-groupe ferm\'e d'un produit de $\mathbb{G}_{a}$, et ceci implique que le morphisme de bidualit\'e
$$M \longrightarrow V(M)$$
est un monomorphisme $G$-\'equivariant de sch\'emas en groupes. Enfin, comme $V(M)$ est lui-m\^eme un sch\'ema en groupes
unipotent annul\'e par le d\'ecalage (car il est isomorphe \`a un produit de $\mathbb{G}_{a}$),
on peut poursuivre le proc\'ed\'e et construire une r\'esolution co-simpliciale
$G$-\'equivariante
$$M \longrightarrow V_{*}(M),$$
o\`u chaque $V_{n}(M)$ est le sch\'ema en groupes dual d'une repr\'esentation lin\'eaire de $G$.

\begin{lem}\label{lt2}
Le morphisme naturel
$$K(G,M,n) \longrightarrow Holim_{p \in \Delta}K(G,V_{p}(M),n)$$
est une \'equivalence dans $SPr(k)$.
\end{lem}

\textit{Preuve:} En consid\'erant le diagramme commutatif suivant
$$\xymatrix{K(G,M,n) \ar[rr] \ar[rd]_{a} & & Holim_{p\in \Delta}K(G,V_{p}(M),n) \ar[dl]^{b} \\
& K(G,1), & }$$
on voit qu'il suffit de montrer que le morphisme induit sur les fibres homotopiques de $a$ et $b$ est
une \'equivalence. En d'autres termes on peut supposer que $G$ est trivial.

Il est facile de voir que le morphisme $K(M,n) \longrightarrow Holim^{trvi}_{p\in \Delta}K(V_{p}(M),n)$ est une \'equivalence.
De plus, comme chaque $V_{p}(M)$ est un produit de faisceaux quasi-coh\'erents, il est acyclique, et donc
tous les pr\'efaisceaux simpliciaux $K(V_{p}(M),n)$ sont des champs. Le lemme se d\'eduit donc du lemme \ref{l2}.
\hfill $\Box$ \\

Le lemme pr\'ec\'edent permet de se ramener au cas o\`u $M$ est le $\mathcal{O}$-module dual d'une repr\'esentation
lin\'eaire de $G$. Comme une telle repr\'esentation est limite inductive de ses sous-repr\'esentations de dimension finie, on peut
\'ecrire $M$ sous la forme d'une limite projective de repr\'esentations lin\'eaires de dimension
finie
$$M\simeq Lim_{i}M_{i}.$$
Par un argument assez similaire \`a celui utilis\'e lors de la preuve du lemme pr\'ec\'edent on montre qu'il
existe une \'equivalence
$$K(G,M,n)\simeq Holim_{i}K(G,M_{i},n)$$
(on utilisera par exemple que le foncteur $Lim_{i}$ est exact dans la category des $k$-espaces vectoriels lin\'eairement
compacts). Ceci permet donc de se ramener au cas o\`u $M$ est une repr\'esentation lin\'eaire de dimension finie de $G$. Comme
le sch\'ema en groupes
$G$ est unipotent, il existe une suite de sous-repr\'esentations lin\'eaires $0=M_{r} \subset M_{r-1} \subset \dots \subset M_{1} \subset
M_{0}=M$, tel que $G$ op\`ere trivialement sur chaque quotient $M_{i}/M_{i+1}$.
Chaque champ $K(M_{i}/M_{i+1},n)$ est alors la fibre homotopique du morphisme naturel
$K(G,M/M_{i},n) \longrightarrow K(G,M/M_{i+1},n)$. Il existe donc des diagrammes homotopiquement cart\'esiens
dans $SPr(k)$
$$\xymatrix{K(G,M/M_{i+1},n) \ar[r] \ar[d] & K(G,M/M_{i},n) \ar[d] \\
K(G,1) \ar[r] & K(G,M_{i}/M_{i+1},n+1).}$$
On se ram\`ene ainsi par une r\'ecurrence au cas o\`u  $F=K(G,1)$ et au cas o\`u $F=K(G,M,n)$ et
$G$ op\`ere trivialement sur $M$. Dans ce cas, on a
$F=K(G,1)\times K(M,n)$ et il nous reste \`a montrer que chacun des facteurs est affine.

Pour le facteur $K(M,n)$, c'est imm\'ediat car $M$ est isomorphe \`a un produit fini de $\mathbb{G}_{a}$.
Enfin, pour le facteur $K(G,1)$ on
\'ecrit $G$ comme la limite projective de ses quotients de type fini sur $k$, ce qui permet
\`a l'aide du lemme suivant  de supposer que $G$ est  de type fini.

\begin{lem}\label{l''}
Soit $G$ un sch\'ema en groupe affine, et $\{G_{i}\}_{i\in I}$ le syst\`eme projectif de ses quotients
qui sont de type fini sur $k$. Alors, le morphisme naturel
$$K(G,1) \longrightarrow Holim_{i\in I}K(G_{i},1)$$
est une \'equivalence.
\end{lem}

\textit{Preuve:} Soit $E$ (resp. $E_{i}$) un mod\`ele fibrant pour $K(G,1)$ (resp. $K(G_{i},1)$). On peut prendre par exemple
pour $E$ (resp. pour $E_{i})$ le champ des $G$-torseurs (resp. des $G_{i}$-torseurs). Il s'agit alors de montrer que
pour tout sch\'ema affine $X$, le groupo\"{\i}de des $G$-torseurs sur $X$ est \'equivalent \`a la limite homotopique des
groupo\"{\i}des des $G_{i}$-torseurs sur $X$. Ce qui est vrai. \hfill $\Box$ \\

Lorsque $G$ est de type fini on utilise qu'il poss\`ede une suite de composition centrale dont les quotients succesifs sont
des sous-groupes ferm\'es de $\mathbb{G}_{a}$ (voir \cite[$IV, \S 2, Prop. 2.5$]{dg}). Il existe une suite
de composition centrale de sous-groupes ferm\'es $1=G_{r} \subset G_{r-1} \subset \dots \subset G_{0}=G$, et des
diagrammes homotopiquement cart\'esiens
$$\xymatrix{
K(G/G_{i+1},1) \ar[r] \ar[d] & K(G/G_{i},1) \ar[d] \\
\bullet \ar[r] & K(G_{i}/G_{i+1},2),}$$
tel que chaque $G_{i}/G_{i+1}$ soit un sous-groupe ferm\'e de $\mathbb{G}_{a}$.
Par r\'ecurrence il nous suffit de montrer que pour tout sous-groupe ferm\'e $G \subset \mathbb{G}_{a}$, le champ
$K(G,2)$ est affine.
Mais dans ce cas, ou bien $G=\mathbb{G}_{a}$, ou bien
il existe une suite exacte courte de sch\'emas en groupes
$$\xymatrix{0 \ar[r] & G \ar[r] & \mathbb{G}_{a} \ar[r] & \mathbb{G}_{a} \ar[r] & 0.}$$
Dans les deux cas, ceci implique que $K(G,2)$ est affine. \hfill $\Box$ \\

Nous allons maintenant donner une r\'eciproque au th\'eor\`eme \ref{t5}. Pour cela nous allons \'etudier les
faisceaux d'homotopie des champs de la forme $\mathbb{R}Spec\, A$, o\`u $A$ est une $k$-alg\`ebre
co-simpliciale augement\'ee $A \longrightarrow k$, et cohomologiquement connexe (i.e. $H^{0}(A)\simeq k$). Le
th\'eor\`eme suivant peut \^etre interpr\'et\'e comme une version g\'eom\'etrique de l'existence
des mod\`eles minimaux (voir \cite{bg,k}).

\begin{thm}\label{t5'}
Soit $A$ une $k$-alg\`ebre co-simpliciale de $\mathbb{U}$, augment\'ee $A \longrightarrow k$, et cohomologiquement connexe.
Notons $Spec\, k=* \longrightarrow \mathbb{R}Spec\, A=F$ le champ point\'e associ\'e. Alors, pour tout $i>0$, le
faisceau $\pi_{i}(F,*)$ est repr\'esentable par un sch\'ema en groupes affine et unipotent. De plus, le champ
$F$ est connexe (i.e. $\pi_{0}(F)=*$).
\end{thm}

\textit{Preuve:} Il s'agit en r\'ealit\'e de construire une d\'ecomposition de Postnikov de la $k$-alg\`ebre
co-simpliciale $A$. Dans le cas o\`u $A$ est $1$-connexe (i.e. on a de plus $H^{1}(A)=0$), l'existence de
cette d\'ecomposition est d\'emontr\'ee dans \cite[$3.2$]{k}. \\

La preuve consiste \`a montrer qu'il existe une tour de champs point\'es
$$\xymatrix{
F \ar[r] & \dots \ar[r] & F_{n} \ar[r] & F_{n-1} \ar[r] & \dots \ar[r] & F_{1} \ar[r] & F_{0}=*}$$
v\'erifiant les trois propri\'et\'es suivantes.

\begin{itemize}
\item Le champ $F_{n}$ est connexe et $n$-tronqu\'e (i.e. $\pi_{i}(F_{n},*)=0$ pour tout $i>n$). De plus, le
morphisme naturel $F_{n} \longrightarrow F_{n-1}$ induit une isomorphisme de champs
$\tau_{\leq n-1}F_{n} \simeq F_{n-1}$.

\item Pour tout $n$ et tout $i$, les faisceaux $\pi_{i}(F_{n},*)$ sont repr\'esentables par des
sch\'emas en groupes affines et unipotents.

\item Le morphisme $f_{n}^{*} : H^{i}(F_{n},\mathbb{G}_{a}) \longrightarrow H^{i}(F,\mathbb{G}_{a})$
est bijectif pour $i\leq n$, et injectif pour $i=n+1$.

\end{itemize}

En effet, supposons l'existence d'une telle tour et d\'emontrons le th\'eor\`eme.

\begin{lem}
Le champ $Holim_{n}F_{n}$ est affine.
\end{lem}

\textit{Preuve:}  Chaque $F_{n}$ est un champ affine d'apr\`es le th\'eor\`eme \ref{t5}. Ainsi, le lemme
se d\'eduit de la proposition \ref{p2'}. \hfill $\Box$ \\

\begin{lem}
Le morphisme $F \longrightarrow Holim_{n}F_{n}$ est une $\mathcal{O}$-\'equivalence.
\end{lem}

\textit{Preuve:} Comme les champs $F_{n}$, ainsi que $Holim_{n}F_{n}$, sont des champs affines, on a d'apr\`es
le corollaire \ref{c1}
$$\mathbb{L}\mathcal{O}(Holim_{n}F_{n})\simeq Hocolim_{n}\mathbb{L}\mathcal{O}(F_{n}).$$
Comme d'apr\`es le corollaire \ref{c2} on a
$$H^{i}(Hocolim_{n}\mathbb{L}\mathcal{O}(F_{n}))\simeq Colim_{n}H^{i}(F_{n},\mathbb{G}_{a}),$$
le lemme se d\'eduit alors imm\'ediatement de la troisi\`eme propri\'et\'e de la tour $\{F_{n}\}_{n}$. \hfill $\Box$ \\

Les deux lemmes pr\'ec\'edents impliquent que le morphisme $F \longrightarrow Holim_{n}F_{n}$ est
un isomorphism. Remarquons de plus que chaque faisceau $\pi_{i}(F_{n},*)$ est un sch\'ema en groupes unipotent,
et donc de dimension cohomologique inf\'erieure \`a $1$. Le corollaire \ref{p1'} permet alors de conclure. \\

Il nous reste maintenant \`a d\'emontrer l'existence de la tour de champs
$$\xymatrix{
F \ar[r] & \dots \ar[r] &  F_{n} \ar[r] & F_{n-1} \ar[r] & \dots \ar[r] & F_{1} \ar[r] & F_{0}=*}$$
v\'erifiant les trois conditions pr\'ec\'edentes. Nous allons proc\'eder par une double r\'ecurrence. Supponsons
pour cela que nous ayons d\'efini les $F_{i}$ pour $i<n+1$, et d\'efinissons $F_{n+1}$ de la fa\c{c}on suivante (cela
va prendre un certain temps). \\

Commen\c{c}ons par quelques mots sur la notion de faisceaux d'homologie de champs. Pour tout champ $F$, le
groupe $H^{i}(F,\mathbb{G}_{a})$ est muni de l'endomorphisme  induit par l'endomorphisme de Frobenius
$\mathbf{F} : \mathbb{G}_{a} \longrightarrow \mathbb{G}_{a}$. Plus pr\'ecisement, pour tout $Spec\, A \in Aff/k$,
on a $\mathbf{F}(a):=a^{p}$ pour $a \in A=\mathbb{G}_{a}(Spec\, A)$, et o\`u $p$ est l'exposant caract\'eristique
de $k$ (ainsi, $\mathbf{F}=Id$ si $k$ est de caract\'eristique nulle). Nous noterons $B$ l'anneau des
endomorphismes de groupes de $\mathbb{G}_{a}$. Ainsi, si $car\, k=p>0$,
$B\simeq k[\mathbf{F}]$ est l'anneau des polyn\^omes non-commutatifs avec un seul g\'en\'erateur $\mathbf{F}$, et une
seule relation $\mathbf{F}.a:=a^{p}.\mathbf{F}$. En contre partie, si $car\, k=0$, alors $B\simeq k$. Nous verrons
donc les groupes $H^{i}(F,\mathbb{G}_{a})$ comme des $B$-modules.

Rappelons alors que tout $B$-module $M$ d\'efinit un faisceau en groupes ab\'elien
$$\begin{array}{cccc}
H(M) : & (Aff/k)_{fpqc} & \longrightarrow & Ab \\
& Spec\, A & \mapsto & Hom_{B}(M,A),
\end{array}$$
o\`u toute $k$-alg\`ebre $A$ est vue comme un $B$-module par l'action $\mathbf{F}(a)=a^{p}$, pour $a \in A$.
D'apr\`es \cite[$IV$, \S $3$, Cor. $6.7$]{dg}, le foncteur $H$ d\'efinit une \'equivalence entre la cat\'egorie des
$B$-modules appartenant \`a $\mathbb{U}$, et celle des sch\'emas en groupes unipotents ab\'eliens et annul\'es par
le d\'ecalage (i.e. qui sont sous-groupes ferm\'es d'un produit $\mathbb{G}_{a}^{J}$, avec $J \in \mathbb{U}$).

Nous noterons $H_{i}(F):=H(H^{i}(F,\mathbb{G}_{a}))$, que nous appellerons les faisceaux d'homologie
de $F$. Ce sont des sch\'emas en groupes unipotents ab\'eliens d\`es que $H^{i}(F,\mathbb{G}_{a})$ est isomorphe \`a un object
de  $\mathbb{U}$ (e.g.
si $F$ est affine). Un exemple fondamental est $K(H,n)$, avec $H$ un sch\'ema en groupes
unipotent, ab\'elien et annul\'e par le d\'ecalage. Dans ce cas il est facile de v\'erifier que l'on a $H_{n}(K(H,n))\simeq H$.
Ainsi, pour un tel sch\'ema en groupes $H$, on dispose par fonctorialit\'e d'un morphisme
$$H^{n}(F,H)\simeq \pi_{0}\mathbb{R}\underline{Hom}(F,K(H,n)) \longrightarrow Hom(H_{n}(F),H).$$
Nous aurons alors besoin du lemme suivant, reliant homologie et cohomologie.

\begin{lem}\label{lt5'}
Soit $F$ un champ affine, et $H$ un sch\'ema en groupes
unipotent, ab\'elien et annul\'e par le d\'ecalage.
Alors le morphisme naturel
$$H^{n}(F,H) \longrightarrow Hom(H_{n}(F),H)$$
est surjectif.
\end{lem}

\textit{Preuve:} Soit $M$ un $B$-module tel que $H\simeq H(M)$, et
$$\xymatrix{0 \ar[r] & B^{J} \ar[r] & B^{I} \ar[r] & M \ar[r] & 0}$$
une r\'esolution libre du $B$-module $M$. Elle correspond \`a une suite exacte de faisceaux en groupes
$$\xymatrix{0 \ar[r] & H(M) \ar[r] & \mathbb{G}_{a}^{I} \ar[r] & \mathbb{G}_{a}^{J} \ar[r] & 0.}$$
On en d\'eduit un diagramme commutatif
$$\xymatrix{ 0 \ar[r] & Hom(H_{n}(F),H) \ar[r] & Hom(H_{n}(F),\mathbb{G}_{a})^{I} \ar[r] & Hom(H_{n}(F),\mathbb{G}_{a})^{J} \\
\dots \ar[r] & H^{n}(F,H) \ar[u] \ar[r] & H^{n}(F,\mathbb{G}_{a})^{I} \ar[u]^-{a} \ar[r] &
H^{n}(F,\mathbb{G}_{a})^{J} \ar[u]^-{b},}$$
o\`u $a$ et $b$ sont des isomorphismes (car $Hom(-,\mathbb{G}_{a})$ est le foncteur inverse de $H$). Ce diagramme montre que
le morphisme en question est surjectif. \hfill $\Box$ \\

Nous allons commencer par construire une tour de champs affines, point\'es et connexes
$$\xymatrix{F\ar[r] & \dots \ar[r] & F_{n,q} \ar[r] & F_{n,q-1} \ar[r] & \dots \ar[r]
& F_{n,1} \ar[r] & F_{n,0}=F_{n},}$$
puis poser $F_{n+1}:=Holim_{q}F_{n,q}$. Les champs $F_{n,q}$ seront d\'efinis par r\'ecurrence de sorte \`a ce que
la fibre homotopique de $F_{n,q+1} \longrightarrow F_{n,q}$ soit de la forme $K(\pi_{n+1,q+1},n+1)$, avec
$\pi_{n+1,q+1}$ un sch\'ema en groupes ab\'elien et unipotent. On montrera de plus par r\'ecurrence que le morphisme
induit $H^{n+1}(F_{n,q},\mathbb{G}_{a}) \longrightarrow H^{n+1}(F,\mathbb{G}_{a})$ est injectif. Le lecteur pourra comparer
la consrtuction des $F_{n,q}$ d\'ecrite ci-dessous avec la construction du mod\`ele minimal dans \cite[Prop. $7.7$]{bg}.

Supposons $F \longrightarrow F_{n,q} \longrightarrow F_{n}$ construit, et consid\'erons la suite exacte
courte
$$\xymatrix{0 \ar[r] & H^{n+1}(F_{n,q},\mathbb{G}_{a}) \ar[r] & H^{n+1}(F,\mathbb{G}_{a}) \ar[r] &
C \ar[r] & 0.}$$
L'endomorphisme de Frobenius $\mathbf{F} : \mathbb{G}_{a} \longrightarrow \mathbb{G}_{a}$ induit un endomorphisme
de la suite exacte pr\'ec\'edente. Cette suite exacte peut donc \^etre consid\'er\'ee comme une suite
exacte de $B$-modules, ou encore \`a travers le foncteur $H$ comme une suite exacte courte de sch\'emas en groupes
ab\'eliens et unipotents
$$\xymatrix{ 0 \ar[r] & H(C) \ar[r] & H_{n+1}(F) \ar[r] & H_{n+1}(F_{n,q}) \ar[r] & 0.}$$
Cette suite exacte de faisceaux ab\'eliens
est classifi\'ee par un morphisme dans la cat\'egorie d\'eriv\'ee des faisceaux ab\'eliens sur $(Aff/k)_{fpqc}$
$$e : H_{n+1}(F_{n,q})[n+1] \longrightarrow H(C)[n+2].$$
En oubliant la structure de groupes, ce morphisme induit \`a travers la correspondance de Dold-Kan un morphisme
de champs
$$e : K(H_{n+1}(F_{n,q}),n+1) \longrightarrow K(H(C),n+2).$$
De plus, le lemme \ref{lt5'} montre que l'on peut trouver un morphisme de champs $F_{n,q} \longrightarrow K(H_{n+1}(F_{n,q}),n+1)$
qui rel\`eve l'identit\'e de $H_{n+1}(F_{n,q})$. Compos\'e avec le morphisme d'extension $e$, on en d\'eduit
un morphisme de champs
$$F_{n,q} \longrightarrow K(H(C),n+2),$$
dont la fibre homotopique sera not\'ee $\widetilde{F_{n,q}} \longrightarrow F_{n,q}$. Le morphisme
$F \longrightarrow F_{n,q} \longrightarrow K(H(C),n+2)$ \'etant trivial dans $Ho(SPr(k))$, on en d\'eduit un
morphisme $F \longrightarrow \widetilde{F_{n,q}}$ qui rel\`eve $F \longrightarrow F_{n,q}$.
Il est alors facile
de v\'erifier que le morphisme induit
$$H^{n+1}(\widetilde{F_{n,q}},\mathbb{G}_{a}) \longrightarrow H^{n+1}(F,\mathbb{G}_{a})$$
est un isomorphisme. \\

Consid\'erons maintenant la suite exacte de $B$-modules
$$\xymatrix{ 0\ar[r] & K \ar[r] & H^{n+2}(\widetilde{F_{n,q}},\mathbb{G}_{a}) \ar[r] & H^{n+2}(F,\mathbb{G}_{a}),}$$
et la suite de sch\'emas en groupes ab\'eliens unipotents associ\'ee
$$\xymatrix{H_{n+2}(F) \ar[r] & H_{n+2}(\widetilde{F_{n,q}}) \ar[r] & H(K) \ar[r] &0.}$$
A l'aide du lemme \ref{lt5'} choississons un morphisme $\widetilde{F_{n,q}} \longrightarrow K(H(K),n+2)$ induisant la
projection $H_{n+2}(\widetilde{F_{n,q}}) \longrightarrow H(K)$, et d\'efinissons $F_{n,q+1}$ comme \'etant sa
fibre homotopique. Comme pr\'ec\'edemment, le morphisme $F \longrightarrow \widetilde{F_{n,q}}$ se factorise
par $F_{n,q+1}$. De plus, par construction il est facile de voir que le morphisme induit
$H^{n+2}(F_{n,q+1},\mathbb{G}_{a}) \longrightarrow H^{n+2}(F,\mathbb{G}_{a})$ envoit $K\subset
H^{n+2}(\widetilde{F_{n,q}},\mathbb{G}_{a})$ sur $0 \in H^{n+2}(F,\mathbb{G}_{a})$.
Pour finir, le lemme suivant montre
que le morphisme induit 
$H^{n+1}(F_{n,q+1},\mathbb{G}_{a}) \longrightarrow H^{n+1}(F,\mathbb{G}_{a})$ est injectif, et ceci termine
la construction du champ $F_{n,q+1}$.

\begin{lem}
Soit $p : F' \longrightarrow F$ un morphisme entre deux champs affines point\'es et connexes, dont la fibre
homotopique est isomorphe \`a $K(\pi,n+1)$. Alors, le morphisme induit en cohomologie
$$p^{*} : H^{n+1}(F,\mathbb{G}_{a}) \longrightarrow H^{n+1}(F',\mathbb{G}_{a})$$
est injectif.
\end{lem}

\textit{Preuve:} En utilisant les corollaires \ref{c1} et \ref{c2}, le morphisme $p$ est induit par un morphisme
de $k$-alg\`ebres co-simpliciales augment\'ees et cohomologiquement connexes $u : A \longrightarrow B$.
Soit $C:=k\otimes^{\mathbb{L}}_{A}B$ la cofibre homotopique de $u$ le long de l'augmentation $A \longrightarrow k$.
Alors, comme $K(\pi,n+1)\simeq \mathbb{R}Spec\, C$, on a $H^{i}(C)=0$ pour $0<i<n+1$. Les premiers
termes non-nuls de la suite spectrale d'Eilenberg-Moore (voir \cite[$2.2$]{k}) donnent donc une suite exacte
$$\xymatrix{0 \ar[r] & H^{n+1}(A) \ar[r] & H^{n+1}(B) \ar[r] & \dots .}$$
D'apr\`es le corollaire \ref{c2} ceci implique le lemme. \hfill $\Box$ \\

En appliquant le corollaire \ref{c1}, et le fait que les colimites filtrantes commutent
avec les foncteurs de cohomologie dans, on trouve des isomorphismes
$$H^{i}(F_{n+1},\mathbb{G}_{a})\simeq Colim_{q}H^{i}(F_{n,q},\mathbb{G}_{a}) \simeq
Colim_{q}H^{i}(\widetilde{F_{n,q}},\mathbb{G}_{a}).$$
Comme chaque morphism $H^{n+1}(\widetilde{F_{n,q}},\mathbb{G}_{a}) \longrightarrow H^{n+1}(F,\mathbb{G}_{a})$ est
un isomorphisme, on en d\'eduit que le morphisme
$H^{n+1}(F_{n+1},\mathbb{G}_{a}) \longrightarrow H^{n+1}(F,\mathbb{G}_{a})$ est un isomorphisme.

De plus, si $x \in Ker \left(H^{n+2}(F_{n+1},\mathbb{G}_{a}) \longrightarrow H^{n+2}(F,\mathbb{G}_{a})\right)$,
on peut repr\'esenter $x$ par un \'el\'ement dans le noyau de
$H^{n+2}(\widetilde{F_{n,q}},\mathbb{G}_{a}) \longrightarrow H^{n+2}(F,\mathbb{G}_{a})$. Par construction
de $F_{n,q+1}$, ceci implique que l'image
de $x$ dans $H^{n+2}(F_{n,q+1},\mathbb{G}_{a})$ est nulle. Ainsi, $x=0 \in H^{n+2}(F_{n+1},\mathbb{G}_{a})$.
Le morphisme $H^{n+2}(F_{n+1},\mathbb{G}_{a}) \longrightarrow H^{n+2}(F,\mathbb{G}_{a})$ est
donc injectif. Il nous reste ainsi \`a d\'emontrer que $F_{n+1}$ est $(n+1)$-tronqu\'e, que le morphisme
$\tau_{\leq n}F_{n+1}\longrightarrow F_{n}$ est un isomorphisme, et que $\pi_{n+1}(F_{n+1},*)$ est repr\'esentable par un
sch\'ema en groupes unipotent. \\

Commen\c{c}ons par montrer par r\'ecurrence sur $q$ que $F_{n,q}$ est $(n+1)$-tronqu\'e, et que
$\pi_{n+1}(F_{n,q+1},*) \longrightarrow \pi_{n+1}(F_{n,q},*)$ est un \'epimorphisme de sch\'emas en groupes
unipotents. Tout d'abord, $\widetilde{F_{n,q}}$ \'etant la fibre homotopique d'un morphisme $F_{n,q} \longrightarrow
K(H,n+2)$, on voit que $\widetilde{F_{n,q}}$ est $(n+1)$-tronqu\'e et que $\pi_{n+1}(\widetilde{F_{n,q}},*)$ est
une extension de $\pi_{n+1}(F_{n,q},*)$ par $H$. Comme $H$ est un sch\'ema en groupes unipotent par construction,
on voit par r\'ecurrence que $\pi_{n+1}(\widetilde{F_{n,q}},*)$ est un encore un sch\'ema en groupes unipotent. Ceci montre
aussi que $\tau_{\leq n}\widetilde{F_{n,q}}\simeq F_{n}$. De m\^eme, on montre que $F_{n,q+1}$ est $(n+1)$-tronqu\'e,
que $\tau_{\leq n}F_{n,q+1}\simeq F_{n}$, et que $\pi_{n+1}(F_{n,q+1},*)$ est un sch\'ema en groupes unipotent qui se
surjecte sur $\pi_{n+1}(F_{n,q},*)$.

Enfin, on voit d'apr\`es le lemme \ref{lt1} que la fibre homotopique de $F_{n+1} \longrightarrow F_{n}$ est
de la forme $Holim_{q}K(\pi_{n+1}(F_{n,q},*),n+1)\simeq K(Lim_{q}\pi_{n+1}(F_{n,q},*),n+1)$. Ainsi, le faisceau
$\pi_{n+1}(F_{n+1},*)$ est une $\mathbb{U}$-limite projective de sch\'ema en groupes unipotent et donc
est lui-m\^eme un sch\'ema en groupes unipotent. Ceci montre au passage que $F_{n+1}$ est $(n+1)$-tronqu\'e, et que
$\tau_{\leq n}F_{n+1}\simeq F_{n}$, et termine donc la preuve du th\'eor\`eme \ref{t5'}. \hfill $\Box$ \\

On peut maintenant rassembler les deux principaux r\'esultats de ce paragraphe en le corollaire suivant.

\begin{cor}
Le foncteur $\mathbb{R}Spec : Ho(k-Alg^{\Delta}/k) \longrightarrow Ho(SPr(k)_{*})$
induit une \'equivalence entre la sous-cat\'egorie pleine des $k$-alg\`ebres co-simpliciales augment\'ees
appartenant \`a $\mathbb{U}$ et cohomologiquement connexes, et la sous-cat\'egorie pleine
des champs point\'es et connexes dont les faisceaux d'homotopie sont re\-pr\'e\-sen\-ta\-bles par des sch\'emas en groupes
affines unipotents.
\end{cor}

Un autre corollaire important du th\'eor\`eme \ref{t5'} est la description du groupe fondamental d'une affination
d'un ensemble simplicial.

\begin{cor}\label{c12'}
Soit $X$ un ensemble simplicial point\'e et connexe appartenant \`a $\mathbb{U}$. Alors il existe un isomorphisme
naturel
$$\mathcal{U}(\pi_{1}(X,x))\simeq \pi_{1}((X\otimes k)^{uni},x).$$
\end{cor}

\textit{Preuve:} Nous savons d'apr\`es le th\'eor\`eme \ref{t5'} que $\pi_{1}((X\otimes k)^{uni},x)$ est un sch\'ema en
groupes affine et unipotent. De plus, d'apr\`es le th\'eor\`eme \ref{t5} et la propri\'et\'e universelle
des affinations, on a pour tout sch\'ema en groupes affine et unipotent $H$
$$Hom(\pi_{1}(X,x),H)\simeq \pi_{0}\mathbb{R}\underline{Hom}_{*}(X,K(H,1)) \simeq
\pi_{0}\mathbb{R}\underline{Hom}_{*}((X\otimes k)^{uni},K(H,1))$$
$$ \simeq Hom(\pi_{1}((X\otimes k)^{uni},x),H).$$
On conclut alors d'apr\`es la propri\'et\'e universelle de $\mathcal{U}(\pi_{1}(X,x))$. \hfill $\Box$ \\

\subsection{Comparaison avec l'homotopie rationnelle et $p$-adique}

Les r\'esultats de comparaison que nous pr\'esentons ici sont le th\'eor\`eme \ref{t8} et ses corollaires
\ref{c13} et \ref{c14}.
A vrai dire, il ne s'agit pas
\`a proprement parler de comparer rigoureusement nos conctructions \`a celles pr\'e-existantes, mais plut\^ot de
d\'emontrer des analogues des r\'esultats standards. Cette comparaison montrera que
les champs affines donnent une solution au probl\`eme de la sch\'ematisation tel qu'il
est \'enonc\'e dans l'appendice A. \\

Le th\'eor\`eme suivant est l'analogue des principaux r\'esultats en homotopie rationnelle et $p$-adique (voir \cite{q,s,bg,g,m1}).
Nous nous
restreindrons aux ensembles simpliciaux nilpotents et de type fini (sur $\mathbb{Z}$).
Pour nous, il s'agira des $\mathbb{U}$-ensembles
simpliciaux point\'es $X$, poss\'edant une d\'ecomposition de Postnikov
$\{\tau_{\leq n}X\}_{n}$, telle que pour tout $n$ il existe
une suite finie de diagrammes homotopiquement cart\'esiens (avec $0\leq i \leq m(n)$)
$$\xymatrix{
\tau_{\leq n-1}^{i}X \ar[r] \ar[d] & \tau_{\leq n-1}^{i}X \ar[d] \\
\bullet \ar[r] & K(M_{n}^{i},n+1),}$$
avec $M_{n}$ un groupe ab\'elien de type fini,
$\tau_{\leq n-1}^{0}X=\tau_{\leq n-1}X$ et $\tau_{\leq n-1}^{m(n)}X=\tau_{\leq n}X$.
Le lecteur pourra s'il le d\'esire adapter lui-m\^eme la preuve du th\'eor\`eme suivant
pour traiter le cas plus g\'en\'eral des espaces nilpotents et de $k$-type finis.

Pour les besoins de l'\'enonc\'e rappelons qu'un morphisme d'ensembles simpliciaux $f : X \longrightarrow Y$ est une $k$-\'equivalence,
si les morphismes induits $f^{*} : H^{i}(Y,k) \longrightarrow H^{i}(X,k)$ sont des isomorphismes pour tout $i\geq 0$.

\begin{thm}\label{t8}
\begin{itemize}
\item Si $k=\mathbb{Q}$, alors pour tout $\mathbb{U}$-ensemble simplicial point\'e, connexe, nilpotent et de type fini $X$,
le morphisme d'adjonction
$$X \longrightarrow \mathbb{R}\Gamma((X\otimes k)^{uni})$$
est une $\mathbb{Q}$-\'equivalence.

\item Si $k$ est alg\'ebriquement clos et de caract\'eristique $p>0$, alors pour tout $\mathbb{U}$-ensemble simplicial point\'e, connexe, nilpotent
et de type fini $X$, le morphisme d'adjonction
$$X \longrightarrow \mathbb{R}\Gamma((X\otimes k)^{uni})$$
est une $k$-\'equivalence.
\end{itemize}
\end{thm}

\textit{Preuve:} D'apr\`es la proposition \ref{paff} on sait que
$(X\otimes k)^{uni}\simeq \mathbb{R}Spec\, (k^{X})$.
Un argument standard de
d\'ecomposition de Postnikov et une utilisation du th\'eor\`eme d'Eilenberg-Moore (voir
\cite{lsm}) permettent alors de se ramener au cas \'el\'ementaire
$X=K(\mathbb{Z},n)$ (voir par exemple \cite{bg,m1}). Commen\c{c}ons alors par montrer que
$(X\otimes k)^{uni}\simeq K(\mathcal{U}(\mathbb{Z}),n)$, o\`u $\mathcal{U}(\mathbb{Z})$ est le compl\'et\'e
unipotent du groupe discret $\mathbb{Z}$.

Pour cela, le diagramme homotopiquement cart\'esien
$$\xymatrix{
K(\mathbb{Z},n-1) \ar[r] \ar[d] & \bullet \ar[d] \\
\bullet \ar[r] & K(\mathbb{Z},n),}$$
et une nouvelle application du th\'eor\`eme d'Eilenberg-Moore montrent que l'on a 
$$\mathbb{R}\Omega_{*}(K(\mathbb{Z},n)\otimes k)^{uni})\simeq (K(\mathbb{Z},n-1)\otimes k)^{uni}.$$ Or, d'apr\`es
le th\'eor\`eme \ref{t5'} on sait que
$(K(\mathbb{Z},n)\otimes k)^{uni}$ est connexe, et donc d'apr\`es le th\'eor\`eme \ref{t2} on a
$$(K(\mathbb{Z},n)\otimes k)^{uni})\simeq B\mathbb{R}\Omega_{*}(K(\mathbb{Z},n)\otimes k)^{uni})\simeq
B(K(\mathbb{Z},n-1)\otimes k)^{uni}.$$
Ceci permet donc de se ramener au cas o\`u $n=1$.

\begin{lem}
Le moprhisme naturel $H^{i}(B\mathcal{U}(\mathbb{Z}),\mathbb{G}_{a}) \longrightarrow H^{i}(B\mathbb{Z},\mathbb{G}_{a})$
est un isomorphisme pour tout $i$.
\end{lem}

\textit{Preuve:} Pour $i=0$ c'est imm\'ediat. Pour $i=1$, on a d'apr\`es la propri\'et\'e universelle de
$\mathcal{U}(\mathbb{Z})$
$$H^{1}(B\mathcal{U}(\mathbb{Z}),\mathbb{G}_{a})\simeq
Hom(\mathcal{U}(\mathbb{Z}),\mathbb{G}_{a})\simeq k\simeq H^{1}(B\mathbb{Z},\mathbb{G}_{a}).$$
Enfin, $\mathcal{U}(\mathbb{Z})$ \'etant un sch\'ema en groupes unipotent libre, on sait que
$H^{i}(B\mathcal{U}(\mathbb{Z}),\mathbb{G}_{a})$ est trivail pour $i>1$ (voir par exemple la preuve
de \ref{l15}). \hfill $\Box$ \\

Le lemme pr\'ec\'edent montre que le morphisme
$K(\mathbb{Z},1)\longrightarrow K(\mathcal{U}(\mathbb{Z}),1)$ est une $\mathcal{O}$-\'equivalence. Ainsi, comme
d'apr\`es le th\'eor\`eme \ref{t5}
$K(\mathcal{U}(\mathbb{Z}),1)$ est un champ affine, on a $(K(\mathbb{Z},1)\otimes \mathbb{Q})^{uni}\simeq
K(\mathcal{U}(\mathbb{Z}),1)$. Ceci ach\`eve la preuve que $(K(\mathbb{Z},n)\otimes k)^{uni}\simeq
K(\mathcal{U}(\mathbb{Z}),1)$. Remarquons que l'on n'a toujours pas utilis\'e l'hypoth\`ese sur $k$. \\

Il nous faut encore montrer que le morphisme d'adjonction
$$K(\mathbb{Z},n)\longrightarrow \mathbb{R}\Gamma(K(\mathcal{U}(\mathbb{Z}),n))$$
est une $k$-\'equivalence. Pour cela remarquons tout d'abord que le faisceau en groupes
$\mathcal{U}(\mathbb{Z})$ est acyclique sur $(Aff/k)_{fpqc}$. En effet, si $k$ est de caract\'eristique nulle
on a $\mathcal{U}(\mathbb{Z})\simeq \mathbb{G}_{a}$ qui est un faisceau quasi-coh\'erent et donc acyclique. Dans
le cas o\`u $k$ est de caract\'eristique $p$, on a $\mathcal{U}(\mathbb{Z})\simeq \mathbb{Z}_{p}$, le groupe
pro-fini des entiers $p$-adiques (vu comme un sch\'ema pro-alg\'ebrique, et donc affine sur $Spec\, k$). Mais comme
alors $k$ est alg\'ebriquement clos on sait que
$$H^{i}_{fpqc}(Spec\, k,\mathbb{Z}_{p})\simeq H^{i}_{et}(Spec\, k,\mathbb{Z}_{p})\simeq 0$$
d\`es que $i>0$.

On peut donc conclure que lorsque $k=\mathbb{Q}$ on a
$$\mathbb{R}\Gamma((\mathcal{U}(\mathbb{Z}),n))\simeq K(\mathbb{Q},n),$$
et lorsque $k$ est de caract\'eristique $p$ et alg\'ebriquement clos on a
$$\mathbb{R}\Gamma((\mathcal{U}(\mathbb{Z}),n))\simeq K(\mathbb{Z}_{p},n).$$
Dans les deux cas, il est bien connu que les morphismes d'adjonction
$$K(\mathbb{Z},n) \longrightarrow K(\mathbb{Q},n) \qquad K(\mathbb{Z},n)\longrightarrow K(\mathbb{Z}_{p},n)$$
sont des $k$-\'equivalences. \hfill $\Box$ \\

\textit{Remarque:} Le lecteur pourra par exemple comparer la preuve du point $(2)$ du th\'eor\`eme pr\'ec\'edent avec
la preuve du th\'eor\`eme principal de \cite{m1}, o\`u le point crucial est de montrer que $K(\mathbb{Z}/p,n)$ est
$k$-r\'esolvable (voir \cite[$4$]{m1}). Le fait que $K(\mathbb{Z}/p,n)$ soit $k$-r\'esolvable est en r\'ealit\'e tout \`a fait analogue
\`a la condition que $K(\mathbb{Z}/p,n)$ soit un champ affine, qui dans notre cas est une cons\'equence directe
de la th\'eorie. Il est par ailleurs int\'eressant de constater que la construction
du mod\`ele explicite
pour la $E_{\infty}$-alg\`ebre $C^{*}(K(\mathbb{Z}/p,n),k)$ qui est pr\'esent\'ee dans \cite[$4$]{m1},
est tr\`es proche de la construction
qui exibe le champ $K(\mathbb{Z}/p,n)$ comme la fibre homotopique du morphisme $1-\mathbf{F} : K(\mathbb{G}_{a},n) \longrightarrow
K(\mathbb{G}_{a},n)$, o\`u $\mathbf{F}$ est l'endomorphisme de $\mathbb{G}_{a}$ qui \'el\`eve \`a la puissance $p$. \\

Le corollaire suivant est un corollaire important de la d\'emonstration du th\'eor\`eme pr\'ec\'edent.

\begin{cor}\label{ct8}
Soit $X$ un $\mathbb{U}$-ensemble simplicial point\'e, connexe, nilpotent et de type fini. Alors, pour tout $i>0$,
le faisceau $\pi_{i}((X\otimes k)^{uni},x)$ est repr\'esentable par le compl\'et\'e unipotent du groupe
$\pi_{i}(X,x)$. En particulier, pour $i>1$ on a
$$\pi_{i}((X\otimes k)^{uni},x)\simeq \left\{
\begin{array}{cl}
\pi_{i}(X,x)\otimes_{\mathbb{Z}}\mathbb{G}_{a} & \text{si car k=0} \\
\pi_{i}(X,x)\otimes_{\mathbb{Z}}\mathbb{Z}_{p} & \textit{si car k=p}.
\end{array} \right.$$
\end{cor}

Le $(-\otimes k)^{uni}$ envoient clairement les
$k$-\'equivalences sur des isomorphismes dans $Ho(SPr(k)_{*})$. Il se factorise donc en un foncteur
$$(-\otimes k)^{uni} : Ho_{k}(\mathbb{U}-SEns_{0}) \longrightarrow Ho(SPr(k)_{*}),$$
o\`u $Ho_{k}(\mathbb{U}-SEns_{0})$ est la cat\'egorie localis\'ee de $Ho(\mathbb{U}-SEns_{0})$ le long
le long des $k$-\'equivalences.

\begin{cor}\label{c13}
Lorsque $k$ est ou bien le corps des rationnels $\mathbb{Q}$, ou bien alg\'e\-bri\-que\-ment clos et de caract\'eristique $p>0$, le
foncteur
$$(-\otimes k)^{uni} : Ho_{k}(\mathbb{U}-SEns_{0}) \longrightarrow Ho(SPr(k)_{*})$$
restreint \`a la sous-cat\'egorie pleine des ensembles simpliciaux nilpotents et de type fini, est pleinement fid\`ele.
\end{cor}

\textit{Preuve:} C'est imm\'ediat d'apr\`es le th\'eor\`eme \ref{t8}. \hfill $\Box$ \\

C'est le corollaire  pr\'ec\'edent qui permet d'affirmer que le champ $(X\otimes k)^{uni}$ est un mod\`ele pour l'homotopie
rationnel ou bien $p$-adique, suivant que $k=\mathbb{Q}$ ou bien $k=\overline{\mathbb{F}}_{p}$.

\begin{cor}\label{c14}
Soit $X$ et $Y$ deux $\mathbb{U}$-ensembles simpliciaux point\'es, connexes, nilpotents et de type fini, et supposons que
$k$ soit de caract\'eristique $p>0$. Si $(X\otimes k)^{uni}$ et $(Y\otimes k)^{uni}$ sont isomorphes dans $Ho(SPr(k)_{*})$, alors
$X$ et $Y$ sont isomorphes dans $Ho_{k}(\mathbb{U}-SEns_{0})$.
\end{cor}

\textit{Preuve:} En effet, si $\overline{k}$ est une cl\^oture alg\'ebrique de $k$, alors d'apr\`es
le corollaire \ref{c12} $(X\otimes \overline{k})^{uni}\simeq (Y\otimes \overline{k})^{uni}$. Cet isomorphisme induit d'apr\`es
le th\'eor\`eme \ref{t8} un isomorphisme entre $X$ et $Y$ dans $Ho_{\overline{k}}(\mathbb{U}-SEns_{0})$. Or, il est facile
de v\'erifier qu'un morphisme est une $\overline{k}$-\'equivalence si et seulement si c'est une $k$-\'equivalence.
Ainsi, il existe une \'equivalence de cat\'egories entre $Ho_{\overline{k}}(\mathbb{U}-SEns_{0})$ et  $Ho_{k}(\mathbb{U}-SEns_{0})$. Les ensembles
simpliciaux $X$ et $Y$ sont donc isomorphes dans le cat\'egorie $Ho_{k}(\mathbb{U}-SEns_{0})$. \hfill $\Box$ \\

\subsection{Critique des champs affines}

Pour terminer ce chapitre nous souhaitons faire la critique des champs affines et du proc\'ed\'e
d'affination des types d'homotopie. Nous avons vus que les champs affines sont tous construit \`a partir
d'extensions successives de champs de la forme $K(\mathbb{G}_{a},n)$. Ainsi, d'un point de vue de la th\'eorie
de l'homotopie les champs affines sont des types d'homotopie \textit{pro-nilpotents}. Ceci implique en particulier
que le foncteur d'affination tue toute l'information non nilpotente que peut poss\'eder un type d'homotopie. Par exemple,
si $X$ est un type d'homotopie tel que son groupe fondamental op\`ere de fa\c{c}on semi-simple sur ses groupes
d'homotopie sup\'erieur, le champ $(X\otimes \mathbb{Z})^{uni}$ ne verra pas une partie importante de l'information homotopique
que contient $X$. A titre d'exemple, examinons le cas tr\`es simple suivant.

Soit $X=K(\mathbb{Z}/m,\mathbb{Z}^{r},n)$, o\`u $\mathbb{Z}/m$ op\`ere non trivialement
sur le groupe ab\'elien $\mathbb{Z}^{r}$, et supposons que $k=\mathbb{Q}$.
Cette action s'\'etend en une action de $\mathbb{Z}/m$
sur le sch\'ema en groupes affine $\mathbb{G}_{a}^{r}$, et donne une
d\'ecomposition $\mathbb{G}_{a}^{r}\simeq \mathbb{G}_{a}^{i}\times \mathbb{G}_{a}^{r-i}$, o\`u le premier facteur
est le sous-groupe des points fixes. On peut alors v\'erifier que
$$\pi_{1}((X\otimes \mathbb{Q})^{uni},*)\simeq 0 \qquad \pi_{n}((X\otimes \mathbb{Q})^{uni},*)\simeq \mathbb{G}_{a}^{i},$$
Ainsi, la partie non triviale de l'action de $\mathbb{Z}/m$ sur $\pi_{n}(X)$ disparait dans
le champ $(X\otimes \mathbb{Q})^{uni}$.

De fa\c{c}on beaucoup plus g\'en\'erale, toute l'information concernant les re\-pr\'e\-sen\-ta\-tions non-nilpotentes des groupes
fondamentaux dans les groupes d'homotopie sup\'erieurs sera inaccessible \`a l'aide de la construction
$(-\otimes \mathbb{Z})^{uni}$. \\

L'id\'ee naturelle pour rem\'edier \`a ce probl\`eme est d'\'elargir la d\'efinition de champs affines afin de permettre
tous les sch\'emas en groupes affines comme des groupes fondamentaux potentiels (ce point de vue est bien entendu
reli\'e \`a la dualit\'e de Tannaka). C'est dans cette optique que nous introduirons les $\infty$-gerbes affines et les
types d'homotopie sch\'ematiques dans le chapitre suivant.

\section{$\infty$-Gerbes affines et types d'homotopie sch\'ematiques}

Dans ce dernier chapitre nous introduisons la seconde notion fondamentale de ce travail, celle d'$\infty$-gerbe affine et
de type d'homotopie sch\'ematique. L'id\'ee est d'utiliser la notion de champ affine afin de donner une version homotopique
de la notion de gerbe affine utilis\'ee dans le formalisme Tannakian (voir \cite{sa,d2}). Pour cela, rappelons qu'une
gerbe affine neutralis\'ee peut s'interpr\'eter comme un champ de la forme $K(G,1)$, o\`u $G$ est un sch\'ema en groupes affine
et plat sur l'anneau de base $k$. On pourrait aussi dire qu'il s'agit d'un champ point\'e et connexe dont le champ
des lacets est repr\'esentable par un sch\'ema affine et plat. Ce point de vue nous incite tout naturellement \`a d\'efinir
une $\infty$-gerbe affine point\'ee comme un champ point\'e et connexe dont le champ des lacets est un champ affine et
\textit{plat} sur $k$. C'est la d\'efinition que nous voulons adopter. Cependant, pour des raisons techniques nous ne voulons
pas d\'efinir la notion de platitude requise, et nous nous restreindrons donc au cas o\`u l'anneau de base est un corps.

D'un point de vue technique, le fait que les champs affines soient des objets $\mathcal{O}$-locaux est fondamental pour
que la notion d'affination d\'efinie dans \ref{daff} reste raisonnable (i.e. poss\`ede une propri\'et\'e d'unicit\'e et de
fonctorialit\'e). De m\^eme,
pour pouvoir par la suite d\'efinir une notion raisonnable de \textit{sch\'ematisation}, il est \`a pr\'evoir que les types
d'homotopie sch\'ematiques devront satisfaire \`a une certaine condition de localit\'e pour une th\'eorie cohomologique
convenable. C'est pour cela que nous introduirons les notions de $P$-\'equivalence et d'objets $P$-locaux,
qui g\'en\'eralisent les notions de $\mathcal{O}$-\'equivalences et d'objets $\mathcal{O}$-locaux en ce sens qu'elles
concernent la cohomologie \`a valeurs dans tous les syst\`emes locaux de $k$-espaces vectoriels.
Nous d\'efinirons alors les types
d'homotopie sch\'ematiques comme les $\infty$-gerbes affines point\'ees qui sont de plus des objets $P$-locaux. Cependant,
nous conjecturons que toute $\infty$-gerbe affine point\'ee est un champ $P$-local, et donc que
les types d'homotopie sch\'ematiques point\'es coincident avec les $\infty$-gerbes affines point\'ees (voir \ref{conj}).

Une fois ces notions d\'efinies nous d\'efinirons la notion de sch\'ematisation d'un type d'homotopie, et nous montrerons
que tout type d'homotopie dans $\mathbb{U}$ poss\`ede une sch\'ematisation. Nous montrerons alors \`a travers quelques exemples
en quoi le foncteur de sch\'ematisation diff\`ere de celui d'affination, et en quoi il r\'esoud le probl\`eme soulev\'e
au cours de paragraphe $2.6$. \\

Pour tout ce chapitre $k$ est un corps de caract\'eristique quelconque.

\subsection{D\'efinitions}

Dans ce paragraphe nous donnerons la d\'efinition d'$\infty$-gerbe affine, et de type d'homotopie
sch\'ematique. Nous n'\'etudierons essentiellement aucune de leurs
propri\'et\'es, et nous nous contenterons de donner une m\'ethode permettant de construire des exemples
\`a partir de sch\'emas en groupes affines simpliciaux. Celle-ci sera utilis\'ee
dans par la suite pour construire le foncteur de sch\'ematisation.

Le lecteur remarquera que nous ne donnerons la d\'efinition de type d'homoto\-pie sch\'ematique que lorsque $k$ est un corps.
De plus, la conjecture \ref{conj} pr\'evoit que
les types d'homotopie sch\'ematiques point\'es sont exactement les $\infty$-gerbes affines point\'es. Cependant, comme
un tel r\'esultat n'est certainement pas vraie sur une base g\'en\'erale (pour une g\'en\'eralisation ad\'equate
des notions d'$\infty$-gerbes affines et de types d'homotopie sch\'ematiques), faire la diff\'erence entre les
$\infty$-gerbes affines point\'ees et les types d'homotopie sch\'ematiques point\'es me semble raisonnable. \\

Pour commencer nous allons d\'efinir la notion de $P$-\'equivalence entre pr\'e\-faisc\-eaux simpliciaux point\'es, qui est
une notion plus forte que celle de $\mathcal{O}$-\'equivalence car faisant aussi intervenir les syst\`emes locaux de
$k$-espaces vectoriels de dimension finie.

\begin{df}\label{d12}
\begin{itemize}
\item
\emph{Un morphisme entre pr\'efaisceaux simpliciaux point\'es et connexes $f : G \longrightarrow H$ est une $P$-\'equivalence, si
pour tout sch\'ema en groupes affine et de type fini $K$, toute repr\'esentation lin\'eaire de dimension finie $\mathbb{V}$
de $K$, et tout entier $n>0$, le morphisme induit
$$f^{*} : \mathbb{R}\underline{Hom}_{*}(H,K(K,V,n)) \longrightarrow \mathbb{R}\underline{Hom}_{*}(G,K(K,V,n))$$
est un isomorphisme dans $Ho(\mathbb{V}-SEns)$.}
\item
\emph{Un pr\'efaisceau simplicial point\'e et connexe $F$ est dit $P$-local si pour toute $P$-\'equivalence $f : G \longrightarrow H$,
le morphisme induit
$$f^{*} : \mathbb{R}\underline{Hom}_{*}(H,F) \longrightarrow \mathbb{R}\underline{Hom}_{*}(G,F)$$
est un isomorphisme dans $Ho(\mathbb{V}-SEns)$.}
\end{itemize}
\end{df}

On remarquera qu'un morphisme de pr\'efaisceaux simpliciaux point\'es et connexes qui est une
$\mathcal{O}$-\'equivalence est aussi une $P$-\'equivalence. \\

\textit{Remarque:} Nous avons choisi l'expression $P$-\'equivalence par r\'ef\'erence aux complexes
parfaits. En effet, on peut montrer qu'un morphisme
est une $P$-\'equi\-va\-len\-ce si et seulement si c'est une \'equivalence pour la th\'eorie cohomologique non-ab\'elienne d\'etermin\'ee
par le $1$-champ de Segal des complexes parfaits (voir \cite[$\S 21$]{s1}). C'est ce point de vue plus g\'en\'eral qui
devra \^etre utilis\'e pour donner une d\'efinition raisonnable des types d'homotopie sch\'ematiques sur
des bases plus g\'en\'erales. \\

Rappelons que pour un pr\'efaisceau simplicial $* \longrightarrow F$, on peut d\'efinir son pr\'efaisceau simplicial
des lacets $\Omega_{*}F$. Le foncteur $\Omega_{*}$ \'etant alors un foncteur de Quillen \`a droite de la cat\'egorie
des pr\'efaisceaux simpliciaux point\'es dans elle m\^eme, nous noterons $\mathbb{R}\Omega_{*}$ son foncteur
d\'eriv\'e \`a droite.

\begin{df}\label{d10}
\begin{itemize}
\item \emph{Un pr\'efaisceau simplicial point\'e $s : * \longrightarrow F$ est une $\infty$-gerbe affine point\'ee sur $k$ s'il
v\'erifie les deux conditions suivantes.}
\begin{itemize}
\item \emph{Le pr\'efaisceau simplicial $F$ est connexe (i.e. $\pi_{0}(F)\simeq *$).}

\item \emph{Le pr\'efaisceau simplicial $\mathbb{R}\Omega_{*}F$ est un champ affine sur $k$.}
\end{itemize}

\item \emph{Un type d'homotopie sch\'ematique point\'e est une
$\infty$-gerbe affine  point\'ee qui est de plus $P$-local (en tant que pr\'efaisceau
simplicial point\'e).}

\end{itemize}
\end{df}

\textit{Remarques:} \begin{itemize}

\item Par d\'efinition une $\infty$-gerbe affine point\'ee est un pr\'efaisceau simplicial point\'e et sera donc toujours
consid\'er\'ee comme un objet dans $SPr(k)_{*}$ ou dans $Ho(SPr(k)_{*})$. Ainsi, un morphisme entre deux tels objets sera toujours un morphisme
d'objets point\'es.

\item Soit $F=K(G,1)$ une gerbe neutre sur $Spec\, k$, li\'ee par un sch\'ema en groupes affine $G$ (voir \cite{sa,d2}). Alors, vu comme
pr\'efaisceau simplicial $F$ est une $\infty$-gerbe affine (voir \cite{j2} pour le passage des champs en groupo\"{\i}des
aux pr\'efaisceaux simpliciaux). Ainsi, tout comme la notion de champ affine est une g\'en\'eralisation de celle de sch\'ema affine, la
notion d'$\infty$-gerbe affine g\'e\-n\'e\-ra\-li\-se celle de gerbe affine neutre. Il faut donc penser \`a la th\'eorie des gerbes
affines neutres comme \`a la th\'eorie homotopique sch\'ematique $1$-tronqu\'ee (et connexe).

\item Nous ne d\'efinirons pas la notion de type d'homotopie sch\'ematiques sur un anneau $k$ quelconque, la notion de
$P$-\'equivalences dans ce cas passant par l'existence du \textit{$1$-champ de Segal des complexes parfaits} que nous ne souhaitons
pas d\'efinir ici (voir cependant \cite{s1}).

\end{itemize}

Il est important de noter le fait \'el\'ementaire suivant.

\begin{lem}\label{l8}
Toute $\infty$-gerbe affine point\'ee est un champ sous-affine.
\end{lem}

\textit{Preuve:} Consid\'erons un remplacement cofibrant dans $Pr-\Delta^{o}-SPr(k)$, $H \longrightarrow \mathbb{R}\Omega_{*}F$,
ainsi que l'objet co-simplicial de $k-Alg^{\Delta}$
$$\begin{array}{cccc}
\mathcal{O}(H_{\bullet}) : & \Delta^{o} & \longrightarrow & k-Alg^{\Delta} \\
& [n] & \mapsto & \mathcal{O}(H_{n}).
\end{array}$$
Soit $A^{(\bullet)} \longrightarrow \mathcal{O}(H_{\bullet})$ un remplacement cofibrant de ce diagramme, avec $A^{(0)}=H_{0}=k$.
Alors, comme chaque $H_{n}$ est un champ
affine (car \'equivalent \`a $\mathbb{R}\Omega_{*}(F)^{n}$), le morphisme d'adjonction
$$\mathbb{R}\Omega_{*}F \longrightarrow Spec\, A^{(\bullet)}$$
est un isomorphisme dans $Ho(Pr-\Delta^{o}-SPr(k))$. En passant aux classifiants, on trouve un isomorphisme dans $Ho(SPr(k)_{*})$
$$F\simeq B\mathbb{R}\Omega_{*}F \simeq BSpec\, A^{(\bullet)}.$$
Comme le pr\'efaisceau simplicial $BSpec\, A^{(\bullet)}$ est repr\'esent\'e par le sch\'ema affine simplicial $[n] \mapsto Spec\, A^{(n)}_{n}$,
$F$ est bien un champ sous-affine. \hfill $\Box$ \\

Une $\infty$-gerbe affine point\'ee $F$ est la m\^eme chose qu'un $H_{\infty}$-champ dont le champ sous-jacent
est un champ affine. Comme la notion de $H_{\infty}$-champ est une notion diagrammatique, on en d\'eduit
une notion duale au niveau des $k$-alg\`ebres co-simpliciales. Ce n'est pas une notion que nous utiliserons, mais
nous donnerons sa d\'efinition \`a titre indicatif.

\begin{df}\label{hopf}
\emph{Une $H_{\infty}$-alg\`ebre de Hopf sur $k$ est la donn\'ee d'un foncteur
$$A^{(\bullet)} : \Delta \longrightarrow k-Alg^{\Delta}$$
tel que les trois propri\'et\'es suivantes soient satisfaites.}
\begin{itemize}
\item $A^{(0)}=k$

\item \emph{Pour tout $n>0$, le morphisme naturel}
$$\underbrace{A^{(1)}\otimes^{\mathbb{L}}_{k}A^{(1)}\otimes^{\mathbb{L}}_{k} \dots \otimes^{\mathbb{L}}_{k}A^{(1)}}_{n \; fois}
 \longrightarrow A^{(n)}$$
\emph{est un isomorphisme dans $Ho(k-Alg^{\Delta})$.}

\item \emph{Vue dans la cat\'egorie homotopique $Ho(k-Alg^{\Delta})$, la loi naturelle de co-mono\"{\i}de induite sur $A_{1}$
est une loie de co-groupe.}
\end{itemize}
\end{df}

Bien entendu, il existe une \'equivalence de la cat\'egorie homotopique des $H_{\infty}$-alg\`ebres de Hopf sur $k$
et appartenant \`a $\mathbb{U}$
et de la cat\'egorie des $\infty$-gerbes affines point\'ees sur $k$. Cette \'equivalence ce d\'eduit imm\'ediatement
du corollaire \ref{c1}. \\

Nous terminerons ce chapitre par un proc\'ed\'e g\'en\'eral de construction d'$\infty$-gerbes affines. Il permet de
donner de tr\`es nombreux exemples, mais nous ne savons pas s'ils sont tous obtenus de cette fa\c{c}on. Nous montrerons
par la suite que cette construction fournit en r\'ealit\'e des types d'homotopie sch\'ematiques, ce qui nous
permettra de construire le foncteur de sch\'ematisation. \\

Consid\'erons un sch\'ema en groupes affine simplicial $G$, et notons encore $G \in GpSPr(k)$ le pr\'efaisceau en groupes simplicial
qu'il repr\'esente. La loi de groupe $\mu : G\times G \longrightarrow G$ induit par fonctorialit\'e une co-multiplication
$$\mu^{*} : \mathcal{O}(G) \longrightarrow \mathcal{O}(G\times G)\simeq \mathcal{O}(G)\otimes \mathcal{O}(G),$$
qui fait de $\mathcal{O}(G)$ une $k$-alg\`ebre de Hopf commutative co-simpliciale.

On peut alors montrer que
$\mathcal{O}(G)\otimes \mathcal{O}(G)$ est isomorphe dans $Ho(k-Alg^{\Delta})$ au produit tensoriel
d\'eriv\'e $\mathcal{O}(G)\otimes^{\mathbb{L}} \mathcal{O}(G)$. En particulier, ceci montre que le morphisme naturel
$$\mathbb{R}Spec\, \mathcal{O}(G\times G) \longrightarrow \mathbb{R}Spec\, \mathcal{O}(G)\times \mathbb{R}Spec\,
\mathcal{O}(G)$$
est une \'equivalence. Plus g\'en\'eralement, ceci implique que $\mathcal{O}(G)$, vu comme objet dans la cat\'egorie homotopique
$Ho(k-Alg^{\Delta})$ est un objet en co-groupes (i.e. un objet en groupes dans la cat\'egorie oppos\'ee). Le diagramme
simplicial
$$\begin{array}{cccc}
\mathbb{R}Spec\, \mathcal{O}(G) : & \Delta^{o} & \longrightarrow & SPr(k) \\
& [n] & \mapsto & \mathbb{R}Spec\, \mathcal{O}(G^{n}),
\end{array}$$
est donc un $H_{\infty}$-champ au sens de la d\'efinition \ref{d5}, que nous noterons symboliquement
$\mathbb{R}Spec\, \mathbb{L}\mathcal{O}(G)$.
De fa\c{c}on duale, on peut aussi affirmer que le diagramme
co-simplicial
$$\begin{array}{cccc}
\mathbb{L}\mathcal{O}(G) : & \Delta^{o} & \longrightarrow & k-Alg^{\Delta} \\
& [n] & \mapsto & \mathbb{L}\mathcal{O}(G^{n}),
\end{array}$$
est une $H_{\infty}$-alg\`ebre de Hopf (voir Def. \ref{hopf}).

Le th\'eor\`eme \ref{t2} implique alors que $F=B\mathbb{R}Spec\, \mathbb{L}\mathcal{O}(G)$ est une $\infty$-gerbe affine
point\'ee, avec $\mathbb{R}\Omega_{*}(F)\simeq \mathbb{R}Spec\, \mathcal{O}(G)$. Remarquer de plus qu'il existe un morphisme
naturel dans $Ho(SPr(k)_{*})$ induit par l'adjonction entre les foncteurs $\mathcal{O}$ et $Spec$
$$BG \longrightarrow F.$$

L'int\'er\^et de cette construction est qu'elle donne un objet poss\'edant une propri\'et\'e universelle.

\begin{lem}\label{l9}
Sous les hypoth\`eses pr\'ec\'edentes, le morphisme naturel 
$$BG \longrightarrow F=B\mathbb{R}Spec\, \mathcal{O}(G)$$
est universel dans $Ho(SPr(k)_{*})$ pour les morphismes de $BG$ vers des $\infty$-gerbes affines point\'ees.
\end{lem}

\textit{Preuve:} Par d\'efinition des  $\infty$-gerbes affines point\'ees il suffit de montrer que le morphisme d'adjonction
$G \longrightarrow \mathbb{R}Spec\, \mathcal{O}(G)$, vu comme un morphisme dans $Ho(Pr-\Delta^{o}-SPr(k))$,
est universel vers les $H_{\infty}$-champs dont le champ sous-jacent est affine. Mais, ceci n'est qu'une version en famille
du corollaire \ref{c1}. \hfill $\Box$ \\

\textit{Remarque:} La question de savoir si toutes les $\infty$-gerbes affines point\'ees  sont obtenues
par cette construction est en fait \'equivalente \`a un probl\`eme de \textit{strictification des $H_{\infty}$-alg\`ebres de Hopf en
alg\`ebre de Hopf co-simpliciales}.

\subsection{Exemples d'$\infty$-gerbes affines et de types d'homotopie sch\'ematiques}

Etant donn\'ee leurs d\'efinitions il est assez difficile de produire des exemples explicites de $\infty$-gerbes affines. Cela
provient principalement du fait que l'on connait assez mal leurs faisceaux d'homotopie en g\'en\'eral. Pour rem\'edier
\`a cela nous allons d\'emontrer un crit\`ere de reconnaissance des types d'homotopie sch\'ematiques analogues au th\'eor\`eme
\ref{t5}. Ceci permet de donner de tr\`es nombreux exemples de types d'homotopie sch\'ematiques, dont
en particulier les champs tr\`es pr\'esentables introduits par C. Simpson (voir
\cite{s3} et \cite{s4}). \\

Nous commen\c{c}erons par citer deux lemmes utiles pour construire des types d'homotopie sch\'ematiques.

\begin{lem}\label{l10}
Tout champ affine point\'e et connexe est un type d'homotopie sch\'e\-ma\-ti\-que point\'e.
\end{lem}

\textit{Preuve:} Soit $A$ une $k$-alg\`ebre co-simpliciale de $\mathbb{U}$ telle que $F$ soit \'equivalent \`a $\mathbb{R}Spec\, A$, et
$A \longrightarrow k$ l'augmentation correspondant au point de $F$. Il est alors facile de voir qu'il existe un
isomorphisme dans la cat\'egorie homotopique $Ho(SPr(k))$
$$\mathbb{R}\Omega_{*}F\simeq \mathbb{R}Spec\, (k\otimes^{\mathbb{L}}_{A}k).$$
Comme $k$ est corps ceci montre bien que $\mathbb{R}\Omega_{*}F$ est un champ affine, et donc que $F$ est une
$\infty$-gerbe affine. De plus, il est imm\'ediat qu'un objet $\mathcal{O}$-local est aussi un objet $P$-local (car une
$P$-\'equivalence est une $\mathcal{O}$-\'equivalence). Par d\'efinition, $F$ est bien un type d'homotopie sch\'ematique
point\'e. \hfill $\Box$ \\

\begin{lem}\label{l10'}
Tout champ point\'e et connexe qui est une $\mathbb{U}$-limites homotopique de types d'homotopie
sch\'e\-ma\-ti\-ques point\'es et lui-m\^eme un type d'homotopie sch\'e\-ma\-ti\-que point\'e et connexe.
\end{lem}

\textit{Proof:} En appliquant \ref{p2'} aux champs des lacets on montre que le lemme est vrai pour
les $\infty$-gerbes affines point\'ees. De plus, il est clair que les objets $P$-locaux sont stables par
$\mathbb{U}$-limites homotopiques. \hfill $\Box$ \\

\begin{lem}\label{l11}
Soit $F$ un pr\'efaisceau simplicial point\'e et connexe, et supposons que $\mathbb{R}\Omega_{*}F$ soit sous-affine et que $F$
soit une $\mathbb{V}$-limite homotopique (dans la cat\'egorie des pr\'efaisceaux simpliciaux point\'es) de types d'homotopie sch\'ematiques point\'es.
Alors $F$ est un type d'homotopie sch\'ematique point\'e.
\end{lem}

\textit{Preuve:} Comme $k$ est un corps tout les champs affines sont plats sur $k$. De plus, le foncteur $\mathbb{R}\Omega_{*}$ \'etant de
Quillen \`a droite il pr\'eserve les limites homotopiques. Le lemme est donc une cons\'equence du corollaire \ref{c4},
et du fait que les objets $P$-locaux sont stables par $\mathbb{V}$-limites homotopiques. \hfill $\Box$ \\

Le th\'eor\`eme suivant est un analogue du th\'eor\`eme \ref{t5}.

\begin{thm}\label{t6}
Soit $* \longrightarrow F$ un pr\'efaisceau simplicial point\'e et connexe. Supposons que les deux assertions
suivantes soient satisfaites.
\begin{itemize}
\item Le faisceau en groupes $\pi_{1}(F,*)$ est repr\'esent\'e par un sch\'ema en groupe affine.
\item Pour tout entier $i>1$ le faisceau en groupes $\pi_{i}(F,*)$ est repr\'esent\'e par un sch\'ema en groupes affine
et unipotent.
\end{itemize}
Alors $F$ est est un type d'homotopie sch\'ematique point\'e.
\end{thm}

\textit{Preuve:}  Commen\c{c}ons par montrer que
$\mathbb{R}\Omega_{*}F$ est affine (i.e. que $F$ est une $\infty$-gerbe affine point\'ee).

Le m\^eme argument que celui utilis\'e dans la preuve du th\'eor\`eme \ref{t5} montre que $F$ satisfait aux hypoth\`eses
de la proposition \ref{p1}.
Un argument de d\'e\-com\-po\-si\-tion de Postnikov et une utilisation du lemme \ref{l10'} montre alors qu'il suffit de traiter le cas \'el\'ementaire
o\`u $F=K(G,M,n)$, avec $G$ un sch\'ema en groupes affine, et $M$ un sch\'ema en groupe affine et unipotent.
Mais dans ce cas on a une \'equivalence de pr\'efaisceaux simpliciaux (mais pas de $H_{\infty}$-champs)
$\mathbb{R}\Omega_{*}K(G,M,n)\simeq G\times K(M,n-1)$. Comme le faisceau $G$ est repr\'esent\'e par un sch\'ema affine c'est un champ affine.
De plus, $M$ \'etant un sch\'ema en groupes affine et unipotent on sait d'apr\`es le th\'eor\`eme \ref{t5} que
$K(M,n)$ est un champ affine. Le champ $F$ est ainsi un produit de champs affines, et est donc lui-m\^eme affine.

Enfin, pour voir que $F$ est $P$-local, on utilise que les champs $P$-locaux sont stables par les $\mathbb{V}$-limites homotopiques.
En utilisant le m\^eme argument que pour la preuve du th\'eor\`eme \ref{t5} on peut supposer que $F$ est de la forme $K(G,M,n)$,
avec $G$ et $M$ de type fini sur $k$. On applique alors les m\^emes technique que lors de la preuve du th\'eor\`eme
\ref{t4} (en particulier les lemmes \ref{lt1} et \ref{lt2}) pour se ramener au cas o\`u $M$ est une repr\'esentation lin\'eaire
de dimension finie de $G$. Mais dans ce cas, $K(G,M,n)$ est tautologiquement $P$-local. \hfill $\Box$ \\

La d\'efinition suivante est une adapation de la d\'efinition de champs tr\`es pr\'esentables de \cite{s3,s4}. Elle diff\`ere un peu
de la d\'efinition originale dans le sens o\`u l'on ne consid\`ere que des pr\'efaisceaux simpliciaux point\'es et connexes, et o\`u
l'on ne se restreint pas \`a la caract\'eristique nulle.

\begin{df}\label{d11}
\emph{Un champ tr\`es pr\'esentable (point\'e et connexe) est un pr\'efaisceau simplicial point\'e et connexe $* \longrightarrow F$, satisfaisant les trois
conditions suivantes.}
\begin{itemize}
\item \emph{Il existe un entier $n$ tel que le pr\'efaisceau simplicial $F$ soit $n$-tronqu\'e (i.e. $F$ est \'equivalent \`a
$\tau_{\leq n}F$).}
\item \emph{Le faisceau en groupes $\pi_{1}(F,*)$ est repr\'esent\'e par un sch\'ema 
en groupes affine et de type fini sur $k$.}
\item \emph{Pour tout entier $i>0$ le faisceau en groupes $\pi_{i}(F,*)$ est repr\'esent\'e par un sch\'ema en groupes affine
unipotent et de type fini sur $k$.}
\end{itemize}
\end{df}

Par d\'efinition un champ tr\`es pr\'esentable est un objet point\'e et sera donc toujours consid\'er\'e comme un objet dans
$SPr(k)_{*}$. \\

On tire imm\'ediatement du th\'eor\`eme \ref{t6} le corollaire important suivant.

\begin{cor}
Tout champ tr\`es pr\'esentable est un type d'homotopie sch\'ematique point\'e.
\end{cor}

Le th\'eor\`eme \ref{t6} permet de d\'emontrer le corollaire suivant. Pour cela, rappelons qu'un paragraphe pr\'ec\'edent
nous avons associ\'e \`a tout sch\'ema en groupes affine simplicial $G$ une $\infty$-gerbe
affine $B\mathbb{R}Spec\, \mathcal{O}(G)$.

\begin{cor}\label{ct6'}
Soit $G$ un sch\'ema en groupes affine simplicial et consid\'erons 
$F:=B\mathbb{R}Spec\, \mathcal{O}(G)$  l'$\infty$-gerbe affine point\'ee
qui lui est associ\'ee. Alors, le faisceau $\pi_{i}(F,*)$ est re\-pr\'e\-sen\-ta\-ble par un sch\'ema
en groupes affine, qui est de plus unipotent d\`es que $i>1$.

En particulier, pour tout sch\'ema en groupes affine simplicial $G$,
le champ point\'e $B\mathbb{R}Spec\, \mathcal{O}(G)$ est un type d'homotopie sch\'ematique point\'e.
\end{cor}

\textit{Preuve:} Notons $A:=\mathbb{L}\mathcal{O}(\Omega_{*}F)\simeq \mathcal{O}(G)$
la $k$-alg\`ebre co-simpliciale de cohomologie
de $\Omega_{*}F\simeq G$. Comme $F$ est une $\infty$-gerbe affine, on sait que $A$ est une $H_{\infty}$-alg\`ebre de Hopf
(voir Def. \ref{hopf}). Ceci implique en particulier que $H^{*}(A)$ est une $k$-alg\`ebre de Hopf gradu\'ee.

Notons $B$ la $k$-alg\`ebre co-simpliciale d\'efinie par le diagramme homotopiquement co-cart\'esien suivant
$$\xymatrix{
H^{0}(A) \ar[r] \ar[d] & A \ar[d] \\
k \ar[r] & B,}$$
o\`u $H^{0}(A) \longrightarrow A$ est le morphisme naturel, et $H^{0}(A) \longrightarrow k$ l'augmentation correspondant
\`a la co-unit\'e de $H^{*}(A)$.

\begin{lem}\label{lconj}
Pour tout entier $n$, $H^{n}(A)$ est un $H^{0}(A)$-module libre. De plus,
on a $H^{n}(B)\simeq H^{n}(A)\otimes_{H^{0}(A)}k$.
\end{lem}

\textit{Preuve:} La seconde assertion est clairement impliqu\'ee par la premi\`ere.

Notons $G=Spec\, H^{0}(A)$ le sch\'ema en groupes affine correspondant \`a la $k$-alg\`ebre de Hopf $H^{0}(A)$.
Pour tout $n$, $H^{n}(A)$ est naturellement un $H^{0}(A)$-module, et correspond \`a un faisceau quasi-coh\'erent
sur $G$, dont on cherche \`a montrer qu'il est libre.
La structure d'alg\`ebre de Hopf induit sur $H^{n}(A)$
une structure de faisceau quasi-coh\'erent $G$-\'equivariant sur $G$ (pour l'action de $G$ sur lui-m\^eme par translations).
Ainsi, le faisceau quasi-coh\'erent $H^{n}(A)$ est trivial sur $G$ (i.e. de la forme $V\otimes_{k} \mathcal{O}_{G}$,
o\`u $V$ est un $k$-espace vectoriel), et donc $H^{n}(A)$ est un $H^{0}(A)$-module libre. \hfill $\Box$ \\

Consid\'erons le diagramme homotopiquement cart\'esien de champs
$$\xymatrix{
\mathbb{R}Spec\, B \ar[r] \ar[d] & \mathbb{R}Spec\, A \ar[d] \\
Spec\, k \ar[r] & Spec\, H^{0}(A).}$$
Comme le lemme pr\'ec\'edent montre que $B$ est une $k$-alg\`ebre co-simpliciale augment\'ee et connexe, le
th\'eor\`eme \ref{t5'} et la suite exacte longue en homotopie impliquent que
les faisceaux $\pi_{i}(F,*)$ sont des sch\'emas en groupes unipotents pour $i>1$. De plus, le
morphisme $\pi_{1}(F,*) \longrightarrow Spec\, H^{0}(A)$ est alors un monomorphisme. Il nous reste donc
\`a montrer que c'est un \'epimorphisme de faisceaux.

Pour cela on consid\`ere les morphismes naturels 
$$G \longrightarrow \mathbb{R}Spec\, A \longrightarrow Spec\, H^{0}(A).$$ 
Par construction $Spec\, H^{0}(A)\simeq \pi_{0}(G)$ est le sch\'ema en groupes
\'equaliseur des morphismes $d_{1},d_{0} : G_{1} \longrightarrow G_{0}$. Ainsi, $G_{0} \longrightarrow
Spec\, H^{0}(A)$ est un quotient de sch\'emas en groupes affines et donc est un morphisme
fid\`element plat de sch\'emas affines. Ceci montre que le morphisme $\mathbb{R}Spec, A \longrightarrow
Spec, H^{0}(A)$ poss\`ede une section apr\`es le changement de bases fid\`element plat $G_{0} \longrightarrow
Spec\, H^{0}(A)$. Ainsi, le morphisme de faisceaux
$$\pi_{1}(F,*)\simeq \pi_{0}(\mathbb{R}Spec\, A) \longrightarrow Spec\, H^{0}(A)$$
est un \'epimorphisme. \hfill $\Box$ \\

Le th\'eor\`eme \ref{t6} poss\`ede aussi une r\'eciproque.

\begin{thm}\label{t6'}
\begin{enumerate}
\item Soit $F$ un type d'homotopie sch\'ematique point\'e. Alors pour tout
$i>0$ le faisceau $\pi_{i}(F,*)$ est repr\'esentable par un sch\'ema en groupes
affine qui est de plus unipotent pour $i>1$.
\item
Soit $F$ une $\infty$-gerbe affine point\'ee. Alors les deux assertions suivantes sont
v\'erifi\'ees.
\begin{itemize}
\item Pour tout $i>1$, le faisceau $\pi_{i}(F,*)$ est repr\'esentable par un sch\'ema en groupes affine et unipotent.

\item Le faisceau $\pi_{1}(F,*)$ est un sous-faisceau d'un faisceau repr\'esentable par un sch\'ema en groupes
affine.
\end{itemize}
\end{enumerate}
\end{thm}

\textit{Preuve:} Le point $(1)$ ne sera pas d\'emontr\'e dans ce travail. Nous
renvoyons \`a \cite{small} pour une preuve. \\

D\'emontrons la seconde assertion.
Soit $A$ une $H_{\infty}$-alg\`ebre de Hopf telle que $F\simeq B\mathbb{R}Spec\, A$.
Il existe un morphisme naturel de $H_{\infty}$-alg\`ebres de Hopf
$$H^{0}(A) \longrightarrow A$$
qui induit un morphisme bien d\'efini $\mathbb{R}\Omega_{*}F \longrightarrow Spec\, H^{0}(A)$. Le lemme \ref{lconj}
implique
alors que la fibre homotopique de ce morphisme est de la forme $\mathbb{R}Spec\, B$, o\`u $B$ est une
$k$-alg\`ebre co-simpliciale augment\'ee et cohomologiquement connexe. Le th\'eor\`eme
\ref{t5'} et la suite longue des faisceaux d'homotopie permet alors de conlcure. \hfill $\Box$ \\

Nous terminerons ce paragraphe par la conjecture suivante.

\begin{conj}\label{conj}
Toute $\infty$-gerbe affine point\'ee est un type d'homotopie sch\'e\-ma\-ti\-que
point\'e.
\end{conj}

\subsection{Sch\'ematisation des types d'homotopie}

Pour tout $\mathbb{U}$-ensemble simplicial point\'e $X$, nous noterons encore $X \in SPr(k)_{*}$ le pr\'efaisceau simplicial
point\'e
constant qu'il d\'efinit. Nous supposerons toujours que de tels ensembles simpliciaux
sont point\'es connexes. Il s'ensuit que les pr\'efaisceaux simpliciaux $X$ seront donc point\'es et connexes.
Nous identifierons alors la cat\'egorie des $\mathbb{U}$-ensembles simpliciaux point\'es comme la
sous-cat\'egorie pleine de $SPr(k)_{*}$ form\'ee des pr\'efaisceaux constants.
Ceci nous permet d'appliquer les d\'efinitions de
$\mathcal{O}$-\'equivalences et de $P$-\'equivalences aux morphismes d'ensembles simpliciaux point\'es.

\begin{df}\label{d13}
\emph{Soit $X$ un ensemble simplicial point\'e et connexe de $\mathbb{U}$.
Une sch\'ematisation
de $X$ sur $k$ est la donn\'ee d'un type d'homotopie sch\'ematique point\'e $(X\otimes k)^{sch}$,
et d'un morphisme dans la cat\'egorie homotopique $Ho(SPr(k)_{*})$, $u : X\longrightarrow (X\otimes k)^{sch}$,
qui soit universel pour les morphismes vers des types d'homotopie sch\'ematiques point\'es.}
\end{df}

Remarquer que si une sch\'ematisation existe alors elle est unique \`a isomorphisme unique pr\`es dans la cat\'egorie
homotopique. De m\^eme, tout morphisme $X \longrightarrow Y$ induira automatiquement un morphisme
sur les sch\'ematisations (si elles existent) $(X\otimes k)^{sch} \longrightarrow (Y\otimes k)^{sch}$.
Ces propri\'et\'es seront tr\`es pratiques
par la suite car elles permettent de n\'egliger les aspects fonctoriels des constructions.

Remarquer aussi qu'une sch\'ematisation est aussi une $P$-localisation (i.e. un morphisme universel vers les
objets $P$-locaux). Ainsi, pour d\'emontrer que $(X\otimes k)^{sch}$ existe il suffit de trouver un
type d'homotopie sch\'ematique $(X\otimes k)^{sch}$ est une $P$-\'equivalence $X \longrightarrow (X\otimes k)^{sch}$. \\

Le lemme suivant nous sera utile pour d\'eterminer les $P$-\'equivalences entre ensembles simpliciaux.

\begin{lem}\label{l14}
Un morphisme entre deux $\mathbb{U}$-ensembles simpliciaux point\'es $f : X \longrightarrow Y$ est une $P$-\'equivalence si et
seulement si les deux conditions suivantes sont satisfaites.
\begin{enumerate}
\item Le morphisme de groupes $f_{*} : \pi_{1}(X,*) \longrightarrow \pi_{1}(Y,*)$ induit un isomorphisme sur les
compl\'et\'es affines
$$f_{*} : \mathcal{A}(\pi_{1}(X,*)) \simeq \mathcal{A}(\pi_{1}(Y,*)).$$

\item Pour tout syst\`eme local de $k$-espaces vectoriels de dimension finie $L$ sur $Y$, et tout entier $i\geq 0$, le morphisme induit
$$f^{*} : H^{i}(Y,L) \longrightarrow H^{i}(X,f^{*}L)$$
est un isomorphisme.

\end{enumerate}
\end{lem}

\textit{Preuve:} En appliquant directement la d\'efinition il est facile de voir qu'il suffit de d\'emontrer que pour tout sch\'ema en groupes
affine et de type fini
$G$, $V$ une de ses repr\'esentations lin\'eaires de dimension finie, et tout $\mathbb{U}$-ensemble simplicial point\'e et connexe $X$, il existe
des isomorphismes naturels
$$[X,K(G,V,n)]_{SPr(k)_{*}}\simeq [X,K(G(k),V(k),n)]_{\mathbb{U}-SEns_{*}}.$$

D'apr\`es le th\'eor\`eme \ref{t2} on a
$$[X,K(G,V,n)]_{SPr(k)_{*}}\simeq
[\mathbb{R}\Omega_{*}X,G\times_{\rho}K(V,n-1)]_{Pr-\Delta^{o}-SPr(k)},$$
o\`u $G\times_{\rho}K(V,n-1)$ est le pr\'efaisceau en groupes simpliciaux produit semi-direct de $G$ par
$K(V,n-1)$. L'adjonction de Quillen
$$Cst : Pr-\Delta^{o}-SEns \longrightarrow Pr-\Delta^{o}-SPr(k)_{*}$$
$$\Gamma : Pr-\Delta^{o}-SPr(k)_{*}  \longrightarrow Pr-\Delta^{o}-SEns$$
donne
$$[\mathbb{R}\Omega_{*}X,G\times_{\rho}K(V,n-1)]_{Pr-\Delta^{o}-SPr(k)}\simeq
[\mathbb{R}\Omega_{*}X,\mathbb{R}\Gamma(G\times_{\rho}K(V,n-1))]_{Pr-\Delta^{o}-SEns}.$$
Cependant, le pr\'efaisceau simplicial sous-jacent \`a $G\times_{\rho}K(V,n-1)$ est $G\times K(V,n-1)$.
De plus, $V$ \'etant un faisceau en groupes sous-jacent \`a un faisceau quasi-coh\'erent il est acyclique
sur $(Aff/k)_{fpqc}$. Ainsi, $G\times K(V,n-1)$ est un champ, et donc $G\times_{\rho}K(V,n-1)$ est un objet
fibrant dans $Pr-\Delta^{o}-SPr(k)$. Ceci implique que
$$\mathbb{R}\Gamma(G\times_{\rho}K(V,n-1))\simeq G(k)\times_{\rho}K(V(k),n-1).$$
On en d\'eduit donc
$$[X,K(G,V,n)]_{SPr(k)_{*}}\simeq
[\mathbb{R}\Omega_{*}X,G(k)\times_{\rho}K(V(k),n-1)]_{Pr-\Delta^{o}-SEns}$$
$$\simeq [X,K(G(k),V(k),n)]_{\mathbb{U}-SEns_{*}},$$
ce qu'il fallait d\'emontrer. \hfill $\Box$ \\

\begin{cor}\label{c6}
Soit $X$ un $\mathbb{U}$-ensemble simplicial point\'e et connexe, $G$ un sch\'ema en groupe affine et $V$ une repr\'esentation
lin\'eaire de $G$. Alors, il existe un isomorphisme naturel
$$\mathbb{R}\underline{Hom}_{*}(X,K(G,V,n))\simeq \mathbb{R}\underline{Hom}_{*}(X,K(G(k),V(k),n)).$$
\end{cor}

\textit{Preuve:} C'est en r\'ealit\'e un corollaire de la preuve du lemme \ref{l14}. \hfill $\Box$ \\

Le th\'eor\`eme principal de ce paragraphe est le suivant.

\begin{thm}\label{t7}
Tout $\mathbb{U}$-ensemble simplicial point\'e et connexe $(X,x)$ poss\`ede une sch\'ematisation $(X\otimes k)^{sch}$.
\end{thm}

\textit{Preuve:} Nous allons donner une construction explicite du type d'homotopie sch\'ematique point\'e
$(X\otimes k)^{sch}$, et nous montrerons que le morphisme $X \longrightarrow (X\otimes k)^{sch}$ est une $P$-\'equivalence.
Comme
les types d'homotopie sch\'ematiques point\'es sont des objets $P$-locaux,
ceci impliquera que ce morphisme est une sch\'ematisation. \\

Consid\'erons le foncteur classifiant $B : \mathbb{U}-SGp \longrightarrow \mathbb{U}-SEns_{*}$, qui \`a un
$\mathbb{U}$-groupe simplicial $G$ associe son ensemble simplicial point\'e classifiant $BG$. On sait que ce foncteur induit un foncteur au niveau
des cat\'egories homotopiques $B : Ho(\mathbb{U}-SGp) \longrightarrow Ho(\mathbb{U}-SEns_{*})$, qui est pleinement fid\`ele et dont l'image essentielle
est form\'ee des ensembles simpliciaux connexes (voir \cite[Prop. $1.5$]{to1} pour une preuve de ce fait bien connue).
Il existe donc un $\mathbb{U}$-groupe simplicial $G$ tel que $BG$ soit \'equivalent \`a $X$. De plus, quitte
\`a prendre une r\'esolution libre de $G$ on pourra supposer que chaque groupe $G_{n}$ est un groupe libre de $\mathbb{U}$.

Consid\'erons alors $\mathcal{A}(G)$ le sch\'ema en groupes affine simplicial d\'eduit de $G$ en appliquant la construction
de compl\'etion affine (voir Def. \ref{d6}). On d\'efinit le champ $(X\otimes k)^{sch}$ comme \'etant
$B\mathbb{R}Spec\, \mathcal{O}(h_{\mathcal{A}(G)}) \in Ho(SPr(k)_{*})$, dont la construction est expliqu\'ee dans le paragraphe $2.3$.
Nous savons donc d\'ej\`a que $(X\otimes k)^{sch}$ est une $\infty$-gerbe affine point\'ee, et que le morphisme
naturel $B\mathcal{A}(G) \longrightarrow (X\otimes k)^{sch}$ est de plus universel pour les morphismes vers des $\infty$-gerbes
affines point\'ees (voir Lem. \ref{l9}). De plus, le corollaire \ref{ct6'} montre que $(X\otimes k)^{sch}$ est un
type d'homotopie sch\'ematique point\'e.
Il nous reste \`a montrer que le morphisme naturel $X \longrightarrow (X\otimes k)^{sch}$
est une $P$-\'equivalence. \\

Nous avons d\'efini
$(X\otimes k)^{sch}$ comme \'etant $B\mathbb{R}Spec\, \mathcal{O}(\mathcal{A}(G))$. On dispose donc des morphismes naturels
$$\xymatrix{
X\simeq BG \ar[r]^-{a} & B\mathcal{A}(G) \ar[r]^-{b} & B\mathbb{R}Spec\, \mathcal{O}(h_{\mathcal{A}(G)})=:(X\otimes k)^{sch},}$$
et il nous suffit donc de montrer que $a$ et $b$ sont des $P$-\'equivalences. Pour le morphisme $b$
c'est imm\'ediat d'apr\`es sa propri\'et\'e universelle (voir Lem. \ref{l9}), et le fait
que $K(K,V,N)$ est une $\infty$-gerbe affine point\'ee d\`es que $K$ est un sch\'ema en groupes
affine et $V$ une de ses repr\'esentations lin\'eaire de dimension finie.
Pour terminer montrons que le morphisme $a : BG \longrightarrow B\mathcal{A}(G)$ est une $P$-\'equivalence. \\

Tout d'abord en utilisant les \'equivalences 
$$BG\simeq Hocolim_{[n] \in \Delta^{o}}BG_{n} \qquad 
B\mathcal{A}(G)\simeq Hocolim_{[n] \in \Delta^{o}}B\mathcal{A}(G_{n})$$
on voit qu'il suffit de montrer que pour tout $n$
le morphisme $BG_{n} \longrightarrow B\mathcal{A}(G_{n})$ est une $P$-\'equivalence.

Soit $H$ un sch\'ema en groupes affine et de type fini sur $k$, $V \in Rep_{k}(H)$
une repr\'esentation lin\'eaire de dimension finie et $F=K(H,V,n)$. D'apr\`es le corollaire \ref{c6} il faut montrer que
$$\mathbb{R}\underline{Hom}_{*}(BG_{n},F(k))\simeq \mathbb{R}\underline{Hom}_{*}(B\mathcal{A}(G_{n}),F).$$
Or, il existe un diagramme commutatif
$$\xymatrix{
\mathbb{R}\underline{Hom}_{*}(BG_{n},F(k)) \ar[r]^{u} \ar[d]^{v} & \mathbb{R}\underline{Hom}_{*}(B\mathcal{A}(G_{n}),F) \ar[d]^{w} \\
\mathbb{R}\underline{Hom}_{*}(BG_{n},K(G(k),1)) \ar[r]_{x} & \mathbb{R}\underline{Hom}_{*}(B\mathcal{A}(G_{n}),K(G,1)).}$$
De plus, la propri\'et\'e universelle du morphisme $G_{n} \longrightarrow \mathcal{A}(G_{n})$ implique clairement que
$x$ est un isomorphisme.
Il nous suffit donc de montrer que $u$ induit des isomorphismes sur toutes les fibres homotopiques
des morphismes $v$ et $w$. Or, si $\rho : G_{n} \longrightarrow G(k)$ est un point dans $\pi_{0}(\mathbb{R}\underline{Hom}_{*}(BG_{n},K(G(k),1)))$,
et si on note $v^{-1}(\rho)$ et $w^{-1}(\rho)$ ces fibres homotopiques, on a
$$\pi_{i}(v^{-1}(\rho))\simeq H^{n-i}(BG_{n},V(k)) \qquad
\pi_{i}(w^{-1}(\rho))\simeq H^{n-i}(B\mathcal{A}(G_{n}),V).$$
On se ram\`ene donc \`a montrer le lemme suivant.

\begin{lem}\label{l15}
Pour tout entier $n$, et toute repr\'esentation lin\'eaire de dimension finie $V$ de $\mathcal{A}(G_{n})$, on a
$$H^{i}(B\mathcal{A}(G)_{n},V)\simeq H^{i}(BG_{n},V)$$
pour tout entier $i\geq 0$.
\end{lem}

\textit{Preuve:} Commen\c{c}ons par montrer que $H^{i}(B\mathcal{A}(G)_{n},V)\simeq H^{i}(BG_{n},V)\simeq 0$
pour $i>1$. Pour ce qui est de l'annulation de $H^{i}(BG_{n},V)$ cela provient imm\'ediatement du fait que $G_{n}$ soit
libre et du corollaire \ref{c6}. Ainsi,
d'apr\`es le lemme \ref{l4} il nous suffit de montrer que pour toute repr\'esentation lin\'eaire $W$ (\'eventuellement de dimensions infinie) on a
$H^{2}(B\mathcal{A}(G_{n}),W)=0$. Comme une telle repr\'esentation est limite inductive de ses sous-repr\'esentations
de dimension finie on peut \`a l'aide du corollaire \ref{c'} se restreindre au cas o\`u $W$ est de dimension finie.

Il est facile de voir que pour montrer que $H^{2}(B\mathcal{A}(G_{n}),W)=0$
il suffit de montrer que toute extension de faisceaux en groupes
$$\xymatrix{0 \ar[r] & W \ar[r] & E \ar[r]^-{p} & \mathcal{A}(G_{n}) \ar[r] & 1}$$
poss\`ede une section. Cependant, comme $E$ est un $W$-torseur sur $\mathcal{A}(G_{n})$ qui est un sch\'ema affine, le morphisme
$p$ poss\`ede une section dans la cat\'egorie des sch\'emas $s : \mathcal{A}(G_{n}) \longrightarrow E$. Remarquons aussi
que $E$ est repr\'esentable par un sch\'ema en groupes affine.
Ainsi, si le groupe $G_{n}$ est un groupe libre sur un $\mathbb{U}$-ensemble $I$, le morphisme de pr\'efaisceaux
$$\xymatrix{I \ar[r] & G_{n} \ar[r] & \mathcal{A}(G_{n}) \ar[r]^-{s} & E}$$
se rel\`eve en un morphisme de groupes $G_{n} \longrightarrow E$.
Ce morphisme induit alors un morphisme de sch\'emas en groupes $\mathcal{A}(G_{n}) \longrightarrow E$, qui est la section cherch\'ee. \\

Il nous reste donc \`a v\'erifier que $H^{i}(B\mathcal{A}(G_{n}),V)\simeq H^{i}(BG_{n},V)$ pour $i<2$. Mais ceci est vrai pour
tout groupe $G_{n}$, libre ou pas, et se d\'eduit ais\'ement de la propri\'et\'e universelle des compl\'etions affines (voir par exemple
\cite[$\S 4$]{lm}). \hfill $\Box$ \\

Ce lemme permet de conlure que le morphisme $a$ est une \'equivalence, et termine donc la preuve du th\'eor\`eme.
\hfill $\Box$ \\

Notons $Ho(THS/k)$  la sous-cat\'egorie pleine de $Ho(SPr(k)_{*})$ form\'ee
des types d'homotopie sch\'ematiques point\'es. Consid\'erons le foncteur d\'eriv\'e des sections globales
$$\mathbb{R}\Gamma : Ho(THS/k)\subset Ho(SPr_{*}(k)) \longrightarrow Ho(\mathbb{V}-SEns_{*}).$$
En se restreignant \`a la composante connexe qui contient le point distingu\'e on en d\'eduit un foncteur
$$\mathbb{R}\Gamma : Ho(THS/k) \longrightarrow Ho(\mathbb{V}-SEns_{0}),$$
o\`u $Ho(\mathbb{V}-SEns_{0})$ est la cat\'egorie homotopique des $\mathbb{V}$-ensembles simpliciaux point\'es et connexes.
En utilisant le th\'eor\`eme \ref{t2} et le corollaire \ref{c1}, on voit que le foncteur
$$\mathbb{R}\Gamma : Ho(THS/k) \longrightarrow Ho(\mathbb{V}-SEns_{0})$$
se factorise en un foncteur
$$\mathbb{R}\Gamma : Ho(THS/k) \longrightarrow Ho(\mathbb{U}-SEns_{0}).$$

\begin{cor}\label{c10}
Le foncteurs d\'eriv\'e du foncteur des sections globales
$$\mathbb{R}\Gamma : Ho(THS/k) \longrightarrow Ho(\mathbb{U}-SEns_{0})$$
poss\`ede un adjoint \`a gauche
$$(-\otimes k)^{sch} : Ho(\mathbb{U}-SEns_{0}) \longrightarrow Ho(THS/k).$$
\end{cor}

\textit{Preuve:} C'est une autre fa\c{c}on d'\'enoncer le th\'eor\`eme \ref{t7}. \hfill $\Box$ \\

\begin{cor}\label{c8}
Pour tout ensemble simplicial point\'e et connexe $(X,x)$ appartenant \`a $\mathbb{U}$, il existe un isomorphisme naturel
$$\pi_{1}((X\otimes k)^{sch},x)\simeq \mathcal{A}(\pi_{1}(X,x)).$$
\end{cor}

\textit{Preuve:} On sait d'apr\`es le th\'eor\`eme
\ref{t6'} que $\pi_{1}(X\otimes ,k,x)$ est un sch\'ema en groupes affine.
De plus, comme $X \longrightarrow (X\otimes k)^{sch}$
est une $P$-\'equivalence, \ref{t6} implique que'on a pour tout sch\'ema en groupes affine $G$
$$Hom(\pi_{1}(X,x),G) \simeq Hom(\pi_{1}((X\otimes k)^{sch},x),G).$$
On conclut alors \`a l'aide de la propri\'et\'e universelle de $\mathcal{A}(\pi_{1}(X,x))$. \hfill $\Box$ \\

\begin{cor}\label{c11}
Soit $X$ un $\mathbb{U}$-ensemble simplicial point\'e et connexe. Si le champ $(X\otimes k)^{sch}$ est simplement connexe
(i.e. $\pi_{1}(X\otimes k)\simeq 0$), alors il existe un isomorphisme naturel
$$(X\otimes k)^{sch} \longrightarrow (X\otimes k)^{uni}.$$
\end{cor}

\textit{Preuve:} Comme $(X\otimes k)^{sch}$ est simplement connexe, c'est un champ affine d'apr\`es \ref{t5} et \ref{t7}. De plus,
$X \longrightarrow (X\otimes k)^{sch}$ \'etant une $P$-\'equivalence, c'est aussi une $\mathcal{O}$-\'equivalence. C'est donc
une affination. \hfill $\Box$ \\

Remarquer qu'en g\'en\'eral il existe toujours un morphisme naturel
$$(X\otimes k)^{sch} \longrightarrow (X\otimes k)^{uni}.$$

\begin{cor}
Soit $X$ un $\mathbb{U}$-ensemble simplicial point\'e, simplement connexe et de type fini. Alors on a
$$\pi_{i}((X\otimes k)^{sch},x)\simeq \left\{
\begin{array}{cl}
\pi_{i}(X,x)\otimes_{\mathbb{Z}}\mathbb{G}_{a} & \text{si car k=0} \\
\pi_{i}(X,x)\otimes_{\mathbb{Z}}\mathbb{Z}_{p} & \textit{si car k=p}.
\end{array} \right.$$
\end{cor}

\textit{Preuve:} Ceci se d\'eduit des corollaires \ref{c11} et \ref{ct8}. \hfill $\Box$ \\

\subsection{Exemples et contre-exemples}

Dans ce paragraphe nous donnons une liste d'exemples et de contre-exemples du comportement du foncteur de
sch\'ematisation. Nous ne donnerons pas de preuves \'etant donn\'e que les arguments
se d\'eduisent ou bien de la propri\'et\'e
universelle, ou bien des r\'esultats des deux paragraphes pr\'ec\'edents.

\begin{itemize}

\item Reprenons l'exemple cit\'e dans le paragraphe $2.6$, d'un type d'homotopie $X=K(\mathbb{Z}/m,\mathbb{Z}^{r},n)$.
Il n'est pas difficile de voir dans ce cas que
$$(X\otimes \mathbb{Q})^{sch}\simeq K(\mathbb{Z}/m,\mathbb{G}_{a}^{r},n),$$
et o\`u l'action de $\mathbb{Z}/m$ sur $\mathbb{G}_{a}^{r}$ est la $\mathbb{Q}$-lin\'earis\'ee de l'action
sur $\mathbb{Z}^{r}$.

\item Il est bien connu qu'il est g\'en\'eralement impossible de donner une formule pour les faisceaux d'homotopie du champ
$(X\otimes k)^{uni}$ lorsque $X$ n'est pas nilpotent. Il est donc aussi impossible de donner une formule g\'en\'erale pour
les faisceaux d'homotopie du champ $(X\otimes k)^{sch}$. Cela peut en r\'ealit\'e d\'ej\`a se remarquer sur le cas $1$-tronqu\'e.

Si $\Gamma$ est un $\mathbb{U}$-groupe, alors on sait que $\pi_{1}((K(\Gamma,1)\otimes k)^{sch},*)$ est le compl\'et\'e affine
de $\Gamma$ sur $k$. Par contre, il n'est pas vrai en g\'en\'eral que $(K(\Gamma,1)\otimes k)^{sch}\simeq K(\mathcal{A}(\Gamma),1)$.
La raison en est que le groupe $\Gamma$ peut tr\`es bien ne pas avoir la m\^eme cohomologie que $\mathcal{A}(\Gamma)$.
Le ph\'enom\`ene qui se produit ici est de nature tout \`a fait \'equivalente \`a ce qu'il se produit lors des compl\'etions
pro-finies (voir \cite{am}). On peut donc introduire un analogue alg\'ebrique de la notion de groupe \textit{bon}.

\begin{df}\label{d14}
\emph{Le groupe $\Gamma$ est bon sur $k$ si le morphisme naturel
$$(K(\Gamma,1)\otimes k)^{sch}\longrightarrow K(\mathcal{A}(\Gamma),1)$$
est un isomorphisme dans $Ho(SPr(k)_{*})$.}
\end{df}

Il est par exemple facile de voir que les groupes libres sont bons sur n'importe quel corps $k$ (ceci
provient du lemme \ref{l15}). De m\^eme, les groupes ab\'eliens de types finis, les groupes
fondamentaux de surfaces de Riemann et les extensions succ\'esives de groupes
libres de types finis sont tous des exemples de groupes bons (voir \cite{small}).
Par contre un groupe ab\'elien qui n'est
pas de type fini n'est g\'en\'eralement pas bon, et ce m\^eme si $k=\mathbb{Q}$.

On en d\'eduit par exemple que
$$(K(\mathbb{Z},1)\otimes k)^{sch}\simeq K(\mathcal{A}(\mathbb{Z}),1).$$
On peut de plus calculer explicitement le sch\'ema groupe affine $\mathcal{A}(\mathbb{Z})$.
Pour cela, consid\'erons $\overline{k}^{*}$, le groupe discret multiplicatif de la cl\^oture alg\'ebrique du corps $k$.
Il est muni d'une action naturelle du groupe de galois $Gal(\overline{k}/k)$. Formons alors le sch\'ema en groupes
de type multiplicatif $D(\overline{k}^{*})$, tel que son faisceau des caract\`eres soit isomorphe \`a $\overline{k}^{*}$, muni de son
action galoisienne. On peut alors montrer qu'il existe un isomorphisme (voir par exemple \cite[App. $A$]{js})
$$\mathcal{A}(\mathbb{Z})\simeq D(\overline{k}^{*})\times \mathcal{U}(\mathbb{Z}),$$
et donc que la sch\'ematisation de $K(\mathbb{Z},1)$ est donn\'ee par
$$(K(\mathbb{Z},1)\otimes k)^{sch}\simeq K(D(\overline{k}^{*}),1)\times K(\mathcal{U}(\mathbb{Z}),1).$$
Ceci implique que le morphisme naturel
$$(K(\mathbb{Z},1)\otimes k)^{sch} \longrightarrow (K(\mathbb{Z},1)\otimes k)^{uni}$$
poss\`ede une fibre homotopique \'equivalente \`a $K(D(\overline{k}^{*}),1)$, qui est un champ \'enorme. On voit ainsi, que m\^eme dans les cas
les plus simples, lorsque $X$ n'est pas simplement connexe le champ $(X\otimes k)^{sch}$ est g\'en\'eralement beaucoup plus
gros que le champ $(X\otimes k)^{uni}$. D\'ej\`a au niveau du groupe fondamental, le morphisme
$\pi_{1}((X\otimes k)^{sch},*) \longrightarrow \pi_{1}((X\otimes k)^{uni},*)$ poss\`ede un tr\`es gros noyau. Il s'ensuit \`a fortiori
que la fibre homotopique de $(X\otimes k)^{sch} \longrightarrow (X\otimes k)^{uni}$ est aussi \'enorme.

Une autre cons\'equence remarquable du fait que $\pi_{1}((X\otimes k)^{sch},*)\simeq \mathcal{A}(\pi_{1}(X,*))$ est que la formation
de la sch\'ematisation $(X\otimes k)^{sch}$ ne commute pas avec les changements de bases en g\'en\'eral. Par exemple, il n'est pas
vrai que $D(\overline{\mathbb{Q}}^{*})\times_{Spec\, \overline{\mathbb{Q}}}Spec\, \mathbb{C} \simeq
D(\mathbb{C}^{*})$. Ceci est une diff\'erence importante avec ce qui se passe pour les affinations (voir Cor. \ref{c12}). \\

\item Lorsque $k$ est de caract\'eristique nulle, la sch\'ematisation $(X\otimes k)^{sch}$ poss\`ede une d\'ecomposition
naturelle en une partie r\'eductive et une partie unipotente, analogue \`a la d\'ecomposition d'un sch\'ema en groupes
affine (voir \cite[$4$]{hm}).

Consid\'erons la projection naturelle 
$$(X\otimes k)^{sch} \longrightarrow K(\pi_{1}((X\otimes k)^{sch},*),1),$$ 
ainsi que le
quotient maximal
r\'eductif 
$$\pi_{1}((X\otimes k)^{sch},*) \longrightarrow \pi_{1}((X\otimes k)^{sch},*)^{red}.$$ 
Par d\'efinition,
$\pi_{1}((X\otimes k)^{sch},*)^{red}$
est la limite projective des quotients r\'eductifs et de type fini de $\pi_{1}((X\otimes k)^{sch},*)$, et se trouve \^etre isomorphe
\`a l'enveloppe r\'eductive de $\pi_{1}(X,*)$.
Le noyau de la projection $\pi_{1}((X\otimes k)^{sch},*)\longrightarrow
\pi_{1}((X\otimes k)^{sch},*)^{red}$
est un sch\'ema en groupe affine et unipotent.
Ainsi, la fibre homotopique du morphisme
$$(X\otimes k)^{sch} \longrightarrow K(\pi_{1}((X\otimes k)^{sch},*)^{red},1)$$
est un champ connexe dont tous les faisceaux d'homotopie sont repr\'esentables par des sch\'emas en groupes unipotents.
D'apr\`es le th\'eor\`eme \ref{t5}, c'est donc un champ affine. Il existe donc une $k$-alg\`ebre co-simpliciale $A$ dans $\mathbb{U}$, une action de
$\pi_{1}((X\otimes k)^{sch},*)^{red}$ sur le champ affine $(X\otimes k)^{o}:=\mathbb{R}Spec\, A$, tel que
$(X\otimes k)^{sch}$ soit le quotient homotopique de $(X\otimes k)^{o}$ par l'action de $\pi_{1}((X\otimes k)^{sch},*)$.
De plus, le sch\'ema en groupes $\pi_{1}((X\otimes k)^{sch},*)^{red}$ \'etant r\'eductif et $k$ de caract\'eristique nulle,
sa cohomologie \`a valeurs dans
des repr\'esentations lin\'eaires est nulle, et on a donc pour toute repr\'esentation lin\'eaire
$\rho : \pi_{1}((X\otimes k)^{sch},*) \longrightarrow Gl(V)$
$$H^{i}((X\otimes k)^{sch},V)\simeq H^{i}((X\otimes k)^{o},V)^{\pi_{1}((X\otimes k)^{sch},*)^{red}}.$$
La $k$-alg\`ebre $A$, munie de l'action de $\pi_{1}((X\otimes k)^{sch},*)^{red}$ sur $\mathbb{R}Spec\, A$, r\'epond
ainsi \`a une question pos\'ee par A. Beilinson, T. Pantev et L. Katzarkov, et qui avait \'et\'e
abord\'ee par des m\'ethodes diff\'erentes
lorsque $X$ \'etait l'espace topologique sous-jacent \`a une vari\'et\'e alg\'ebrique lisse sur $\mathbb{C}$.
Cette construction est aussi li\'ee \`a un mod\`ele explicite pour le champ $(X\otimes k)^{sch}$ d\'ecrit \`a l'aide
d'une notion de \textit{champs affines \'equivariants} (voir \cite{kt}).

\item Nous terminerons ces exemples par un exemple qui montre que bien que le champ $(X\otimes k)^{sch}$ contienne beaucoup d'information
homotopique il ne suffit pas \`a reconstruire $X$. Pour cela nous allons construire un morphisme $f : X \longrightarrow Y$
entre deux $\mathbb{U}$-ensembles simpliciaux point\'es, qui induit des \'equivalences $(X\otimes k)^{sch} \longrightarrow
(Y\otimes k)^{sch}$
pour tout $\mathbb{U}$-corps $k$, mais qui n'est pas une \'equivalence faible. Cet exemple m'a \'et\'e communiqu\'e par C. Simpson.

D'apr\`es \cite{bm} il existe un groupe $\Gamma$ qui est de pr\'esentation finie, infini et simple. Posons alors
$X=K(\Gamma,1)$ et prenons pour $Y$ la construction $+$ de Quillen $Y:=X^{+}$, muni du morphisme naturel
$f : X \longrightarrow Y$. L'ensemble simplicial $Y$ est donc simplement connexe et le morphisme $f$ est une \'equivalence
de cohomologie (\`a coefficients dans $\mathbb{Z}$ et donc dans tout corps $k$). Remarquer d\`es \`a pr\'esent que
$f$ n'est pas une \'equivalence.
D'apr\`es le corollaire \ref{c8}, on
a pour tout corps $k$
$$\pi_{1}((X\otimes k)^{sch})\simeq \mathcal{A}(\pi_{1}(X,x))\simeq \mathcal{A}(\Gamma).$$
Or, le groupe $\Gamma$ \'etant simple, infini et de pr\'esentation fini il ne
poss\`ede aucune repr\'esentation lin\'eaire sur
$k$ qui soit
non triviale (car tout sous-groupe de pr\'esentation finie de $Gl_{n}(k)$ poss\`ede des sous-groupes
d'indices finis non triviaux). Le champ $(X\otimes k)^{sch}$ est donc simplement connexe, et on a donc $(X\otimes k)^{sch}
\simeq (X\otimes k)^{uni}$
(voir Cor. \ref{c11}).
Or, comme $f$ est une \'equivalence de cohomologie le morphisme induit
$$(X\otimes k)^{sch}\simeq (X\otimes k)^{uni} \longrightarrow (Y\otimes k)^{sch} \simeq (Y\otimes k)^{uni}$$
est une \'equivalence.

\end{itemize}

\subsection{Types d'homotopie des vari\'et\'es alg\'ebriques}

Le but de ce paragraphe est d'\'esquisser quelques applications
de la notion de types d'homotopie sch\'ematiques dans le cadre de la g\'eom\'etrie
alg\'ebrique. Nous ne donnerons pas de d\'etails ni de preuves qui apparaitront dans
des travaux ult\'erieurs. Notons au passage
que les constructions ci-dessous peuvent \^etre facilement obtenues
\`a l'aide des techniques utilis\'ees dans \cite{kt,ol}, bas\'ees sur la notion
de \textit{champs affine \'equivariants}. Elles peuvent \^etre aussi obtenues, tout
au moins conjecturalement,
en appliquant le formalisme des cat\'egories de Segal Tannakiennes (voir
\cite{to2} et \cite{to3} pour plus de d\'etails). \\

On rappelle que $Ho(THS/k)$ est le sous-cat\'egorie pleine de
$Ho(SPr_{*}(k))$ form\'ee des types d'homotopie sch\'ematiques
point\'es sur le corps $k$.

\subsubsection{Type d'homotopie sch\'ematiques et th\'eorie de Hodge}

Soit $X$ une vari\'et\'e lisse, connexe et projective sur $\mathbb{C}$, et $x \in X$.
On peut associer \`a $X$ plusieurs types d'homotopie sch\'ematiques
sur $\mathbb{C}$, chacun d'eux associ\'es aux th\'eories cohomologiques
de Betti, de de Rham et de Dolbeault. La th\'eorie de Hodge non-ab\'elienne
de \cite{s1} et la correspondence de Riemann-Hilbert permettent alors
de construire des isomorphismes de comparaison entre ces types d'homotopie, et
de construire une certaine d\'ecomposition de Hodge qui englobe les
d\'ecompositions de Hodge usuelles sur la cohomologie, l'homotopie
rationnelle et le groupe fondamental. 

\begin{itemize}

\item Notons $X^{top}$ l'espace topologique des points complexes de $X$ muni
de la topologie analytique. On dispose alors
du type d'homotopie sch\'ematique point\'e $(X^{top}\otimes\mathbb{C})^{sch}$, sch\'ematisation
du type d'homotopie de $X^{top}$. La propri\'et\'e universelle de la sch\'ematisation
nous dit que la cohomologie \`a coefficients locaux de $(X^{top}\otimes\mathbb{C})^{sch}$
calcule la cohomologie de Betti de $X$ (i.e. la cohomologie de l'espace $X^{top}$
\`a valeurs dans des syst\`emes locaux de $\mathbb{C}$-vectoriels de dimension finis).

\item On associe \`a $X$ un pr\'efaisceau $X_{DR}$ d\'efini par la formule
$$X_{DR}(A):=X(A_{red}),$$
pour toute $\mathbb{C}$-alg\`ebre $A$, et o\`u $A_{red}$ est la
$\mathbb{C}$-alg\`ebre r\'eduite (voir \cite{s3}). On consid\`ere
$X_{DR}$ comme un objet de $Ho(SPr_{*}(\mathbb{C}))$, et donc comme un champ
point\'e sur $\mathbb{C}$. On peut alors montrer (par exemple en utilisant les
techniques utilis\'ees dans \cite{kt}) que le foncteur
$$\begin{array}{ccc}
Ho(THS/\mathbb{C}) & \longrightarrow & Ens \\
 F & \mapsto & [X_{DR},F]_{Ho(SPr_{*}(\mathbb{C}))}
\end{array}$$
est co-repr\'esentable par un type d'homotopie sch\'ematique point\'e, que l'on 
notera $(X,x)^{DR}$, et qui est appel\'e \textit{le type d'homotopie sch\'ematique
de de Rham de $(X,x)$}.

Par d\'efinition
il existe une \'equivalence tensorielle entre la cat\'egorie tensorielle
des repr\'esentations lin\'eaires de dimension finies de
$\pi_{1}((X,x)^{DR},*)$ et celle des fibr\'es vectoriels alg\'ebriques munis de connexions
int\'egrables sur $X$. En d'autres termes, le sch\'ema en groupes
$\pi_{1}((X,x)^{DR},*)$ est le dual de Tannaka de la cat\'egorie
des fibr\'es plats sur $X$.
De plus, pour une telle repr\'esentation $V$, la cohomologie
du champ $(X,x)^{DR}$ \`a valeurs dans $V$ s'identifie naturellement avec
la cohomologie de de Rham alg\'ebrique de $X$ \`a valeurs dans la connexion
correspondante.
$$H^{*}((X,x)^{DR},V)\simeq H^{*}_{DR}(X,V).$$

\item La correspondence de Riemann-Hilbert interpr\'et\'ee de fa\c{c}on ad\'equate permet de construire
un isomorphisme naturel de champs point\'es
$$\rho_{RH} : (X^{top}\otimes \mathbb{C})^{sch} \simeq (X,x)^{DR}.$$
En d'autres termes, il existe un morphisme naturel de champs point\'es
$$(X,x) \longrightarrow (X,x)^{DR},$$
que l'on pourrait appel\'e \textit{application des p\'eriodes}, et qui
induit par propri\'et\'e universelle de la sch\'ematisation l'isomorphisme $\rho_{RH}$.
Cette \'equivalence est une g\'en\'eralisation de l'isomorphisme
de comparaison bien connu entre cohomologie de de Rham et cohomologie de Betti au cas
des types d'homotopie.

\item \`A $X$, on peut aussi associer le champ $X_{Dol}$, qui par d\'efinition
est le champ classifiant du $X$-sch\'ema en groupes formel $\hat{TX}$, compl\'et\'e
formel du fibr\'e tangent (vu comme sch\'ema en groupes sur $X$, voir
\cite{s3}). On consid\'ere
$X_{Dol}$ comme un champ point\'e (en $x\in X(\mathbb{C})$), et donc comme
un objet de la cat\'egorie $Ho(SPr_{*}(\mathbb{C}))$. Comme pr\'ec\'edemment, on d\'efinit un foncteur
$$\begin{array}{ccc}
Ho(THS/\mathbb{C}) & \longrightarrow & Ens \\
 F & \mapsto & [X_{Dol},F]_{Ho(SPr_{*}(\mathbb{C}))},
\end{array}$$
et l'on consid\`ere un sous-foncteur d\'efini comme suit. \`A un
morphism de champs point\'es $X_{Dol} \longrightarrow F$, o\`u $F$
est un type d'homotopie sch\'ematique point\'e, on peut associer
le morphism induit par composition avec la projection $F
\longrightarrow \tau_{\leq 1}F\simeq K(H,1)$, o\`u $H=\pi_{1}(F,*)$.
Ce morphisme correspond \`a une classe d'isomorphisme de
$H$-torseur sur le champ $X_{Dol}$, o\`u de mani\`ere
\'equivalente \`a un \textit{$H$-fibr\'e principal de Higgs sur $X$}
(comme d\'efini dans \cite{s6}), que nous noterons $P$. Pour tout repr\'esentation
lin\'eaire de dimension finie $V$ de $H$, ce $H$-fibr\'e principal de Higgs
donne lieu \`a un fibr\'e de Higgs $P\times^{H}V$. Nous dirons alors que le morphisme
$X_{Dol} \longrightarrow F$ est \textit{semi-stable de degr\'e $0$}
si pour tout repr\'esentation lin\'eaire de dimension finie $V$ le fibr\'e
de Higgs $P\times^{H}V$ est semi-stable et de degr\'e $0$ au sens
de \cite{s6}. Si l'on note par $[X_{Dol},F]^{ss,0}_{Ho(SPr_{*}(\mathbb{C}))}$
le sous-ensemble de $[X_{Dol},F]_{Ho(SPr_{*}(\mathbb{C}))}$ form\'e des morphismes
semi-stables de degr\'e $0$ on obtient ainsi un foncteur
$$\begin{array}{ccc}
Ho(THS/\mathbb{C}) & \longrightarrow & Set \\
 F & \mapsto & [X_{Dol},F]^{ss,0}_{Ho(SPr_{*}(\mathbb{C}))}.
\end{array}$$
On montre alors (par exemple \`a l'aide des techniques de constructions
de \cite{kt}) que ce foncteur est co-repr\'esentable par un type
d'homotopie sch\'ematique point\'e $(X,x)^{Dol}$, et appel\'e le
\textit{type d'homotopie de Dolbeault de $X$}.

\item La correspondence fondamentale de la th\'eorie de Hodge non-ab\'elienne, telle
qu'\'enonc\'ee dans \cite[Lem. 2.2]{s6}, permet de construire un isomorphisme naturel
de champs point\'es
$$\rho_{Hod} : (X,x)^{Dol} \simeq (X,x)^{DR}.$$
Il existe donc un diagramme d'isomorphismes naturels de champs point\'es
$$\xymatrix{(X,x)^{Dol} \ar[r]^-{\rho_{Hod}} & (X,x)^{DR} & \ar[l]_-{\rho_{RH}} (X^{top}\otimes\mathbb{C})^{sch}.}$$
Le groupe (discret) $\mathbb{C}^{\times}$ op\`ere par
homoth\'eties sur le $X$-sch\'ema en groupes $TX$, et donc par
fonctorialit\'e sur le champ $X_{Dol}$. Ceci induit, \`a travers
les isomorphismes ci-dessus une action de $\mathbb{C}^{\times}$ sur le
champ $(X^{top}\otimes \mathbb{C})^{sch}$. Cette action n'est pas uniquement une action
dans la cat\'egorie homotopique des champs point\'es, mais existe de fac\c{c}on naturelle
en tant qu'objet de la cat\'egorie homotopique des pr\'efaisceaux simpliciaux point\'es
$\mathbb{C}^{\times}$-\'equivariants (voir \cite{kt} pour plus de d\'etails). Cette action
de $\mathbb{C}^{\times}$ sur $(X^{top}\otimes \mathbb{C})^{sch}$ est par d\'efinition
la \textit{d\'ecomposition de Hodge de $(X^{top}\otimes \mathbb{C})^{sch}$}. Cette terminologie
est justifi\'ee par le th\'eor\`eme principal de \cite{kt} qui affirme que
l'action de $\mathbb{C}^{\times}$ permet de retrouver les structures de Hodge
sur la cohomologie, l'homotopie rationnelle et le groupe fondamental
de $X$.

\item Le champ $(X^{top}\otimes \mathbb{C})^{sch}$, muni de sa d\'ecomposition
de Hodge me semble un invariant digne d'int\'er\^et. Par exemple, dans
\cite{kt} on montre comment on peut obtenir des exemples de types d'homotopie
qui ne sont pas r\'ealisables par des vari\'et\'es complexes lisses et projectives, et dont
l'obstruction \`a la r\'ealisabilit\'e se trouve dans des invariants d'homotopie
sup\'erieure (essentiellement l'action du groupe fondamental sur les groupes
d'homotopie). Une question qui reste \`a \'etudier est le probl\`eme de type
Torelli, c'est \`a dire la reconstruction de certain type de vari\'et\'es
alg\'ebriques \`a partir du champ $(X^{top}\otimes \mathbb{C})^{sch}$
muni de sa d\'ecomposition de Hodge.

\item Pour finir avec la th\'eorie de Hodge, signalons qu'il existe aussi une variante
de la construction pr\'ec\'edente, qui \`a $(X,x)$ associe
un \textit{type d'homotopie de Hodge} $X^{\mathbb{H}}$. Il s'agit d'un type
d'homotopie sch\'ematique point\'e sur $\mathbb{C}$, muni d'un morphisme naturel
$$X^{\mathbb{H}} \longrightarrow K(\mathbb{G}_{m},1).$$
La fibre homotopique de ce morphisme est not\'ee $\overline{X}^{\mathbb{H}}$,
et appel\'e le \textit{type d'homotopie de Hodge g\'eom\'etrique} de $X$.
Les syst\`emes locaux sur $\overline{X}^{\mathbb{H}}$ correspondent
aux syst\`emes locaux qui peuvent \^etre munis d'une
structure de variation de structures de Hodge complexes polarizabes sur $X$. De plus
la cohomologie de $\overline{X}^{\mathbb{H}}$ \`a coefficients dans un tel
syst\`eme local s'identifie \`a la cohomologie de $X$ \`a valeurs dans la
variation de structure de Hodge correspondante. En particulier, pour un tel
syt\`eme local, l'action de
$\mathbb{G}_{m}$ induite sur $H^{n}(\overline{X}^{\mathbb{H}},L)$
correspond \`a la d\'ecomposition de Hodge
$$H^{n}(X,L)\simeq \bigoplus_{p+q=n}H^{p}(X,L\otimes\Omega^{q}),$$
en op\'erant par poids $q$ sur la composante $H^{p}(X,L\otimes\Omega^{q})$.
Il existe de plus, un morphisme naturel de types d'homotopie sch\'ematiques
$$(X\otimes \mathbb{C})^{sch} \longrightarrow \overline{X}^{\mathbb{H}},$$
qui est $\mathbb{C}^{\times}$-\'equivariant. D'une certain fa\c{c}on, 
ce morphisme est le quotient maximal de $(X\otimes \mathbb{C})^{sch}$
sur lequel le groupe $\mathbb{C}^{\times}$ op\`ere alg\'ebriquement.

Il existe aussi des versions plus fines de $X^{\mathbb{H}}$ correspondants aux
variations de structures de Hodge r\'eelles, voir rationnelles. Les types
d'homotopie sch\'ematiques sont alors d\'efinis sur $\mathbb{R}$ ou sur
$\mathbb{Q}$. Le groupe $\mathbb{G}_{m}$ est alors remplac\'e par
le groupe de Tannaka de la cat\'egorie des structures de Hodge
r\'eelles ou rationnelles. Les syst\`emes locaux sur
$X^{\mathbb{H}}$ correspondent alors aux variations de structures de Hodge
(d\'efinies sur $\mathbb{R}$ ou $\mathbb{Q}$) et sa cohomologie
calcule la cohomologie absolue au sens de A. Beilinson et P. Deligne.
Il doit aussi exister des versions Hodge mixtes de ces constructions, controlant
les variations de structures de Hodge mixtes et leur cohomologie (y compris
pour des vari\'et\'es ouvertes singuli\`eres).

\end{itemize}

\subsubsection{Type d'homotopie (iso-)cristallin}

Dans \cite{ol}, M. Olsson construit des analogues cristallins des constructions
pr\'esent\'ees dans le paragraphe pr\'ec\'edent.
Dans le court paragraphe ci-dessous je r\'esume la situation. Il va sans dire
que les r\'esultats de ce paragraphe sont tous dus \`a M. Olsson. \\

Dans \cite{ol}, \`a une vari\'et\'e lisse et projective $X$ sur un corps alg\'ebriquement clos $k$ de caract\'eristique
$p>0$, une sous-cat\'egorie Tannakienne $\mathcal{C}$ de la cat\'egorie $Isoc(X)$
des isocristaux sur $X$, et un point $x\in X(k)$, M. Olsson associe un type
d'homotopie sch\'ematique point\'e $X_{\mathcal{C}}$ sur $K=Frac(W(k))$ (le corps
des fractions de l'anneau des vecteurs de Witt sur $k$). Ce type d'homotopie sch\'ematique
poss\`ede une cat\'egorie des syst\`emes locaux en $K$-espaces vectoriels de dimension finie
\'equivalente \`a la cat\'egorie $\mathcal{C}$, et sa cohomologie \`a coefficients
dans un tel syst\`eme local est naturellement isomorphe \`a la cohomologie cristalline
de $X$ \`a coefficients dans l'isocristal correspondant. Lorsque
$\mathcal{C}=Isoc(X)$, le champ $X_{\mathcal{C}}$ est not\'e
$(X,x)^{isoc}$, appel\'e le \textit{type d'homotopie (iso-)cristallin
de $X$}, et est l'analogue cristallin du type d'homotopie de de Rham
pr\'esent\'e pr\'ecedemment. Un autre exemple int\'eressant est celui
o\`u $\mathcal{C}$ est engendr\'ee, comme cat\'egorie Tannakienne sur $K$, par
les isocristaux provenant de $F$-isocrystaux. Dans ce cas,
le champ $X_{\mathcal{C}}$ est not\'e $(X,x)^{F-isoc}$ et appel\'e
le \textit{type d'homotopie $F$-(iso-)cristallin de $X$}. \\

Il existe un isomorphisme naturel
$F^{*}((X,x)^{isoc}) \simeq (X,x)^{isoc}$, o\`u $F^{*}$ d\'esigne
l'image r\'eciproque par le Frobenius. Cet
isomorphisme est l'analogue cristallin de l'action
de $\mathbb{C}^{\times}$ sur la sch\'ematisation
d'une vari\'et\'e complexe pr\'esent\'ee pr\'e\-c\'e\-dem\-ment. 
Cependant, tout comme dans le cas complexe l'action de $\mathbb{C}^{\times}$
n'est pas une action alg\'ebrique, l'isomorphisme
$F^{*}((X,x)^{isoc}) \simeq (X,x)^{isoc}$ ne peut pas s'interpr\'eter
raisonablement comme une structure de $F$-isocristal sur 
le champ $(X,x)^{isoc}$. 

En contre partie, le champ
$(X,x)^{F-isoc}$ est lui muni d'une structure
de $F$-isocristal, et le morphisme
$(X,x)^{isoc} \longrightarrow (X,x)^{F-isoc}$ est l'analogue
du morphisme 
$(X\otimes \mathbb{C})^{sch} \longrightarrow \overline{X}^{\mathbb{H}}$.
On peut d\'ecrire la structure de $F$-isocristal 
sur $(X,x)^{F-isoc}$ \`a l'aide de la notion
de type d'homotopie sch\'ematique non-neutre\footnote{Un type d'homotopie sch\'ematique
(non-neutre) sur un corps $k$ est un pr\'efaisceau simplicial sur $Aff/k$, tel qu'il existe une extension de
corps $k \subset K$, et un point $* \in F(K)$ faisant de $F$ un type d'homotopie
sch\'ematique point\'e sur $K$. Ils sont consid\'er\'es comme
des champs sur $k$ (i.e. des objets de $Ho(k)$).} de la fa\c{c}on suivante.

Notons $\mathbb{H}_{F}$ la gerbe affine (sur $\mathbb{Q}_{p}$)
dual de la cat\'egorie
Tannakienne des $F$-isocristaux sur $Spec\, k$ (voir par exemple
\cite[VI \S 4.3]{sa}). Comme
$k$ est alg\'ebriquement clos le lien de $\mathbb{H}_{F}$ est un sch\'ema en groupes diagonalisable
dont le groupe des caract\`eres est $\mathbb{Q}$ (ceci correspond \`a la d\'ecomposition
des $F$-isocristaux \`a l'aide des pentes).
La structure de $F$-isocristal sur $(X,x)^{F-isoc}$ peut s'interpr\'eter, \`a travers
la correspondence entre $F$-isocristal sur $Spec\, k$ et repr\'esentations
lin\'eaires de $\mathbb{H}_{F}$, comme un type d'homotopie sch\'ematique (non-neutre car la gerbe
$\mathbb{H}_{F}$ elle m\^eme est non-neutre)
$F^{0}$ sur $\mathbb{Q}_{p}$
muni d'une action de la gerbe $\mathbb{H}_{F}$. Ce champ $\mathbb{H}_{F}$-\'equivariant
s'ins\'ere dans une suite de fibrations de type d'homotopie sch\'ematiques (non point\'es)
$$F^{0} \longrightarrow [F^{0}/\mathbb{H}_{F}] \longrightarrow \mathbb{H}_{F}.$$
Le champ du milieu est reli\'e de tr\`es pr\`es au
\textit{type d'homotopie \'etale $p$-adique de $X$}. En effet, les syst\`emes locaux
de $\mathbb{Q}_{p}$-espaces vectoriels de dimension fini sur
$[F^{0}/\mathbb{H}_{F}]$ correspondent aux
$F$-isocrystaux sur $X$, et incluent donc les repr\'esentations continues
$p$-adiques de $\pi_{1}^{et}(X,x)$ de dimension finies. De plus, pour un syst\`eme
local $L$ sur $[F^{0}/\mathbb{H}_{F}]$ associ\'e \`a une
repr\'esentation $p$-adique $V$, la cohomologie de $[F^{0}/\mathbb{H}_{F}]$ \`a valeurs
dans $L$ s'identifie avec la cohomologie \'etale $p$-adique 
de $X$ \`a valeurs dans $L$.
La relation pr\'ecise entre $[F^{0}/\mathbb{H}_{F}]$ et le type d'homotopie
\'etale $p$-adique de $X$ peut aussi s'expliciter et l'on peut montrer que
$F^{0}$ muni de son action de $\mathbb{H}_{F}$ d\'etermine celui-ci a \'equivalence pr\`es
(voir \cite{ol} pour plus de d\'etails. D'une certaine fa\c{c}on, le type d'homotopie
\'etale $p$-adique de $X$ est la \textit{partie de pente $0$ de $F^{0}$}).

\subsubsection{Type d'homotopie $l$-adique}

Pour terminer je signale l'existence d'un type d'homotopie $l$-adique
d'un sch\'ema $X$ sur lequel $l$ est inversible. Il s'agit d'une construction pour
laquelle il n'existe pas de r\'ef\'erences d\'etaill\'ees (on trouvera
tout de m\^eme quelques lignes dans \cite{to2,to3}), et je me contenterais donc
de donner un bref apper\c{c}u des principales propri\'et\'es. S'il le d\'esire, le lecteur pourra
consid\'erer que les lignes qui suivent d\'ecrivent une situation conjecturale. \\

On fixe un sch\'ema $X$ et $x \in X(\overline{k})$ un point
g\'eom\'etrique. On peut alors construire un type d'homotopie
sch\'ematique point\'e sur $\mathbb{Q}_{l}$, $(X^{et}\otimes
\mathbb{Q}_{l})^{sch}$, appel\'e le \textit{type d'homotopie
$l$-adique de $X$} et jouissant des propri\'et\'es suivantes. Les
syst\`emes locaux de $\mathbb{Q}_{l}$-espaces vectoriels de
dimension finie sur $(X^{et}\otimes \mathbb{Q}_{l})^{sch}$ sont en
correspondance avec les repr\'esentations $l$-adiques continues de
$\pi_{1}^{et}(X,x)$. De plus, pour un tel syst\`eme local $L$,
la cohomologie de $(X^{et}\otimes \mathbb{Q}_{l})^{sch}$ s'identifie naturellement
avec la cohomologie de $X$ \`a valeurs dans le syst\`eme local \'etale
$l$-adique correspondant. Le champ
$(X^{et}\otimes \mathbb{Q}_{l})^{sch}$ est une version $l$-adique
de la sch\'ematisation de l'espace topologique sous-jacent \`a une vari\'et\'e
alg\'ebrique complexe.

Consid\`erons maintenant le cas o\`u $X$ est une vari\'et\'e sur
un corps $k$, et soit $\overline{X}=X\times_{Spec\, k}Spec\,
\overline{k}$ son extension \`a la clot\^ure s\'eparable de $k$. Le groupe
de Galois $H=Gal(\overline{k}/k)$ op\`ere par fonctorialit\'e
sur le champ \textit{non-point\'e}
$(\overline{X}^{et}\otimes \mathbb{Q}_{l})^{sch}$. Cette action
n'est pas continue, de la m\^eme fa\c{c}on que l'action
de $\mathbb{C}^{\times}$ donnant la d\'ecomposition de Hodge
d\'ecrite pr\'ec\'edemment n'est pas 
une action de $\mathbb{G}_{m}$.
On peut cependant consid\'erer un certain quotient
$(\overline{X}^{et}\otimes \mathbb{Q}_{l})_{ar}^{sch}$ de 
$(\overline{X}^{et}\otimes \mathbb{Q}_{l})^{sch}$, 
sur lequel $H$ va op\`erer de fa\c{c}on continue. Par d\'efinition, 
la cat\'egorie des syst\`emes locaux 
sur $(\overline{X}^{et}\otimes \mathbb{Q}_{l})_{ar}^{sch}$
est la sous-cat\'egorie Tannakienne des syst\`emes locaux 
$l$-adiques sur $\overline{X}$ engendr\'e par ceux qui proviennent
de syst\`emes locaux $l$-adiques sur $X$. De plus, pour $L$ un tel syst\`eme
local la cohomologie de $(\overline{X}^{et}\otimes \mathbb{Q}_{l})_{ar}^{sch}$
\`a coefficients dans $L$ s'identifie naturellement
avec la cohomologie de $\overline{X}$ \`a valeurs dans le syst\`eme local
$l$-adique correspondant.

Il existe une action naturelle de $H=Gal(\overline{k}/k)$
sur le type d'homotopie sch\'ematique non-point\'e 
$(\overline{X}^{et}\otimes \mathbb{Q}_{l})_{ar}^{sch}$, et cette
action
est \textit{continue} en un sens ad\'equat (e.g. on doit avoir que l'action de $H$
induite sur les $\mathbb{Q}_{l}$-espaces vectoriels (lin\'eairement compacts)
$\pi_{i}((\overline{X}^{et}\otimes \mathbb{Q}_{l})_{ar}^{sch})$ est continue
pour la topologie $l$-adique). Lorsque $k$ est par exemple un corps de nombre,
le champ $H$-\'equivariant
$(\overline{X}^{et}\otimes \mathbb{Q}_{l})_{ar}^{sch}$ est un analogue arithm\'etique du champ $\overline{X}^{\mathbb{H}}$ muni de sa d\'ecomposition
de Hodge, ou encore de $(X,x)^{F-Isoc}$ muni de sa structure
de $F$-isocristal (voir $\S 3.5.1$ et $\S 3.5.2$).

Par fonctorialit\'e, un point
rationnel $x \in X(k)$ induit un point fixe homotopique
de $H$ sur $(\overline{X}^{et}\otimes \mathbb{Q}_{l})_{ar}^{sch}$, 
ou en d'autres termes
un morphismes de champs $H$-\'equivariants
$* \longrightarrow (\overline{X}^{et}\otimes \mathbb{Q}_{l})_{ar}^{sch}$ 
(on rappelle que
l'ensemble simplicial des points fixes homotopiques d'une action de $H$ sur
un champ $F$ peut-\^etre d\'efini comme $\mathbb{R}\underline{Hom}_{H}(*,F)$,
o\`u $\mathbb{R}\underline{Hom}_{H}$ fait r\'ef\'erence \`a la cat\'egorie
de mod\`eles des pr\'efaisceaux simplicial munis d'une action de $H$).
On obtient donc une application
$$NAJ : X(k) \longrightarrow \pi_{0}(((\overline{X}^{et}\otimes \mathbb{Q}_{l})_{ar}^{sch})^{hH}),$$
de l'ensemble des points $k$-rationnels de $X$ vers l'ensemble des composantes
connexes des points fixes homotopiques de $H$. Cette application
s'appelle \textit{l'application d'Abel-Jacobi $l$-adique non-ab\'elienne}.
On justifie cette terminologie en remarquant qu'il existe un diagramme commutatif
$$\xymatrix{
X(k) \ar[r]^-{NAJ} \ar[d] & \pi_{0}(((\overline{X}^{et}\otimes \mathbb{Q}_{l})^{sch})_{ar}^{hH}) \ar[d] \\
CH_{0}(X) \ar[r]_-{\gamma} & H_{0}(X_{et},\mathbb{Q}_{l}),}$$
o\`u $\gamma$ est le morphisme \emph{classe de cycle}, et o\`u le morphisme vertical
de droite est induit par le morphisme d'ab\'elianisation
de $(\overline{X}^{et}\otimes \mathbb{Q}_{l})_{ar}^{sch}$ vers son type
d'homologie. Notons que le morphisme classe de cycles
permet de construire les applications d'Abel-Jacobi $l$-adiques
sup\'erieures (voir par exemple \cite{ja}), \`a l'aide de la suite spectrale de Leray-Serre
(qui d\'eg\'enerre en $E^{2}$ lorsque $X$ est propre et lisse)
$$E_{p,q}^{2}=H^{p}_{cont}(H,H_{q}(\overline{X}_{et},\mathbb{Q}_{l}) \Rightarrow
H_{q-p}(X_{et},\mathbb{Q}_{l}).$$

L'existence de cette application $NAJ$ semble montrer que le champ
$(\overline{X}^{et}\otimes \mathbb{Q}_{l})_{ar}^{sch}$, muni de son action
du groupe $H$, d\'etient des informations int\'eressantes sur
les propri\'et\'es arithm\'etiques de $X$. Il serait par exemple int\'eressant
de trouver des conditions sur la vari\'et\'e $X$ pour que l'application
$NAJ$ soit injective. On aimerait aussi trouver des exemples
de vari\'et\'es $X$ avec deux points $k$-rationnels $x$ et $y$, 
tel que $NAJ(x)\neq NAJ(y)$ mais $[x]=[y]$ dans le groupe de Chow 
des $0$-cycles $CH_{0}(X)$. Ceci impliquerait en particulier que 
l'application $NAJ$ d\'etecte des informations strictement plus fines que
les applications d'Abel-Jacobi sup\'erieures. 

Remarquons enfin que lorsque $X$ est un courbe lisse projective de genre au moins $2$ sur une corps
de nombre $k$, l'application $NAJ$ est une version
$l$-adique de l'application d\'ej\`a consid\'er\'ee
par A. Grothendieck dans son programme de \textit{g\'eom\'etrie anab\'elienne}, et 
qui \`a un point $k$-rationnel de $X$ associe un scindage
de la suite exacte fondamentale
$$\xymatrix{1 \ar[r] & \pi_{1}^{et}(\overline{X}) \ar[r] & \pi_{1}^{et}(X) \ar[r] & H=Gal(\overline{k}/k) \ar[r] &
1.}$$
A. Grothendieck conjecture en particulier que
cette derni\`ere application est bijective. Il parrait donc
raisonable d'essayer d'utiliser l'application $NAJ$ afin
de detecter les points rationnels de vari\'et\'es arithm\'etiques
qui ne sont plus des $K(\pi,1)$
(comme par exemple des sections hyperplanes
de $K(\pi,1)$). J'esp\`ere revenir sur ce sujet 
dans un travail ult\'erieur. 

\section{Champs $\infty$-g\'eom\'etriques}

Dans cette section on supposera que $k$ est un corps. Cette hypoth\`ese
n'est nullement n\'ecessaire mais simplifie quelque peu les
d\'efinitions.

Rappelons que pour un espace vectoriel $V$ de dimension finie, on peut construire
un champ alg\'ebrique (au sens d'Artin) $Alg_{V}$ classifiant les
alg\`ebres associatives et unitaires dont l'espace vectoriel sous-jacent est
isomorphe \`a $V$. Rappelons aussi qu'\'etant donn\'e une $k$-alg\`ebre
commutative $A$ de dimension fini on peut construire une vari\'et\'e
projective $Grass_{A}$ dont les points classifient les
chaines d'id\'eaux
$$I_{n} \hookrightarrow I_{n-1} \dots \hookrightarrow I_{0}=A.$$
Dans ce paragraphe nous allons construire des \textit{champs de modules}
analogues \`a $Alg_{V}$ et $Grass_{A}$ dans le cadre
de l'alg\`ebre homotopique (i.e. lorsque $V$ est un complexe
de $k$-vectoriel et $A$ une dga commutative, et que l'on ne s'int\'eresse
aux structures qu'\`a quasi-isomorphismes pr\`es). Nous montrerons en particulier
comment ces champs de modules sont recouverts par des champs affines
pour en faire des champs $\infty$-g\'eom\'etriques.

\subsection{D\'efinition}

On peut d\'efinir un champ alg\'ebrique (au sens d'Artin, voir par exemple \cite{lamo})
comme le champ quotient associ\'e \`a un groupoide
affine et lisse (pr\'ecisemment ceci ne donne que les champs
alg\'ebriques quasi-compact et \`a diagonale affine, mais nous nous suffirons de ceci).
Dans ce paragraphe nous allons g\'en\'eraliser la notion de champ alg\'ebrique
en consid\'erant des quotients associ\'es \`a des objets en groupoides
lisses (en un sens \`a d\'efinir) dans la cat\'egorie des champs affines. \\

Pour commencer nous aurons besoin de certaines notions
de finitude et de lissit\'e pour des morphismes
entre champs affines. \\

Soit $A \in k-Alg^{\Delta}$ une $k$-alg\`ebre co-simpliciale
(appartenant \`a $\mathbb{V}$ d'apr\`es nos conventions. On dispose
de la cat\'egorie $A/k-Alg^{\Delta}$, des objets sous $A$, munie
de sa structure de cat\'egorie de mod\`eles induite
(fibrations, \'equivalences et cofibrations test\'ees dans
$k-Alg^{\Delta}$), et dont les objets seront appel\'es
des \textit{$A$-alg\`ebres}. La cat\'egorie
$A/k-Alg^{\Delta}$ poss\`ede une structure simpliciale
(pour laquelle le foncteur d'oubli $A/k-Alg^{\Delta}
\longrightarrow k-Alg^{\Delta}$ pr\'eserve les exponentiations
par des ensembles simpliciaux) qui en fait une cat\'egorie
de mod\`eles simpliciale. On peut donc d\'efinir des
$Hom$'s simpliciaux d\'eriv\'es qui seront not\'es
$\mathbb{R}\underline{Hom}_{A-Alg}(-,-)$.

\begin{df}\label{d18}
\begin{itemize}
\item
\emph{Soit $A \in k-Alg^{\Delta}$ une $k$-alg\`ebre co-simpliciale
et $B\in A/k-Alg^{\Delta}$ une $A$-alg\`ebre. Nous dirons que
$B$ est \emph{de pr\'esentation finie sur $A$} (ou encore
que le morphism $A \rightarrow B$ est de pr\'esentation finie)
si pour tout syst\`eme inductif filtrant de $A$-alg\`ebres
$\{C\}_{i\in I}$ (o\`u $I\in \mathbb{V}$) le morphism naturel
$$Colim_{i\in I} \mathbb{R}\underline{Hom}_{A-Alg}(B,C_{i}) \longrightarrow
\mathbb{R}\underline{Hom}_{A-Alg}(B,Colim_{i\in I}C_{i})$$
est une \'equivalence.}

\item \emph{Un morphisme dans $Ho(k-Alg^{\Delta})$ est \emph{de pr\'esentation finie}
si l'on peut le repr\'esenter dans $k-Alg^{\Delta}$ (\`a isomorphisme pr\'es)
par un morphisme de pr\'esentation finie au sens ci-dessus.}

\item \emph{Un morphisme de champs affines est \emph{de pr\'esentation finie} s'il correspond
par Cor. \ref{c1} \`a un morphisme de pr\'esentation finie dans
$Ho(k-Alg^{\Delta})$.}
\end{itemize}
\end{df}

\textit{Remarques:}
\begin{itemize}
\item On remarque ais\'emment que si un morphisme
de $Ho(k-Alg^{\Delta})$ poss\`ede un repr\'esentant dans
$k-Alg^{\Delta}$ qui soit de pr\'esentation finie alors il en est de m\^eme
de tous ses repr\'esentants. La notion de morphismes de pr\'esentation
finie dans $Ho(k-Alg^{\Delta})$ est donc raisonable.
\item
Par adjonction on voit imm\'ediatement qu'un morphisme
de $k$-alg\`ebres $A \longrightarrow B$ (non co-simpliciales) est de pr\'esentation
finie comme morphisme de $k-Alg^{\Delta}$ si et seulement si $B$ est une
$A$-alg\`ebre de pr\'esentation finie au sens usuel. La d\'efinition
\ref{d18} est donc une g\'en\'eralisation de la notion usuelle.
\end{itemize}

Nous passons maintenant \`a la notion de morphismes \'etales et
de morphismes lisses.

\begin{df}\label{d21}
\emph{Soit $f : F \longrightarrow G$ un morphisme de champs affines.}
\begin{itemize}
\item \emph{Nous dirons que $f$ est \emph{formellement \'etale} si pour toute
$k$-alg\`ebre $A\in k-Alg$, et $I\subset A$ un id\'eal de carr\'e nul,
le morphisme induit
$$F(A) \longrightarrow F(A/I)\times^{h}_{G(A/I)}G(A)$$
est une \'equivalence.}

\item \emph{Nous dirons que $f$ est \emph{\'etale} s'il est formellement
\'etale et de pr\'esentation finie}
\item \emph{Nous dirons que $f$ est un \emph{recouvrement \'etale} s'il
est \'etale et si de plus le morphisme de faisceaux
$\pi_{0}(F) \longrightarrow \pi_{0}(G)$ est un \'epimorphisme.}
\end{itemize}
\end{df}

\textit{Remarque:}
La notion de morphisme \'etale de champs affines g\'en\'eralise
la notion usuelle pour les sch\'emas affines. Cependant il est important de noter
qu'il existe des recouvrements \'etales de champs affines $F \longrightarrow Spec\, k$
o\`u $F$ n'est pas repr\'esentable un sch\'ema affine. En effet, on v\'erifie
imm\'ediatement que lorsque $k=\mathbb{F}_{p}$ le champ
$K(\mathbb{Z}/p,1)$ est un tel exemple. \\

Soit $E$ un $k$-espace vectoriel co-simplicial. Nous dirons que $E$ est
\emph{parfait} si $H^{*}(E)$ est de dimension finie sur $k$.

Rappelons que
pour $E$ un $k$-espace vectoriel co-simplicial on dispose de la $k$-alg\`ebre
co-simpliciale libre sur $E$ not\'ee $S^{*}(E)$. Un morphisme
dans $Ho(k-Alg^{\Delta})$ sera dit \emph{parfait} s'il est isomorphe \`a un morphisme
de la forme
$$A \longrightarrow A\otimes_{k}S^{*}(E),$$
o\`u $E$ parfait.
De part la propri\'et\'e universelle
satisfaite par $S^{*}(E)$ on voit qu'un morphisme
parfait est toujours de pr\'esentation fini au sens
de Def. \ref{d18}.
Par le plongement $\mathbb{R}Spec$ de Cor. \ref{c1} on \'etent la notion
de morphismes parfaits au cas des morphismes de champs affines.

\begin{df}\label{d22}
\emph{Un morphisme de champs affines $f : F \longrightarrow G$  est lisse s'il est de pr\'esentation fini
et s'il existe des champs affines $F'$ et $G'$ et un diagramme commutatif de champs
$$\xymatrix{
F' \ar[d]_-{f'} \ar[r]^-{u} & F \ar[d]^-{f} & \\
G' \ar[r]_-{v} & G,}$$
o\`u $f'$ est un morphisme \'etale,
$u$ est un recouvrement \'etale, et $v$
est un morphisme parfait.}
\end{df}

On v\'erifie que les morphismes lisses entre champs affines sont stables par composition
et changements de bases (i.e. produits fibr\'es homotopiques). On v\'erifie aussi
qu'un morphisme entre sch\'emas affines est lisse au sens de la d\'efinition
\ref{d22} si et seulement s'il est lisse au sens usuel. \\

Soit $X_{\bullet} : \Delta^{op} \longrightarrow SPr(k)$ un objet simplicial
dans la cat\'egorie des pr\'e\-fai\-sc\-eaux simpliciaux. Pour tout
entier $n>1$ on dispose du $n$-\'eme morphisme de Segal (bien d\'efini dans $Ho(SPr(k))$)
$$\sigma_{n} : X_{n} \longrightarrow \underbrace{X_{1}\times^{h}_{X_{0}}X_{1}
\dots \times^{h}_{X_{0}}X_{1}}_{n \; fois},$$
induit par les $n$ morphismes $[1] \longrightarrow [n]$ dans $\Delta$ qui envoient
$0$ et $1$ sur $i$ et $i+1$ (pour $0\leq i<n$). De m\^eme, on dispose
du morphisme
$$i : X_{2} \longrightarrow X_{1}\times^{h}_{X_{0}}X_{1},$$
induit par les morphismes $d_{1},d_{2} : [1] \longrightarrow [2]$, bien d\'efini
dans $Ho(SPr(k))$. Lorsque tous les morphismes $\sigma_{n}$ et $i$ sont des \'equivalences
nous dirons que $X_{\bullet}$ est un \emph{groupoide de Segal dans $SPr(k)$}. \\

Nous arrivons enfin \`a la d\'efinition des champs $\infty$-g\'eom\'etriques. La d\'efinition
que nous donnons ici n'est pas la plus g\'en\'erale car elle n'englobe que le
cas des quotients de champs affines par des champs en groupoides affines
et lisses. Quoiqu'il en soit cette d\'efinition sera suffisante pour
traiter les exemples qui nous int\'eressent dans cet article.
Rappelons que pour $X_{\bullet}$ un objet simplicial
de $SPr(k)$, en d'autres termes un pr\'efaisceau bi-simplicial,
nous notons $|X_{\bullet}|$ son pr\'efaisceau simplicial
diagonal. L'objet $|X_{\bullet}|$ est naturellement isomorphe
dans $Ho(SPr(k))$ \`a la colimite homotopique
du diagramme simplicial $n \mapsto X_{n}$.

\begin{df}\label{d23}
\emph{Un champ $F$ est $\infty$-g\'eom\'etrique s'il
est isomorphe \`a un champ de la forme $|X_{\bullet}|$, o\`u $X_{\bullet}$ est
un groupoide de Segal dans $SPr(k)$ qui satisfait aux deux conditions suivantes.}
\begin{itemize}
\item \emph{Les champs $X_{0}$ et $X_{1}$ sont des champs affines.}

\item \emph{Le morphisme $d_{0} : X_{1} \longrightarrow X_{0}$
est un morphisme lisse de champs affines.}
\end{itemize}
\end{df}

Il est imm\'ediat de v\'erifier les propri\'et\'es suivantes.

\begin{itemize}
\item Tout champ affine est un champ $\infty$-g\'eom\'etrique.

\item Un produits fibr\'es homotopique de champs $\infty$-g\'eom\'etriques
est un champ $\infty$-g\'eom\'etrique. En particulier les champs $\infty$-g\'eom\'etriques
sont stables par limites homotopiques finies.

\item Un champ alg\'ebrique au sens d'Artin (voir \cite{lamo}) qui est quasi-compact et dont
la diagonale est un morphisme repr\'esentable et affine
est un champ $\infty$-g\'eom\'etrique.

\item Tout champ tr\'es pr\'esentable au sens de Def. \ref{d11} est un champ $\infty$-g\'eom\'e\-tri\-que.

\end{itemize}

La th\'eorie des champs $\infty$-g\'eom\'etriques peut \'evidemment se poursuivre dans les m\^emes
directions que celles suivit par exemple dans \cite{lamo,s5}. Nous ne le ferons pas ici et nous
renvoyons \`a \cite{hagdag,tv2} pour plus de d\'etails sur une th\'eorie plus g\'en\'erale.

\subsection{Application}

Dans ce paragraphe on pr\'esente un exemple d'application
de la notion de champs $\infty$-g\'eom\'etrique du paragraphe
pr\'ec\'edent. \\

\subsubsection{Le champ des structures multiplicatives}

On fixe un complexe de $k$-espaces vectoriels $E$, \`a cohomologie born\'ee
et de dimension finie. On cherche \`a d\'efinir
le champ classifiant les structures multiplicatives sur $E$ \`a quasi-isomorphisme pr\`es, ou en d'autres
termes les $k$-alg\`ebres diff\'erentielles gradu\'ees (associatives et unitaires par exemple)
dont le complex sous-jacent est quasi-isomorphe \`a $E$.

Pour cela, on rappelle que pout toute anneau $A$ (associatif, commutatif, unitaire),
il existe une structure de cat\'egorie de mod\`eles
sur la cat\'egorie $A-DGA$, des $A$-alg\`ebres diff\'erentielles gradu\'ees
(non born\'ees). Les \'equivalences (resp. les fibrations) pour cette structure
sont les quasi-isomorphismes (resp. les \'epimorphismes). L'existence
de cette structure de mod\`eles est d\'emontr\'ee par exemple dans
\cite{ss}.

Pour toute $k$-alg\`ebre $A \in k-Alg$, on d\'efinit une
cat\'egorie $\underline{Alg}_{E}(A)$ de la fa\c{c}on suivante. Ses
objets sont les $A$-alg\`ebres diff\'erentielles gradu\'ees $B\in A-DGA$,
cofibrantes pour la structures de mod\`eles
d\'ecrite dans \cite{ss}, et telles qu'il existe un morphisme
fid\`element plat $A \longrightarrow A'$ tel que les complexe de $A'$-modules
$B\otimes_{A}^{\mathbb{L}}A'$ et
$E\otimes_{k}^{\mathbb{L}}A'$ soient quasi-isomorphes. Les morphismes dans
$\underline{Alg}_{E}(A)$ sont les \'equivalences (i.e. les quasi-isomorphismes) de
$A$-alg\`ebres diff\'erentielles gradu\'ees.

Pour un morphisme de $k$-alg\`ebres $A \rightarrow A'$, on dispose
d'un foncteur de changement de base
$$\begin{array}{ccc}
\underline{Alg}_{E}(A) & \longrightarrow & \underline{Alg}_{E}(A') \\
 B & \mapsto & B\otimes_{A}A',
\end{array}$$
qui est bien d\'efini car l'on s'est restreint aux $A$-alg\`ebres diff\'erentielles gradu\'ees
cofibrantes, et donc plates sur $A$ (ainsi le foncteur ci-dessus pr\'eserve bien
les \'equivalences). Ceci d\'efinit un pseudo-foncteur
$$\begin{array}{cccc}
\underline{Alg}_{E} : & k-Alg & \longrightarrow & Cat \\
 & A & \mapsto & \underline{Alg}_{E}(A),
\end{array}$$
que l'on s'empresse de strictifier en un vrai foncteur
par le proc\'ed\'e standard (voir par exemple \cite[Thm. 3.4]{ma}).

Finalement, on compose avec le foncteur associant \`a toute cat\'egorie
son nerf pour obtenir un pr\'efaisceau simplicial
$$\begin{array}{cccc}
B\underline{Alg}_{E} : & k-Alg & \longrightarrow & SSet \\
 & A & \mapsto & B\underline{Alg}_{E}(A),
\end{array}$$
et donc un objet de $SPr(k)$. Il se trouve que
$B\underline{Alg}_{E}$ ainsi d\'efini est un champ au sens de Def.
\ref{d1}. La preuve de ce fait, qui semble cependant un r\'esultat
folklorique, est longue et technique (elle suit la d\'emarche
adopt\'ee dans \cite{s1} pour d\'emontrer que le pr\'echamp des
complexes est un champ) et nous ne la donnerons pas. Le th\'eor\`eme suivant
donne un premier exemple non-trivial de champ
$\infty$-g\'eom\'etrique. Nous nous contenterons d'esquisser sa preuve qui
sera reprise en plus grande g\'en\'eralit\'e dans \cite{tv2}.

\begin{thm}\label{t20}
Le champ $B\underline{Alg}_{E}$ est un champ $\infty$-g\'eom\'etrique.
\end{thm}

\textit{Esquisse de preuve:} Nous allons
montrer que $B\underline{Alg}_{E}$ est en r\'ealit\'e le champ quotient
d'un champ affine par un champ en groupes affine et lisse.

On commence par d\'efinir
un champ $B\underline{Mod}_{E}$ de la fa\c{c}on suivante.
Pour toute $k$-alg\`ebre $A$ on d\'efinit une cat\'egorie
$\underline{Mod}_{E}(A)$, dont les objets sont les complexes de $A$-modules
plats
qui sont localement, pour la topologie plate sur $A$, quasi-isomorphes \`a 
$E\otimes_{k}A$,
et les morphismes sont les quasi-isomorphismes de complexes de $A$-modules.
Lorsque
$A$ varie dans $k-Alg$, ceci d\'efinit un pseudo-foncteur
$A \mapsto \underline{Mod}_{E}(A)$ de $k-Alg$
vers $Cat$ que l'on strictifie et \`a qui l'on applique le foncteur
classifiant $B : Cat \longrightarrow SSet$. Ceci nous donne
le champ cherch\'e $B\underline{Mod}_{E}$.

Bien entendu, il existe un morphisme de champs $B\underline{Alg}_{E} \longrightarrow
B\underline{Mod}_{E}$, qui oubli la structure d'alg\`ebre diff\'erentielle
gradu\'ee. De plus, par \cite[2.3]{dk} on voit que le champ
$B\underline{Mod}_{E}$ est isomorphe \`a $K(H,1)$, o\`u $H$ est le
pr\'efaisceau en monoides simpliciaux d\'efini par
$$\begin{array}{cccc}
H : & k-Alg & \longrightarrow & SMon \\
 & A & \mapsto & Aut_{C(A)}(E\otimes_{k}A),
\end{array}$$
o\`u $Aut_{C(A)}(E\otimes_{k}A)$ d\'esigne le sous-ensemble simplicial
du l'ensemble simplicial des morphismes (''mapping spaces'')
$Map_{C(A)}(E\otimes_{k}A,E\otimes_{k}A)$ de la cat\'egorie de mod\`eles
des complexes (non born\'es) de $A$-modules (voir par exemple \cite{hin,ho}
pour la description de la structure de mod\`eles $C(A)$ et la d\'efinition
des mapping spaces) form\'e des quasi-isomorphismes.

\begin{lem}\label{l30}
Le champ $H$ est affine et lisse.
\end{lem}

\textit{Preuve:} On commence par consid\`erer le champ des endomorphismes de $E$
$$\begin{array}{cccc}
\mathbb{R}\underline{End}(E) : & k-Alg & \longrightarrow & SSet \\
 & A & \mapsto & Map_{C(A)}(E\otimes_{k}A,E\otimes_{k}A).
\end{array}$$
Il est facile de voir que ce champ est isomorphe \`a
$\mathbb{R}Spec\, S^{*}(M)$, o\`u $M$ est le $k$-module co-simplicial
associ\'e au complexe tronqu\'e $\tau_{\geq 0}E\otimes_{k}E^{*}$ par la correspondence
de Dold-Puppe. Ainsi, le champ $\mathbb{R}\underline{End}(E)$ est affine.

On regarde maintenant le sous-pr\'efaisceau en monoides
simpliciaux $H$ de $\mathbb{R}\underline{End}(E)$ form\'e
des quasi-isomorphismes. 
Il existe un diagramme homotopiquement cart\'esien de champs
$$\xymatrix{
H \ar[d] \ar[r] & \pi_{0}(H)  \ar[d] \\
\mathbb{R}\underline{End}(E) \ar[r] & \pi_{0}(\mathbb{R}\underline{End}(E)).}$$
De plus, le faisceau en monoides $\pi_{0}(\mathbb{R}\underline{End}(E))$
est isomorphe au produit $\prod_{i}\mathcal{M}_{n_{i}}$, o\'u
$\mathcal{M}_{n_{i}}$ est le sch\'ema affine des matrices car\'ees de rangs
$n_{i}$, et $n_{i}=dim\, H^{i}(E)$. De m\'eme,
$\pi_{0}(H)$ est le sous-faisceau des \'el\'ements inversibles dans
$\pi_{0}(\mathbb{R}\underline{End}(E))$, ou en d'autres termes
est isomorphe au produit $\prod_{i}Gl_{n_{i}}$. Le diagrame homotopiquement
cart\'esien pr\'ec\'edent montre donc que $H$ s'\'ecrit comme un produit
fibr\'e homotopique de champs affines et de pr\'esentation finie (sur 
$Spec\, k$) et donc
est lui-m\^eme un champ affine et de pr\'esentation finie. \\

Il nous reste a voir que $H \longrightarrow Spec\, k$ est lisse au sens de la d\'efinition
Def. \ref{d22}.
Or, nous avons vu que le morphisme $\mathbb{R}\underline{End}(E) \longrightarrow Spec\, k$
est un morphisme parfait et donc lisse. Enfin, il est imm\'ediat de v\'erifier
que le morphisme naturel $H=\mathbb{R}\underline{Aut}(E)\longrightarrow \mathbb{R}\underline{End}(E)$
est \'etale au sens de la d\'efinition Def. \ref{d21}, et donc lisse.
Ceci montre que $H$ est un champ affine lisse.
\hfill $\Box$ \\

On vient donc de construire un morphisme de champs
$$B\underline{Alg}_{E} \longrightarrow K(H,1),$$
o\`u $H$ est un pr\'efaisceau en monoides simpliciaux
dont le champ sous-jacents est affine et lisse.
De plus, par d\'efinition le faisceau en monoides
$\pi_{0}(H)$ est un faisceaux en groupes (i.e.
$H$ est un $H_{\infty}$-champs au sens de Def.
\ref{d5}. On sait alors (voir par exemple \cite[Prop. 1.5]{to1}) que
$H$ est \'equivalent, en tant que pr\'efaisceau en monoides simpliciaux,  \`a
un pr\'efaisceau en groupes simpliciaux.
On peut alors appliquer les techniques
de champs \'equivariants, d\'evelopp\'ees dans \cite{kt}, qui nous apprennent
que le champ $B\underline{Alg}_{E}$ est le champ
quotient $[X/H]$, o\`u $X$ est la fibre homotopique
du morphisme $B\underline{Alg}_{E} \longrightarrow K(H,1)$.
Le champ quotient $[X/H]$ s'\'ecrit aussi de la forme
$|X_{*}|$, o\`u $X_{*}$ est le groupoide classifiant de l'action
de $H$ sur $X$, d\'efini par la formule
$X_{n}:=X\times H^{n-1}$. Ainsi, le lemme \ref{l30} montre
que pour d\'emontrer le th\'eor\`eme \ref{t20} il nous suffit de voir que $X$ est un champ affine.

\begin{lem}\label{l31}
Le champ $X$, fibre homotopique du morphisme
$$B\underline{Alg}_{E} \longrightarrow B\underline{Mod}_{E}$$
est un champ affine.
\end{lem}

\textit{Preuve:} Ce lemme est une cons\'equence
d'un des r\'esultats principal de \cite{re}.

On rappelle qu'il existe une structure de cat\'egorie de mod\`eles
sur la cat\'egorie des op\'erades (unitaires) dans la cat\'egorie
monoidale $C(k)$ des complexes $k$-espaces vectoriels (voir par exemple
\cite{hin}).
Pour toute $k$-alg\`ebre $A\in k-Alg$ on d\'efinit une
op\'erade $\underline{End}(E)_{A}$ dans $C(k)$ par la formule
$$\underline{End}(E)_{A}(n):=\underline{Hom}_{k}(E^{\otimes_{k} n},E\otimes_{k}A),$$
o\`u $\underline{Hom}_{k}$ d\'esigne les $Hom$ internes de la cat\'egorie
des complexes de $k$-espaces vectoriels. Pour un morphisme de $k$-alg\`ebres
$A \longrightarrow A'$ on dispose d'un morphisme d'op\'erades
$\underline{End}(E)_{A} \longrightarrow \underline{End}(E)_{A'}$, et ceci fait
de $A \mapsto \underline{End}(E)_{A}$ un pr\'efaisceau en op\'erades
unitaires sur la cat\'egorie $Aff/k$.

On d\'efinit alors un foncteur
$$\begin{array}{cccc}
\mathbb{R}\underline{Hom}(\mathcal{ASS},\underline{End}(E)) : & k-Alg & \longrightarrow  & SSet \\
 & A & \mapsto & Map_{Op(k)}(\mathcal{ASS},\underline{End}(E)_{A}),
\end{array}$$
o\`u $Map_{Op(k)}$ d\'esigne les mapping spaces de la cat\'egorie de mod\`eles
des op\'erades unitaires dans $C(k)$ (voir par exemple \cite[\S 5.2]{ho}),
et $\mathcal{ASS}$ est l'op\'erade unitaire finale (i.e. celle qui classifie
les alg\`ebres associatives et unitaires).
On obtient ainsi un objet $\mathbb{R}\underline{Hom}(\mathcal{ASS},\underline{End}(E)) \in Ho(SPr(k))$.
Le th\'eor\`eme \cite[Thm. 1.1.5]{re} (ou du moins sa version o\`u la cat\'egorie
de mod\`eles $\mathcal{M}_{R}$ est rempmlac\'ee par $C(k)$) nous dit que le champ $X$
est isomorphe au champ
$\mathbb{R}\underline{Hom}(\mathcal{ASS},\underline{End}(E))$. Il nous reste donc \`a montrer que
le champ $\mathbb{R}\underline{Hom}(\mathcal{ASS},\underline{End}(E))$ est affine.

Pour cela, on \'ecrit $\mathcal{ASS} \simeq Hocolim_{n\in \Delta^{op}}\mathcal{O}_{n}$, o\`u
chaque op\'erade $\mathcal{O}_{n}$ est une op\'erade libre. Ainsi, on a une \'equivalence naturelle
de champs
$$\mathbb{R}\underline{Hom}(\mathcal{ASS},\underline{End}(E)) \simeq
Holim_{n\in \Delta}\mathbb{R}\underline{Hom}(\mathcal{O}_{n},\underline{End}(E)),$$
o\`u $\mathbb{R}\underline{Hom}(\mathcal{O}_{n},\underline{End}(E))$ est d\'efini
comme $\mathbb{R}\underline{Hom}(\mathcal{ASS},\underline{End}(E))$ en rempla\c{c}ant
l'op\'erade $\mathcal{ASS}$ par $\mathcal{O}_{n}$. Par la proposition Prop. \ref{p2'} il nous
faut donc montrer que pour une op\'erade libre $\mathcal{O}$
le champ $\mathbb{R}\underline{Hom}(\mathcal{O},\underline{End}(E))$
est un champ affine. Mais dire que $\mathcal{O}$ est libre signifie qu'il existe
une famille de complexes $D_{n}$ pour $n\geq 1$ et des \'equivalences fonctorielles
en $A$
$$\mathbb{R}\underline{Hom}(\mathcal{O},\underline{End}(E))(A)\simeq
\prod_{n}Map_{C(k)}(E^{\otimes_{k}n}\otimes_{k}D_{n},E\otimes_{k}A)$$
$$\simeq \prod_{n}Map_{C(k)}(E^{\otimes_{k}n}\otimes_{k}E^{*}\otimes_{k}D_{n},A),$$
o\`u $Map_{C(k)}$ sont les mapping spaces de la cat\'egorie de mod\`eles
des complexes de $k$-espaces vectoriels.
Ceci montre qu'il suffit de v\'erifier que pour un complexe fix\'e $C$, le champ d\'efini par
$A \mapsto Map_{C(k)}(C,A)$ est un champ
affine. Mais il est clair que ce champ est isomorphe \`a $\mathbb{R}Spec\, S^{*}(M)$, o\`u
$M$ est le $k$-espace vectoriel co-simplicial obtenu en appliquant la correspondence
de Dold-Puppe au complexe tronqu\'e $\tau_{\geq 0}C$. \hfill
$\Box$ \\

En conclusion, $B\underline{Alg}_{E}$ s'\'ecrit
$[X/H]$, o\`u $X$ est un champ affine et
$H$ est un \emph{champ en groupes affine et lisse} qui op\`ere sur $X$. Ceci
montre que $B\underline{Alg}_{E}$ est un champ $\infty$-g\'eom\'etrique
au sens de la d\'efinition \ref{d23}. \hfill $\Box$ \\

\textit{Remarques:}
\begin{itemize}
\item On pourrait aussi d\'emontrer un th\'eor\`eme analogue
au th\'eor\`eme Thm. \ref{t20} pour le champs classifiants les sturctures
multiplicatives associatives, unitaires et \emph{commutatives}. Dans le cas
o\`u $k$ n'est pas de carat\'eristique nulle il faut cependant utiliser
la notion de $E_{\infty}$-alg\`ebres au lieu de celle d'alg\`ebres
diff\'erentielles gradu\'ees commutatives.
\item Des espaces de modules formels de structures mutiplicatives avaient d\'ej\`a
\'et\'e construits dans \cite{koso,hin2}, et notre champ
$B\underline{Alg}_{E}$ en est une contre-partie globale et
\emph{tronqu\'ee}. La version non-tronqu\'ee, ou encore
\emph{d\'eriv\'ee} du champ $B\underline{Alg}_{E}$ est
d\'efinie et \'etudi\'ee dans \cite{hagdag} (voir aussi \cite{tv2}).
\end{itemize}

\subsubsection{P\'eriodes non-ab\'eliennes}

On suppose maintenant que $k=\mathbb{C}$. Pour une vari\'et\'e lisse projective
complexe $X$, on dispose de sa $\mathbb{C}$-alg\`ebre diff\'erentielle gra\-du\'ee
commutative de cohomologie $C^{*}(X^{top},\mathbb{C})\simeq C^{*}(X,\Omega_{X}^{*})$. Elle est muni
d'une filtration de Hodge
$$F^{i}C^{*}(X,\Omega_{X}^{*}):=C^{*}(X,\Omega_{X}^{>i}),$$
qui est une filtration par des id\'eaux diff\'erentiels gradu\'es. Dans cette section on se propose
d'\'etudier la variation de la filtration de Hodge sur $C^{*}(X^{top},\mathbb{C})$ lorsque
l'on d\'eforme la vari\'et\'e $X$, et ce \`a l'aide d'une \textit{application
des p\'eriodes non-ab\'eliennes}. L'approche que nous proposons dans
ce paragraphe n'est qu'un premi\`ere approximation et ne tient pas compte
de toute la structure (en particulier la transversalit\'e de Griffiths est
n\'eglig\'ee), et demandera certainement \`a \^etre modifi\'ee
et complet\'ee dans un travail ult\'erieur. \\

On rappelle que pour toute $\mathbb{C}$-alg\`ebre $B \in \mathbb{C}-Alg$,
il existe une structure de cat\'egorie de mod\`eles
sur la cat\'egorie $B-CDGA_{+}$, des $\mathbb{C}$-alg\`ebres
diff\'erentielles gradu\'ees commutatives concentr\'ees en degr\'es positifs, o\`u les \'equivalences sont
les quasi-isomorphismes et les fibrations sont les surjections (voir par
exemple \cite{bg} ou encore la preuve de Thm. \ref{t3}). On fixe un objet $A \in \mathbb{C}-CDGA_{+}$ \`a cohomologie
born\'ee et de dimension finie. Nous allons d\'efinir un champ
$B\underline{Filt}^{n}_{A}$, classifiant les filtrations de longueur $n$ sur $A$.

Pour $B\in \mathbb{C}-Alg$, on d\'efinit une cat\'egorie
$\underline{Filt}^{n}_{A}(B)$ de la fa\c{c}on suivante. Les objets
de $\underline{Filt}^{n}_{A}(B)$ sont les diagrames dans $B-CDGA_{+}$
$$\xymatrix{A\otimes_{\mathbb{C}}B=A^{(n+1)} \ar[r] & A^{(n)} \ar[r] &
\dots A^{(1)} \ar[r] & A^{(0)}=0,}$$
v\'erifiant les deux conditions suivantes
\begin{itemize}
\item Pour tout entier $i\leq n$, l'objet $A^{(i)} \in B-CDGA_{+}$ est cofibrant.
\item
Pour tout entier $i$  le morphisme induit
$$H^{*}(A\otimes_{\mathbb{C}}B)\simeq H^{*}(A)\otimes_{\mathbb{C}}B \longrightarrow
H^{*}(A^{(i)})$$
est un \'epimorphisme scind\'e de $B$-modules (en particulier, on voit que tous les
$B$-modules $H^{*}(A^{(i)})$ sont projectifs et de type fini, et tous les morphismes
$H^{*}(A^{(i)}) \longrightarrow H^{*}(A^{(i-1)})$ sont des \'epimorphismes
scind\'es de $B$-modules).
\end{itemize}

Les morphismes dans $\underline{Filt}^{n}_{A}(B)$
sont les diagrames commutatifs dans $B-CDGA_{+}$
$$\xymatrix{ & A^{(n)} \ar[dd] \ar[r] & A^{(n-1)} \ar[r] \ar[dd] &
\dots A^{(1)} \ar[dd] \ar[rd] &  \\
A\otimes_{\mathbb{C}}B \ar[ru] \ar[rd] & & & & 0 \\
 & (A')^{(n)} \ar[r] & (A')^{(n-1)} \ar[r] &
\dots (A')^{(1)} \ar[ru] & }$$
o\`u les morphismes verticaux sont tous des quasi-isomorphismes.

Pour un morphisme de $\mathbb{C}$-alg\`ebres $B \longrightarrow B'$ on dispose
d'un foncteur de changement de base
$\underline{Filt}^{n}_{A}(B)  \longrightarrow   \underline{Filt}^{n}_{A}(B')$,
qui \`a un objet
$$\xymatrix{A\otimes_{\mathbb{C}}B=A^{(n+1)} \ar[r] & A^{(n)} \ar[r] &
\dots A^{(1)} \ar[r] & A^{(0)}=0,}$$ associe l'objet
$$\xymatrix{A\otimes_{\mathbb{C}}B'=A^{(n+1)}\otimes_{\mathbb{B}}B' \ar[r] & A^{(n)}\otimes_{\mathbb{B}}B' \ar[r] &
\dots A^{(1)}\otimes_{B}B' \ar[r] & 0.}$$
Ceci d\'efinit un pseudo-foncteur $B \mapsto \underline{Filt}^{n}_{A}(B)$, que l'on strictifie et au quel
on applique le foncteur classifiant pour obtenir un objet
$B\underline{Filt}^{n}_{A} \in SPr(\mathbb{C})$. \\

La proposition suivante est une g\'en\'eralisation de l'existence des
vari\'et\'es Grassmaniennes. Nous ne la d\'emontrerons pas.

\begin{prop}\label{p30}
Le pr\'efaisceau simplicial $B\underline{Filt}^{n}_{A}$ est un champ $\infty$-g\'eom\'etrique.
\end{prop}

Supposons maintenant que $B$ soit une $\mathbb{C}$-alg\`ebre locale Artinienne et
$p : \mathcal{X} \longrightarrow Spec\, B$ un morphisme projectif et lisse.
Sur le site Zariski de $\mathcal{X}$, on dispose d'un diagrame de faisceaux
de $B-CDGA_{+}$
$$\xymatrix{
\Omega^{*}_{\mathcal{X}/B} \ar[r] & \Omega_{\mathcal{X}/B}^{\leq n-1} \ar[r] &
\Omega_{\mathcal{X}/B}^{\leq n-2} & \dots \ar[r] & \Omega_{\mathcal{X}/B}^{\leq 0}=\mathcal{O}_{\mathcal{X}},}$$
o\`u $\Omega^{*}_{\mathcal{X}/B}$ est le complexe de de Rham (alg\'ebrique) relatif du morphismes $p$.
En utilisant une structure de mod\`eles ad\'equate sur la cat\'egorie des faisceaux
de $B-CDGA_{+}$ (voir par exemple \cite{kt}) on peut prendre l'image directe d\'eriv\'ee
de ce diagrame
$$\xymatrix{
\mathbb{R}p{*}\Omega^{*}_{\mathcal{X}/B} \ar[r] & \mathbb{R}p_{*}\Omega_{\mathcal{X}/B}^{\leq n-1} \ar[r] &
\mathbb{R}p_{*}\Omega_{\mathcal{X}/B}^{\leq n-2} & \dots \ar[r] &
\mathbb{R}p_{*}\mathcal{O}_{\mathcal{X}},}$$
qui est un diagrame dans $B-CDGA_{+}$.

En utilisant la connexion de Gauss-Manin on peut voir qu'il existe une \'equivalence
(non canonique)
dans $B-CDGA_{+}$
$$\mathbb{R}p{*}\Omega^{*}_{\mathcal{X}/B} \simeq
C^{*}(X^{top},\mathbb{C})\otimes_{\mathbb{C}}B,$$
o\`u $X:=\mathcal{X}\otimes_{B}B/m$ est la fibre sp\'eciale de $p$ et
$C^{*}(X^{top},\mathbb{C})$ son alg\`ebre diff\'erentielle gradu\'ee de
cohomologie. Ainsi, si l'on choisit une telle \'equivalence et si l'on pose
$A=C^{*}(X^{top},\mathbb{C})$, on obtient un diagrame
dans $B-CDGA_{+}$
$$\xymatrix{
A\otimes_{\mathbb{C}}B \ar[r] & \mathbb{R}p_{*}\Omega_{\mathcal{X}/B}^{\leq n-1} \ar[r] &
\mathbb{R}p_{*}\Omega_{\mathcal{X}/B}^{\leq n-2} & \dots \ar[r] &
\mathbb{R}p_{*}\mathcal{O}_{\mathcal{X}}.}$$
Enfin, quitte \`a prendre un remplacement cofibrant de ce diagrame, il est bien connu
que ceci d\'efinit un point dans $B\underline{Filt}^{n}_{A}(B)$. Ce point peut aussi se repr\'esenter
par un morphisme de champs
$$Spec\, B \longrightarrow B\underline{Filt}^{n}_{A},$$
qui est appel\'e \textit{application des p\'eriodes non ab\'eliennes de la
famille $\mathcal{X} \longrightarrow Spec\, B$}. Ce morphisme est bien d\'efini \`a un
isomorphisme non canonique pr\`es (qui d\'epend du choix de la trivialisation
$\mathbb{R}p{*}\Omega^{*}_{\mathcal{X}/B} \simeq
C^{*}(X^{top},\mathbb{C})\otimes_{\mathbb{C}}B$). \\

Pour terminer, signalons que la construction pr\'ec\'edente poss\`ede
aussi la version globale suivante.

Soit $p : X \longrightarrow S$ un morphisme projectif et lisse de
$\mathbb{C}$-sch\'emas s\'epar\'es de type fini. On suppose que $S$ est connexe et on choisit
un point ferm\'e $s \in S$. On pose $A=C^{*}(X_{s}^{top},\mathbb{C})$, la $\mathbb{C}$-alg\`ebre
diff\'erentielle gradu\'ee de cohomologie de la fibre de $p$ en $s$.

On dipose sur l'espace $S^{top}$ d'un faisceau de $\mathbb{C}-CDGA_{+}$,
$\mathbb{R}p^{top}_{*}(\mathbb{C})$ ($p^{top}$ est le morphisme
d'espaces topologique $X^{top} \longrightarrow S^{top}$), qui est localement \'equivalent au faisceau constant
de fibre $A$. Il est donc classifi\'e par un morphisme de monoides simpliciaux (morphisme
de monodromie)
$$\rho : G \longrightarrow \mathbb{R}\underline{Aut}_{\mathbb{C}-CDGA_{+}}(A),$$
o\`u $G$ est le groupe simplicial des lacets de $S^{top}$ (i.e.
un groupe simplicial tel que $BG$ soit \'equivalent au type d'homotopie de $S^{top}$).
En d'autres termes le groupe simplicial $G$ op\`ere sur $A$, et donc aussi sur le pr\'efaisceau
simplicial $B\underline{Filt}^{n}_{A}$. On peut donc former le champ
quotient (i.e. le quotient homotopique par $G$)
$D:=[B\underline{Filt}^{n}_{A}/G]$. Le champ $D$, qui n'est plus $\infty$-g\'eom\'etrique
(mais seulement un \emph{champ analytique $\infty$-g\'eom\'etrique}, en un sens
\`a d\'efinir),
peut \^etre appel\'e le \textit{domaine des p\'eriodes
non-ab\'eliennes de la famille $p$}.

On peut alors construire un morphisme de champs analytiques
$$P : S^{an} \longrightarrow D^{an},$$
qui est une globalisation de l'application des p\'eriodes d\'efinies
ci-dessus. Pour cela il nous faudrait parler du proc\'ed\'e \emph{d'analytification} de champs, et
construire $P$ en recollant des applications des p\'eriodes locales pour les quelles
le faisceau $\mathbb{R}p^{top}_{*}(\mathbb{C})$ est \'equivalent au faisceau constant
de fibre $A$. Nous ne le ferons pas dans cet article. \\

L'espoir est bien entendu que les applications des p\'eriodes non-ab\'eliennes detectent des
variations que l'on ne peut d\'etecter \`a l'aide des variations de structures
de Hodge sur la cohomologie ou m\^eme sur les groupes d'homotopie. On pose donc
la question suivante. \\

\textbf{Probl\`eme:} \emph{Trouver un morphisme projectif et lisse $X \longrightarrow S$ o\`u
$S$ est simplement connexe et o\`u l'application des p\'eriodes non-ab\'eliennes n'est pas constante}. \\

Dans \cite{si5}, C. Simpson construit un exemple de famille $X \longrightarrow S$, avec $S$ simplement
connexe, et pour laquelle la variation sur le type d'homotopie de Dolbeault de la fibre est infinit\'esimalement
non triviale. Cet exemple fournit-il une r\'eponse positive \`a la question pr\'ec\'edente ?

\newpage

\begin{appendix}

\section{Le probl\`eme de la sch\'ematisation}

Dans son manuscript \cite{gr}, A. Grothendieck aborde un probl\`eme qu'il appelle
\textit{probl\`eme de la sch\'ematisation des types d'homotopie}. Dans cet appendice
je r\'esume ce que j'ai personellement compris de cette question. Il s'agit bien entendu
d'un point de vue personel et de ce fait subjectif. \\

Tout d'abord, pour tout sch\'ema de base $S$ (on supposera que $S=Spec\, k$ pour un anneau
$k$), il doit exister un notion \emph{d'$\infty$-champs en groupoides sur $S$} (appel\'es par
la suite simplement \emph{champs sur $S$}, ou \emph{champs sur $k$}). On demande que les champs
sur $S$ soient des g\'en\'eralisations des faisceaux et des champs en groupoides (au sens
de \cite{lamo}), par exemple sur le grand site fpqc de $S$.

Si l'on suit le point de vue
de A. Grothendieck, les champs sur $S$ doivent former une \emph{$\infty$-cat\'egorie}. De fa\c{c}on plus
approximative on peut demander qu'il existe des notions de morphismes de champs et
\emph{d'\'equivalences de champs}, qui permettent alors de parler de la
cat\'egorie homotopique des champs sur $S$, not\'ee $Ho(S)$ (obtenue en inversant formellement
les \'equivalences).

Enfin, on s'attend \`a ce que les constructions standards
sur les faisceaux (limites, colimites, $Hom$ internes, images inverses et directes \dots)
s'\'etendent en des constructions au niveaux des champs. Par exemple il doit exister
une notion de \emph{limites et colimites homotopique de champs sur $S$}. Une mani\`ere
efficace d'obtenir l'existence de telles constructions est
de demander de plus que les champs sur $S$ forment une cat\'egorie de mod\`eles
dont $Ho(S)$ en est la cat\'egorie homotopique. Cependant, pour des raisons qui semblent personelles A. Grothendieck n'envisage pas dans \cite{gr} d'utiliser
la th\'eorie des cat\'egories de mod\`eles.

La th\'eorie des pr\'efaisceaux simpliciaux de \cite{jo,j} donne \'evidemment
un exemple d'une telle th\'eorie des champs (voir \cite{hol} pour
plus de d\'etails sur la comparaison entre champs
en groupoides et pr\'efaisceaux simpliciaux). \\

Supponsons maintenant que l'on dispose d'une th\'eorie des champs sur $S$
satisfaisant les attentes d\'ecrites ci-dessus. En utilisant l'existences de colimites homotopiques on
peut d\'efinir une construction $A \mapsto BA$, qui fait d'un \emph{champ en groupes ab\'eliens $A$},
un nouveau champ en groupes ab\'eliens $BA$ dont le champ des lacets est \'equivalents \`a $A$.
Par it\'erations on obtient des champs d'Eilenberg-MacLane
$$K(\mathbb{G}_{a},n):=B(K(\mathbb{G}_{a},n-1)) \qquad K(\mathbb{G}_{a},0):=\mathbb{G}_{a},$$
o\`u $\mathbb{G}_{a}$ est le faisceau en groupes additif sur $S$.

Les champs $K(\mathbb{G}_{a},n)$ semblent d'un importance capitale
aux yeux de Groth\-en\-dieck, car il les consid\`erent comme les
exemples primitifs fondamentaux de \emph{types d'ho\-mo\-to\-pie
sch\'ematiques}, dont la d\'efinition est un des probl\`eme
central de \cite{gr}.  Bien qu'il ne pr\'ecise jamais ce que
doivent \^etre ces types d'homotopie sch\'ematiques, les
champs $K(\mathbb{G}_{a},n)$ en sont toujours des exemples. On
peut aussi d\'eduire des constructions propos\'ees dans \cite{gr}
que les types d'homotopie sch\'e\-ma\-ti\-ques forment une
sous-cat\'egorie pleine de la cat\'egorie des champs $Ho(S)$ qui
est stable par limites homotopiques (ou encore de fa\c{c}on plus
restrictive on peut demander la stabilit\'e par produits fibr\'es homotopiques uniquement, si l'on veut se
restreindre aux objets de \emph{type fini}).

Nous dirons donc
qu'une sous-cat\'egorie pleine
$\mathcal{C}$ de $Ho(S)$ est une \emph{cat\'egorie de types d'homotopie sch\'ematiques sur $S$},
si elle est stable par limites homotopiques et si elle contient les champs
$K(\mathbb{G}_{a},n)$. \\

On se donne maintenant une cat\'egorie de types d'homotopie sch\'ematiques
sur $S$, $\mathcal{C} \subset Ho(S)$, au sens ci-dessus. L'existence d'une th\'eorie raisonable des
colimites homotopiques dans $Ho(S)$ implique que la cat\'egorie $Ho(S)$ est enrichie
dans la cat\'egorie homotopique  des ensembles simpliciaux $Ho(SSet)$. Ainsi, on peut d\'efinir
un foncteur des sections globales
$$\begin{array}{cccc}
\mathbb{R}\Gamma : & Ho(S) & \longrightarrow & Ho(SSet) \\
 & F & \mapsto & \underline{Hom}(*,F).
\end{array}$$
Par restriction on obtient un foncteur
$$\mathbb{R}\Gamma : \mathcal{C} \longrightarrow Ho(SSet).$$
On dira alors que la cat\'egorie $\mathcal{C}$ \emph{v\'erifie les conditions
du probl\`eme de la sch\'e\-ma\-ti\-sa\-tion} si le foncteur
$\mathbb{R}\Gamma$ poss\`ede un adjoint \`a gauche. Nous noterons cet adjoint
par
$$-\otimes k : Ho(SSet) \longrightarrow \mathcal{C}.$$

Par adjonction on d\'eduit la formule de conservation de la cohomologie suivante,
ch\`ere \`a A. Grothendieck
$$H^{n}(X,k)\simeq H^{n}(X\otimes k,\mathbb{G}_{a}):=[X\otimes k,K(\mathbb{G}_{a},n)]_{Ho(S)},$$
pour tout ensemble simplicial $X$. \\

Soit $\mathcal{C} \subset Ho(S)$ une cat\'egorie de types
d'homotopie sch\'ematiques v\'erifiant les conditions du
probl\`eme de la sch\'ematisation. On suppose de plus que
$S=Spec\, k$ o\`u $k$ est un corps. On rappelle que l'on note
$$-\otimes k : Ho(SSet) \longrightarrow Ho(S)$$
l'adjoint \`a gauche du foncteur $\mathbb{R}\Gamma$.
On peut alors \'enoncer le probl\`eme de la sch\'ematisation de la fa\c{c}on suivante. \\

\textsf{Probl\`eme de la sch\'ematisation:} \emph{Soit $k$ un
corps de caract\'eristique nulle (resp. de caract\'eristique
positive $p>0$). Trouver des cat\'egories de types d'homotopie
sch\'ematiques $\mathcal{C} \subset Ho(Spec\, k)$ qui v\'erifient les conditions
du probl\`eme de la sch\'e\-ma\-ti\-sa\-tion au sens ci-dessus, et tel que le foncteur}
$$-\otimes k : Ho(SSet) \longrightarrow \mathcal{C} \subset Ho(Spec\, k)$$
\emph{restreint aux types d'homotopie rationnels 
$1$-connexes et de type fini
(resp. aux types d'homotopie $p$-complets $1$-connexes et de type fini) soit pleinement fid\`ele.} \\

Bien entendu, de part son \'enonc\'e le probl\`eme de la sch\'ematisation poss\`ede
de nombreuses solutions, dont au moins une triviale lorsque
$\mathcal{C}=Ho(S)$ (et $k$ est alg\`ebriquement clos). Nous avons montr\'e dans
la paragraphe \S 2.3 et \S 2.5 que les champs affines sur un corps $k$ forment
une cat\'egorie $\mathcal{C}$ qui donne un solution au probl\`eme pr\'ec\'edent. \\

Pour finir, signalons que l'on peut \'enoncer des versions quelque peu modifi\'ees du probl\`eme
pr\'ec\'edent, en consid\'erant par exemples des champs point\'es, ou encore
des champs connexes et/ou point\'es \dots Nous laissons le soin au lecteur
de modifier les d\'efinitions pr\'ec\'edentes pour qu'elles soient adapt\'ees \`a ces nouveaux contextes.
Nous avons montr\'e dans le paragraphe \S 3.3 que les types d'homotopie sch\'ematiques
donnaient une solution au probl\`eme de la sch\'ematisation dans le cas
point\'e et connexe.

\end{appendix}

\end{document}